\documentclass[preprint,12pt]{elsarticle}

\usepackage[colorlinks, linktocpage=true, pdfpagemode=No,
linkcolor=blue, citecolor=blue, pagecolor=blue]{hyperref}

\usepackage{listings}
\usepackage{multirow}
\usepackage{csml}
\usepackage{float}
\usepackage{subfig}
\usepackage{cleveref}
\usepackage{array}
\usepackage{multirow}
\usepackage{graphicx}
\usepackage{siunitx}

\begin{document}
\title{Design-Informed Generative Modelling using Structural Optimization} 

\author[1]{Sivanantha Sarma Lowhikan}
\author[1]{Chinthaka Mallikarachchi}  
\author[1]{Sumudu Herath \corref{cor1}} \ead{sumuduh@uom.lk}

\cortext[cor1]{Corresponding author}
\address[1]{Department of Civil Engineering, University of Moratuwa, Bandaranayaka Mawatha, Moratuwa 10400, Sri Lanka}
			   
\begin{abstract}		
Although various structural optimization techniques have a sound mathematical basis, the practical constructability of optimal designs poses a great challenge in the manufacturing stage. Currently, there is only a limited number of unified frameworks which output ready-to-manufacture parametric Computer-Aided Designs (CAD) of the optimal designs. From a generative design perspective, it is essential to have a single platform that outputs a structurally optimized CAD model because CAD models are an integral part of most industrial product development and manufacturing stages. This study focuses on developing a novel unified workflow handling topology, layout and size optimization in a single parametric platform, which subsequently outputs a ready-to-manufacture CAD model. All such outputs are checked and validated for structural requirements; strength, stiffness and stability in accordance with standard codes of practice. In the proposed method, first, topology-optimal model is generated and converted to a one-pixel-wide chain model using skeletonization. Secondly, a spatial frame is extracted from the skeleton for its member size and layout optimization. Finally, the CAD model is generated using constructive solid geometry trees and the structural integrity of each member is assessed to ensure structural robustness prior to manufacturing. Various examples presented in the paper showcase the validity of the proposed method across various engineering disciplines.
\end{abstract}


\maketitle
\section{Introduction}
\label{sec:introduction}

\subsection{Background}\label{subsec:background}
Structural optimization is a set of well-established mathematical techniques for finding the optimal design of a structure with minimal stress, weight or compliance for a given amount of materials and constraints. It systematically improves geometric designs by increasing stiffness, lowering material consumption, and reducing production time \cite{dede2019usage, kurdi2015structural, mei2021structural}. Evidently, these optimal designs will lead to low-carbon production and processes with significantly improved material efficiency. Applications of structural optimization methods are diverse: ranging from small robot arms \cite{YIN2020113102, yunfei2016structural} to multi-span bridges \cite{li2022innovative} spanning across numerous engineering disciplines such as mechanical~\cite{bandyopadhyay2019additive, airbus.com}, structural~\cite{kingman2014applications, tsavdaridis2015applications}, bio-medical~\cite{deng2016topology} and aerospace~\cite{zhu2016topology,liu2019topological} covering many industries like military~\cite{hayduketopology}, robotics \cite{yunfei2016structural, wildman2019topology}, construction~\cite{mei2021structural, stoiber2021topology, vantyghem20203d}, automobile~\cite{airbus.com, matsimbitopology, zhang2021integrated, yang1995automotive}, space vehicles~\cite{autodesk_2021, kang2009topology, ashikhmina2019wing} and medical~\cite{liao2007shape, wu2021advances}. In this paper, structural optimization in civil and structural engineering applications is outlined. However, it is noted that the extension of the proposed workflow into other applications is direct and straightforward.

Although structural optimization techniques have produced material efficient designs, the practical manufacturing and implementation of those designs are limited due to the complex geometries of the optimized design outputs~\cite{YIN2020113102}. The current advancements in additive manufacturing partially sidestep this problem predominantly in the field of mechanical engineering \cite{wu2016critical,buchanan2019metal,vafadar2021advances}. However, civil engineering applications are sparse, and there is potential for these techniques to be used, especially in precast and modular constructions. Additive manufacturing in the construction industry is gaining popularity, however to-date, the majority of the structural designs are performed as per the guidelines given in standard codes of practice such as; BS8110-Structural use of Concrete~\cite{bs19978110}, Euro Code 2-Design of Concrete Structures~\cite{standard2004eurocode} for concrete and BS5950-Structural use of Steelwork in Building~\cite{standard2001bs}, Euro Code 3-Design of Steel Structures~\cite{standard2006eurocode} for steel. This research work investigates the possibility of merging structurally optimized geometric designs with structural design requirements to arrive at robust structural designs that are ready-to-assemble or -construct. 

To the best knowledge of the authors, there is only a limited number of unified frameworks which output ready-to-additively manufacture parametric Computer-Aided Design (CAD) models of optimal designs. Commercially available software such as Abaqus FEA \cite{abaqus2009}, Ansys and LS-DYNA \cite{lsdyna2003} do not facilitate the output of a structurally optimized CAD model. Instead, those software identifies the critical load path to consider and regions for material removal. Hence, to output the CAD models, which are the integral parts of most industrial product development and manufacturing \cite{YIN2020113102, smith2016application}, manual intervention is essential. To circumvent these issues, there are several latest developments such as Limitstate:Peregrine \cite{limitstate} and Ameba \cite{YiMinXiesteams} where the structurally efficient designs are produced in the Rhino Grasshopper platform. However, they do not evaluate the structural integrity of the optimal designs against structural design standards. The proposed work in this paper attempts to mitigate all the above-identified weaknesses and non-implementations in existing works.
\begin{figure}[H]
\centering
\includegraphics[width=0.99\textwidth]{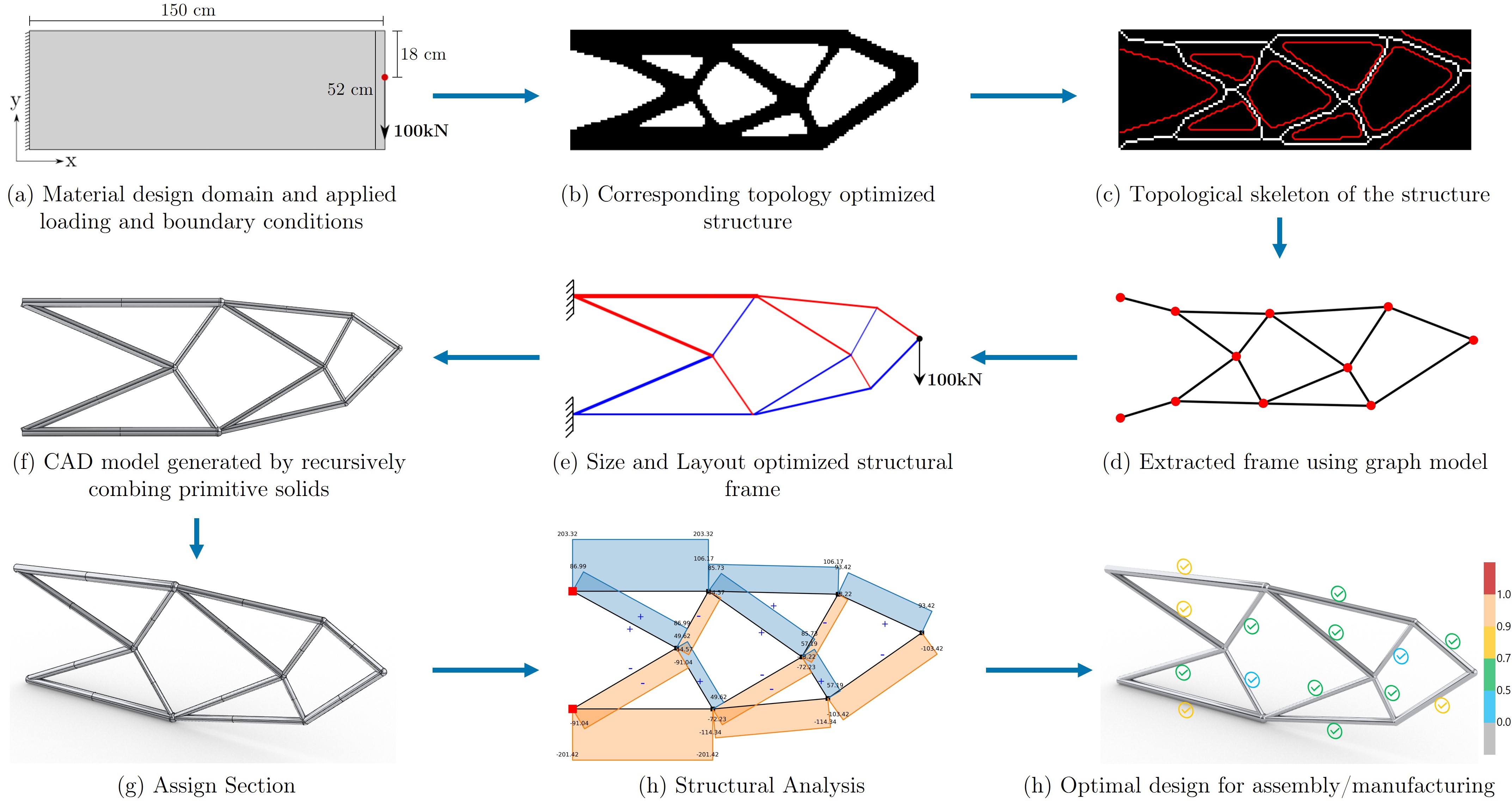}
\caption{Structurally robust CAD model generation workflow from topology, size and layout optimal structure for a cantilevered beam example}
\label{fig:overall_illustration}
\end{figure}
This paper presents a novel approach to producing structurally sound, design guideline complaint, ready-to-manufacture CAD models of structurally optimized designs. In the proposed method, first, any two-dimensional (2D) surface is discretized into some finite elements (equivalent to finite-size pixels in 2D) and a topology-optimized 2D binary image is generated. Then, the obtained pixel model is converted into a one-pixel-wide chain model using a skeletonization algorithm. From the obtained skeleton, a spatial frame structure is extracted and represented by a graph model. Thereafter, the sizes and layouts of the frame members are iteratively optimized. Finally, the CAD model is generated using constructive solid geometry trees and its structural performance is assessed. In addition, the industry-standard structural sections are assigned to the optimal structure and designed under standard codes of practice. Herein, EN 1993-1-1: Eurocode 3: Design of steel structures~\cite{standard2006eurocode} is used as the code of practice. Figure \ref{fig:overall_illustration} illustrates this procedure briefly using a cantilever example. The proposed workflow is entirely implemented in Rhinoceros CAD software utilizing the visual programming tool Grasshopper in a Python coding environment.

This section briefly highlighted the background of this research. Section \ref{subsec:related_work} presents a few of the related works and Section \ref{subsec:review_stropt} discusses the theoretical aspects of the optimization techniques. In Section \ref{sec:Opt_CAD_Gen}, the proposed methodology adopted to obtain the optimized models is detailed. The new implementations are explained through the aid of a simple cantilever beam. Section \ref{sec:design_informed_stropt} elaborates the structural analysis and design implementation adopted for the optimized models. The applications of using the developed workflow are presented in Section \ref{sec:applications}. Finally, \ref{sec:conclusion} concludes the paper by outlining the key contributions of this study to the related research and industrial communities.

\subsection{Related Work}\label{subsec:related_work}
Several recent studies have focused on finding the optimal designs of trusses and frames within various engineering and design domains. This section presents a review of some of the most relevant and recent works found in the literature.

Structural optimization of trusses has been a widely researched area in recent years. One of the most recent and interesting works is the adaptive layout optimization of trusses by the University of Sheffield \cite{he2019python}. The study presented a simple Python implementation of a truss layout optimization using an adaptive ‘member adding’ scheme. While previous research studies focused on ground structure-based layout optimization, the results were deemed impractical from an assembly/manufacturing perspective. However, a series of research works \cite{gilbert2003layout, he2015rationalization, he2019conceptual, liu2014efficient} provided more efficient ways to produce optimum truss structures. Further, CAD model generation for additive manufacturing techniques has also been presented, and a grasshopper optimization plugin named “Peregrine” \cite{limitstate} was released as a result of the deliverables. 

Another inspiring work in this area is the topologically robust CAD model generation workflow presented by Yin et al. \cite{YIN2020113102}. The study presented a fully automated workflow for the conversion of a topology-optimized model into a CAD model using staged skeletonization and frame extraction, which has been a major motivation for this study (see Fig. \ref{fig:overall_illustration}. However, it should be noted that this study was focused more on framework development rather than design applications in Civil engineering design and applications. Moreover, in this paper, the limitations highlighted by Yin et al. \cite{YIN2020113102} are also addressed in relation to framed truss structures. An alternative approach for the sequential size and layout optimization of structures is the nodal-based evolutionary optimization of frame structures presented by X. Zhang et al. \cite{zhang2022nodal}. The study proposed an optimization algorithm to minimize the volume and compliance of frame structures, where the size and nodal coordinates are optimized in parallel. 
\subsection{Review of Structural Optimization}\label{subsec:review_stropt}
This section summarizes each structural optimization technique used in this study, especially its formulation and theoretical background. The standard density-based topology optimization using SIMP (Solid Isotropic Material with Penalization) method is used in this study. For size and layout optimization, bi-dimensional frame elements which can transfer forces and moments are considered. The stiffness matrix of the frame element is obtained by combining the stiffness of the truss/bar element and the Bernoulli beam element \cite{ferreira2009matlab}. In this study, sequential size and layout optimizations are performed to obtain both size and layout optimized structure. Parallelization of the two is also a possibility, as outlined in Zhang at el. \cite{zhang2022nodal}. However, this study uses a staged-sequential approach where size and layout optimizations are iterated until a global convergence is attained. The following sections are aimed only to outline the formulations of structural optimization methods. For interested readers, detailed information about these optimization methods can be found in \cite{bendsoe2003topology, christensen2008introduction, ohsaki2016optimization,rozvany2014shape}. 
\subsubsection{Topology Optimization}\label{subsec:TopOpt_Review}
Topology optimization is focused on finding the optimal material distribution within a specified design domain \cite{bendsoe1988generating,bendsoe2003topology,sigmund200199}. The topology optimization problem of a finite element discretized solid is given by \cite{YIN2020113102,bendsoe2003topology,sigmund200199,herath2021topologically},
\begin{equation}\label{eq:TopOpt1}
    \mathrm{minimize}\ C_{(\bm \rho )}= \bm f^{T} \bm u_{(\bm \rho )}
\end{equation}
\begin{equation}\label{eq:TopOpt2}
    \mathrm{subject\ to}\  \bm K_{(\bm \rho )}\bm u = \bm f
\end{equation}
\begin{equation}\label{eq:TopOpt3}
    \mathrm{design\ constraint:}\ \frac{V_{(\bm \rho )}}{\overline{V}}\le V_f
\end{equation}
\begin{equation}\label{eq:TopOpt4}
    \mathrm{design\ variable:}\ 0\le \bm \rho \le 1
\end{equation}
where $C_{(\bm \rho )}$ is the objective function, $\bm \rho $ is the vector of relative element densities, $\bm u$ is the displacement vector, $\bm K_{(\bm \rho )}$ is the global stiffness matrix, $\bm f$ is the global external force vector, $V_{(\bm \rho )}$is the material volume, $\overline{V}$ is the design volume and $V_f$ is the prescribed volume fraction. The relative density of each element is constrained to be $\ 0\le {\rho }_i\le 1$ with $i\ \epsilon \ \left\{1,2,3\dots n_e\right\}$ where $n_e$ is the number of elements in the domain. The global stiffness matrix $\bm K$ and vector $\bm f$ are assembled from the $n_e$ element contributions $\bm K_i$ and $\bm f_i$ respectively. In this paper, we discretize the design domain with a structured grid and hexahedral linear elements. In each element, the material is isotropic and homogeneous and Young's modulus $E$ is penalized depending on the relative density ${\bm \rho }_i$ according to, 
\begin{equation}
    E_{{(\rho }_i)}=E_{min}+{\rho }^p_i(\overline{E}-E_{min})
\end{equation}
where $\overline{E}$ is the prescribed Young's modulus of the solid material, $\rho $ and $E_{min}$ are two algorithmic parameters. The penalization parameter $\rho \le 3$ ensures that elements with densities close to ${\rho }_i=0$ (void) and ${\rho }_i=1$ (solid) are favoured. An acceptable penalization power value such as 3 or higher is suggested for topology optimization problems according to Bends{\o}e \& Sigmund, 2003 \cite{bendsoe2003topology}. The small Young's modulus $E_{min}\cong {10}^{-9}$ of the void material prevents ill-conditioning of the global stiffness matrix when ${\rho }_i=0$. Each element stiffness matrix $\bm K_i$ is computed using the corresponding Young's modulus $E_{{(\rho }_i)}$ with the relative density ${\rho }_i$ which is constant within an element.

Further, in order to prevent checkerboarding and mesh dependency instabilities, the topology optimization problem \cref{eq:TopOpt1,eq:TopOpt2,eq:TopOpt3,eq:TopOpt4} is regularized by filtering out the element densities \cite{bruns2001topology}. Filtered densities can be found by using the kernel function. A basic filter density function is given as,
\begin{equation}
   \widehat{\rho }_i=\frac{\sum_j{H_{(i,j)}v_j{\rho }_i}}{\sum_j{H_{(i,j)}v_j}} 
\end{equation}
where ${\widehat{\rho }}_i$ is filtered density,  $v_j$ is the volume of element $j$, $H\left(i,j\right)$ is weight factor and the summations are over the elements inside the filter radius of element $i.$ The weight factor $H\left(i,j\right)$ can be defined as a function of the distance between neighbouring elements, 
\begin{equation}
    H\left(i,j\right)=\left\{ \begin{array}{c}
    R-dist\ \left(i,j\right)\ \ \ \ \ \ \ \ \mathrm{\ \ if\ }dist\left(i,j\right)<R \\ 
    0\ \ \ \ \ \ \ \ \ \mathrm{\ \ \ \ \ \ \ \ \ \ \ \ \ \ \ \ otherwise} \end{array}
    \right.
\end{equation}
where the operator $dist\ \left(i,j\right)$ is the distance between the center of element $i$ and the center of element $j$, and $R$ is the prescribed filter radius.

Gradient-based optimization requires the derivatives of the objective function as well as constraints with respect to the design variable. In density-based topology optimization the relative density ${\rho }_i$ is the design variable. The relative densities ${\rho }_i$ in the topology optimization problem \cref{eq:TopOpt1,eq:TopOpt2,eq:TopOpt3,eq:TopOpt4} are substituted with the filtered relative densities in the following ${\widehat{\rho }}_i$ \cite{bendsoe2003topology,bruns2001topology,YIN2020113102}. The derivative of the compliance objective function according to the chain rule is computed as,
\begin{align}
    \frac{\partial C_{(\widehat{\bm \rho })}}{\partial {\rho }_i} & =\frac{\partial C_{(\widehat{\bm \rho })}}{\partial {\widehat{\rho }}_j\mathrm{\ }}\frac{\partial {\widehat{\rho }}_j}{\partial {\rho }_i} \\
    \frac{\partial C_{(\widehat{\bm \rho })}}{\partial {\widehat{\rho }}_j} & = \bm f^T\frac{\partial \bm u_{(\widehat{\bm \rho })}}{\partial {\widehat{\rho }}_j}. \label{11}
\end{align}
Similary, differentiating the equilibrium constraint,
\begin{align}
    \frac{\partial \bm u_{(\widehat{\bm \rho })}}{\partial {\widehat{\rho }}_j} & =  - \bm K^{-1}\frac{\partial \bm K_{(\widehat{\bm \rho })}}{\partial {\widehat{\rho }}_j} \bm u \\
    \frac{\partial C_{(\widehat{\bm \rho })}}{\partial {\widehat{\rho }}_j} & = -\bm u^T\frac{\partial \bm K_{(\widehat{\bm \rho })}}{\partial {\widehat{\rho }}_j} \bm u \nonumber \\
    & =-\sum^{n_{ele}}_{l=1}{{\bm u_l}^T\frac{\partial \bm K_{l(\widehat{\bm \rho })}}{\partial {\widehat{\rho }}_j} \bm u_l} \nonumber \\
    & =-p{\rho }^{p-1}_j\frac{\left(\overline{E}-E_{\mathrm{min}}\right)}{\overline{E}}{\bm u_j}^{T}{\overline{\bm K}_j}{\bm u_j}
\end{align}
The derivative of the filtered densities with respect to relative density is,
\begin{equation}\label{14}
    \frac{\partial {\widehat{\rho }}_j}{\partial {\rho }_i}\ =\frac{H_{(i, j)}v_i}{\sum_k{H_{(j,k)}v_k}}
\end{equation}
where the summation in the denominator is over the elements inside the filter radius of element $j$. By substituting equations \eqref{14} and \eqref{15} in \eqref{11}, the derivative of the objective function with respect to unfiltered relative density can be obtained as,
\begin{equation}\label{15}
    \frac{\partial C_{(\widehat{\bm \rho })}}{\partial {\rho }_i}=\sum^{n_{ele}}_j{\left(-p{\rho }^{p-1}_j\frac{\left(\overline{E}-E_{\mathrm{min}}\right)}{\overline{E}}{\bm u_j}^T{\overline{\bm K}}_j \bm u_j\frac{H_{(i, j)}v_i}{\sum_k{H_{(j,k)}v_k}}\right)}
\end{equation}
Similarly, the derivative of volume constraint with respect to relative density ${\rho }_i$ can be found as below.
\begin{equation}
    V_{\left(\widehat{\bm \rho }\right)}=\ \sum_i{{\widehat{\rho }}_jv_i}\le V_f\overline{V}    
\end{equation}
\begin{equation}
    \frac{\partial V_{\left(\widehat{\bm \rho }\right)}}{\partial {\rho }_i}=\ \frac{\partial V_{\left(\widehat{\bm \rho }\right)}}{\partial {\widehat{\rho }}_j\mathrm{\ }}\frac{\partial {\widehat{\rho }}_j}{\partial {\rho }_i}=v_i\sum^{n_{ele}}_j{\left(\frac{H_{(i, j)}v_i}{\sum_k{H_{(j,k)}v_k}}\right)}  
\end{equation}
\subsubsection{Size Optimization}\label{subsec: SizeOpt_Review}
In size optimization, the optimization variable is a vector of cross-section areas or any other parameters describing the cross-section \cite{kurdi2015structural}. Size optimization follows a similar formulation process as the topology optimization briefed in Section \ref{subsec:TopOpt_Review}. The related equations take the forms,
\begin{equation}
    \mathrm{minimize} ~C_{(\bm A)}= \bm f^T \bm u 
\end{equation}
\begin{equation}
    \mathrm{subject\ to}\ \bm K_{(\bm A)} \bm u = \bm f
\end{equation}
\begin{equation}
    \mathrm{design\ constraint:}\ \frac{V_{(\bm A)}}{\overline{V}}\le V_f
\end{equation}
\begin{equation}
    \mathrm{design\ variable:}\ \bm A_l\le \bm A\le \bm A_u
\end{equation}
where the cross-section areas considered here is denoted as $\bm A = (A_1,A_2,A_3 \dots A_{n_{elm}})$, where a component $A_i$ represents the area of member $\ i$ and $n_{elm}$ denotes the number of members in the frame model. The lower bound and upper bound of the cross-sectional area, which constrains the member sizes to fit the manufacturing requirements is represented as $\bm A_l$ and $\bm A_u$, respectively. The derivatives of the objective function and constraints with respect to cross-sectional area $\bm A$ are expressed as,
\begin{equation}
    \frac{{\partial C}_{\left(\bm A\right)}}{\partial A_i}=-\bm u^T\frac{\partial \bm K_{\left(\bm A\right)}}{\partial A_i} \bm u=-\int^{n_{elm}}_{j=1}{\bm u^T_j\frac{\partial {\bm K_j}_{\left(\bm A\right)}}{\partial A_i} \bm u_j}
\end{equation}
where$\ i,j=(1,2,3,\dots \dots \dots , n_{ele})$ are element indices and $u_j$ is the vector of nodal displacements and rotation of element $j$. It should be noted that,  $\frac{\partial {K_j}_{\left(A\right)}}{\partial A_i}=0\ $ when  $A_i$ is not the area of element $j$.

The volume constraint in both size and layout optimizations can be expressed as,
\begin{equation}\label{eq:vol_cons}
    V_{\left(\bm {A,s}\right)}=\int^{n_{elm}}_i{A_iL_i\le }V_f\overline{V}
\end{equation}
where $L_i$ is the length of member $i$. Hence, the derivative of the volume constraint in size optimization is found as,
\begin{equation}
    \frac{{\partial V}_{(\bm {A,s})}}{\partial A_i}=L_i 
\end{equation}
where $i\mathrm{=(1,2,3\dots \dots \dots }n_{elem}\mathrm{)}$ is the index of elements.
\subsubsection{Layout Optimization}\label{subsec: LayOpt_Review}
Herein, the optimization variables are the frame nodal coordinates. Layout optimization follows the same computational procedure as the topology and size optimizations where the formulation is succinctly written as,
\begin{equation}
    \mathrm{minimize}\ C_{(\bm s)}= \bm f^T \bm u
\end{equation}
\begin{equation}
    \mathrm{subject\ to}\ \bm K_{(\bm s)} \bm u = \bm f
\end{equation}
\begin{equation}
    \mathrm{design\ constraint:}\ \frac{V_{(\bm s)}}{\overline{V}}\le V_f
\end{equation}
\begin{equation}
    \mathrm{design\ variable:}\ \bm s_l\le \bm s\le \bm s_u
\end{equation}
The nodal coordinates considered here are denoted as $\bm s=(n_1,n_2,n_3 \dots \\ n_{n_{nodes}})$, where a component $s_i$ represents the node of an element in a structure and $n_{nodes}$ represent the number of nodes in the structure. In a 2D Cartesian coordinate system, the coordinate of a node $n_i$ is simply expressed as $n_i=(x^{(i)},y^{(i)})$. Furthermore, the lower bound and upper bound of the nodal coordinate, $\bm s_l$ and $\bm s_u$ are used to vary the coordinates inside the design domain. These bounds can be assigned to each node individually or can be allowed to vary within the whole design domain. The derivative of the objective function with respect to nodal coordinate can be expressed as,
\begin{align}
    & \frac{{\partial C}_{\left(\bm s\right)}}{\partial s_i} =- \bm u^T\frac{\partial \bm K_{\left(\bm s\right)}}{\partial s_i}\bm u=-\int^{n_{elm}}_{j=1}{\bm u^T_j\frac{\partial {\bm K_j}_{\left(\bm s\right)}}{\partial s_i}\bm u_j} \\
    & =-\int^{n_{elm}}_{j=1}{\bm u^T_j\left(\frac{{\partial \bm T}^T_j}{\partial s_i}{\bm K_j}^l \bm T_j\ +\ \bm T^T_j\frac{\partial {\bm K_j}^l}{\partial s_i}\bm T_j+ \bm T^T_j{\bm K_j}^l\frac{\partial \bm T_j}{\partial s_i}\ \right) \bm u_j}
\end{align}
where$\mathrm{\ }i\mathrm{=(1,2,3,\dots }{,n}_{nodes}\mathrm{)}$ and $j\mathrm{=(1,2,3\dots }n_{elem}\mathrm{)}$. Here $T_j$ denotes the transformation matrix of element $j$ and  ${K_j}^l$ indicates the local stiffness matrix of element $j$. It should be noted that,  $\frac{\partial \bm T_j}{\partial s_i}\mathrm{=0\ }$ and $\frac{\partial {{\bm K_j}^l}_{\left(s\right)}}{\partial s_i}\mathrm{=0\ }$ when $s_i$ is not the node on element $j$. 

Also, the derivative of the volume constraint \eqref{eq:vol_cons} function with respect to nodal coordinates takes the form,
\begin{equation}
\frac{{\partial V}_{(\bm A,s)}}{\partial s_i}=\int^{n_{elm}}_{j=1}{A_j\frac{\partial L_{j(\bm s)}}{\partial s_i}}
\end{equation}
where $i =(1,2,3,\dots {n}_{nodes})$ and $j=(1,2,3\dots n_{elem})$. It should be noted that,  $\frac{\partial L_{j(\bm s)}}{\partial s_i}\mathrm{=0\ }$ when $s_i$ is not the node on element $j$.

Chapter 1 provides the necessary background for the research, including a review of related works and an overview of various optimization techniques. This chapter establishes the context for the subsequent chapters, highlighting the need for a methodology to generate optimized models for continuum topologies. Chapter 2 focuses on presenting the structurally optimized CAD model generation workflow for discrete models converted from continuum topology-optimized models. The chapter provides a detailed explanation of the conversion process and the steps involved in generating an optimized model suitable for analysis.


\newpage 

\section{Optimal CAD Model Generation}
\label{sec:Opt_CAD_Gen}
In this Section, the proposed workflow to convert a continuum topology-optimized model into a discrete frame model is presented. Herein, a cantilever plate example is used to illustrate the proposed workflow with implementation details. First, the topology optimization of the selected example is explained in Section \ref{subsec:TopOpt}. Next, Section \ref{subsec:Skeletonization} provides the fundamental concepts of skeletonization for a topology-optimized model. Then, Section \ref{subsec:frame_skel} presents a procedure to generate a frame model with the help of graph theory. The sequential size and layout optimization of a generated frame model is provided in Section \ref{subsec:seq_size&layout_opt}. Finally, the CAD model generation using Constructive Solid Geometry (CSG) is explained in Section \ref{subsec:CAD_modelgen}.
\subsection{Topologically optimal continuum model}\label{subsec:TopOpt}
\vspace{-18pt}
\begin{figure}[H] 
	\centering
	\subfloat[Choosen Cantilever Example]{
		\includegraphics[width=0.47\textwidth]{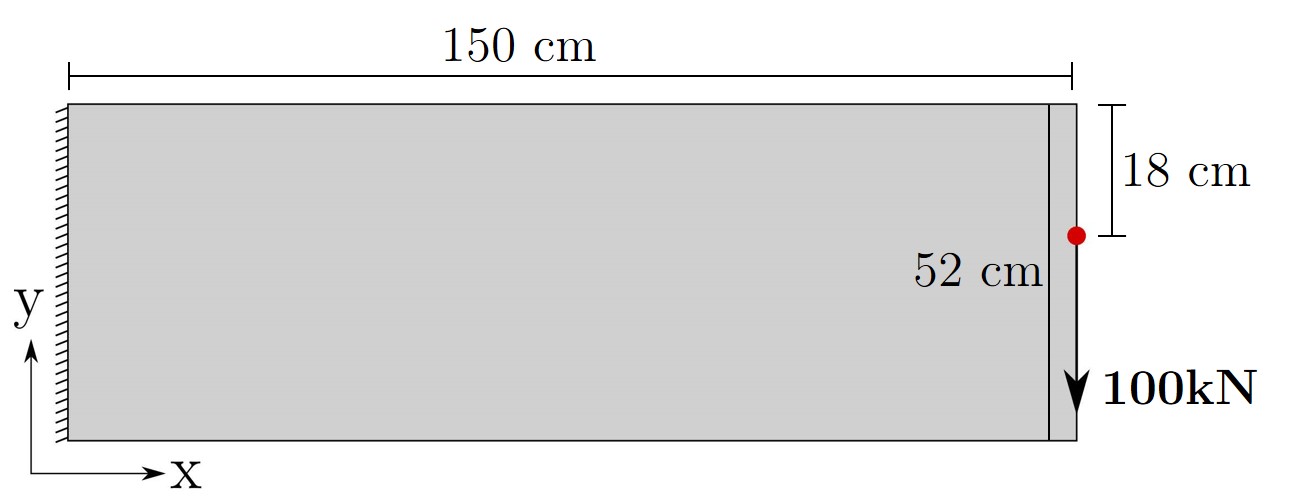}
		\label{fig:TopOpt_CantPlate_DifVol_a}
	}
	\centering
	\subfloat[$V_f$ = $0.7$]{
		\includegraphics[width=0.47\textwidth]{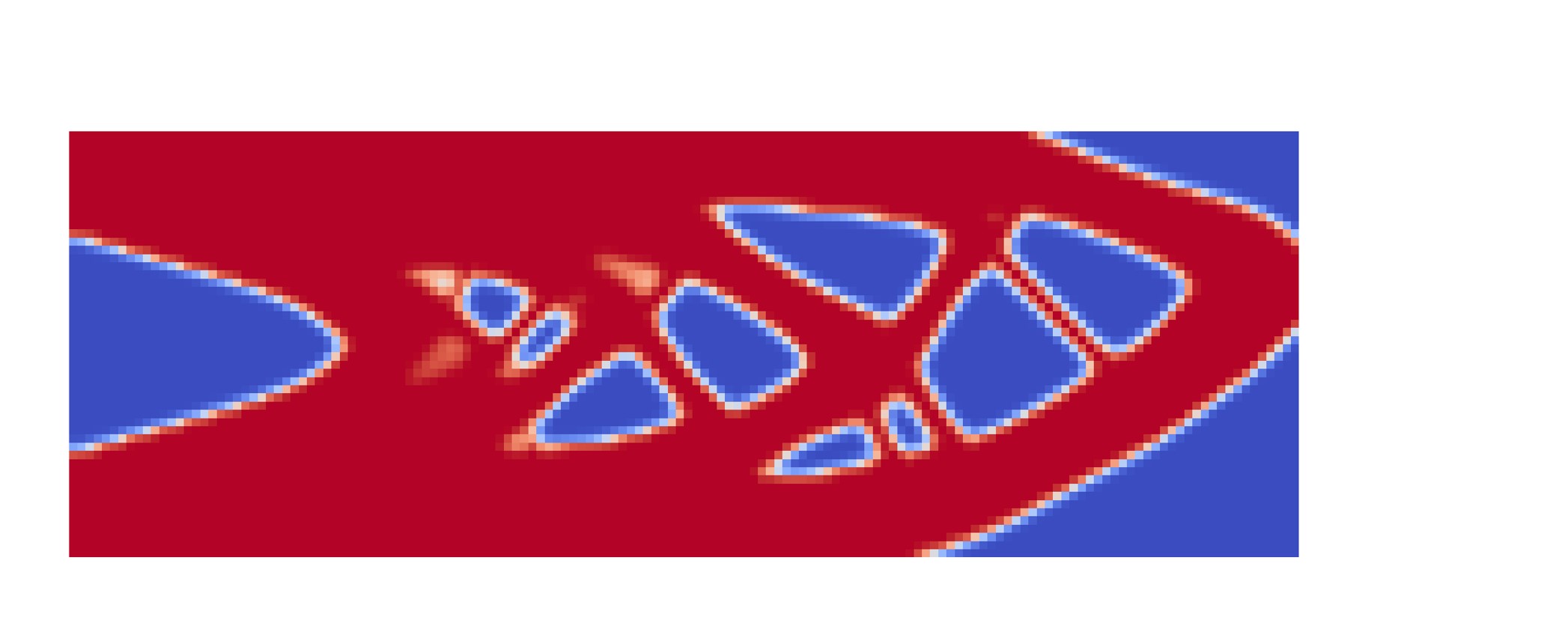}
		\label{fig:TopOpt_CantPlate_DifVol_b}
	} \\[-3pt]
        \centering
	\subfloat[$V_f$ = $0.5$]{
		\includegraphics[width=0.47\textwidth]{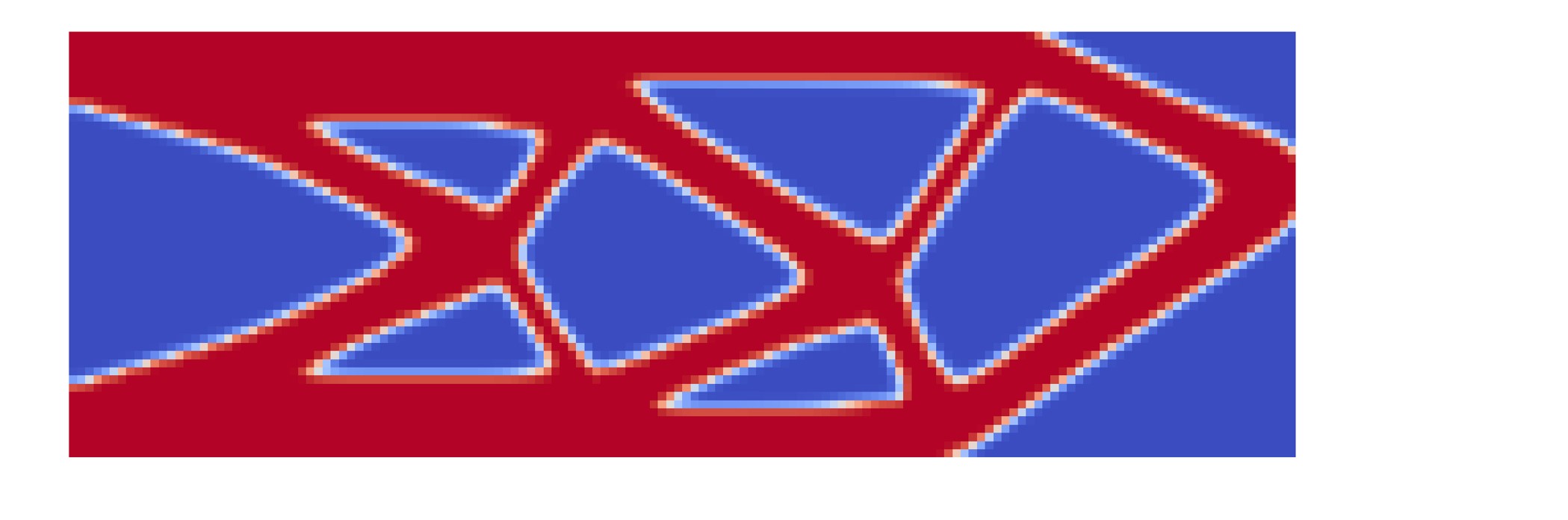}
		\label{fig:TopOpt_CantPlate_DifVol_c}
	}
	\centering
	\subfloat[$V_f$ = $0.3$]{
		\includegraphics[width=0.47\textwidth]{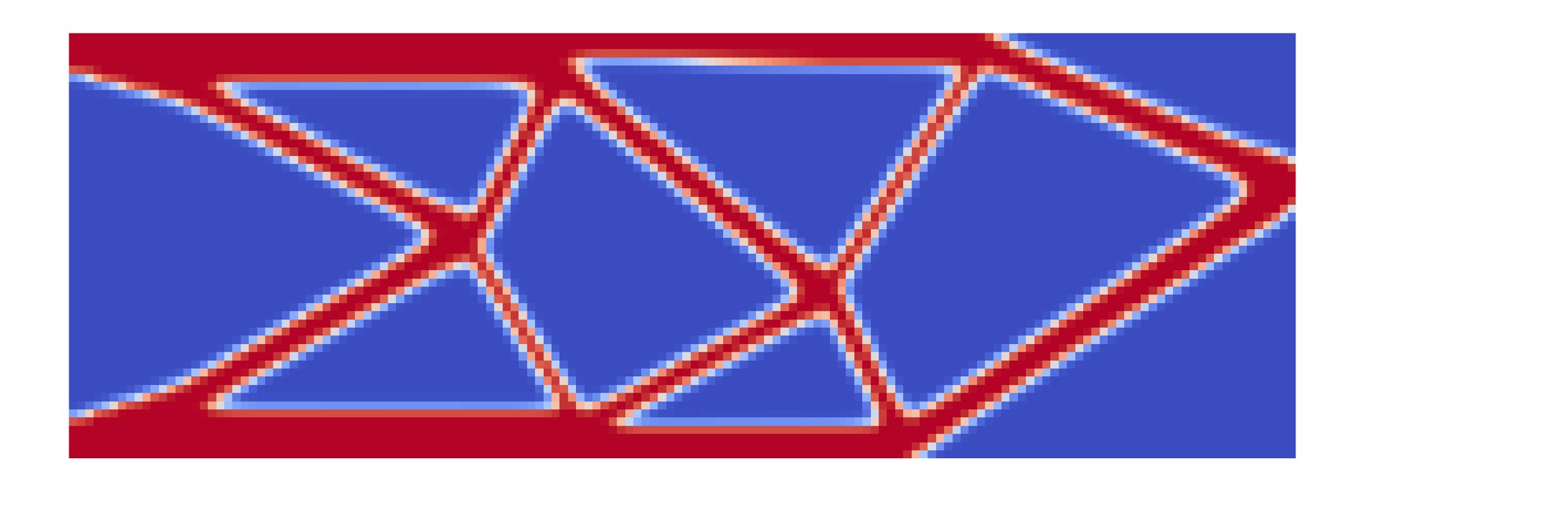}
		\label{fig:TopOpt_CantPlate_DifVol_d}
	} \\[-3pt] 
	\centering
	\caption{Topology optimization of a cantilever plate. Design, material, and optimization parameters are   $\overline{E}$ = $2.1\ \times {10}^5\ {N/mm}^2$, ${E}_{min}$ = ${1\ \times 10}^{-9}\ {N/mm}^2$, $\nu$ = $0.3$, $p$ = $3$, $R$ = $1.2$, Maximum iterations = $200$.}
	\label{fig:TopOpt_CantPlate_DifVol}
\end{figure}

To perform the topology optimization, first, the design problem and optimization parameters should be decided. Herein, a cantilever plate example is selected to illustrate the proposed workflow. As shown in Figure \ref{fig:TopOpt_CantPlate_DifVol}, a 150 cm x 52 cm plate is selected, and one edge of the domain is fixed. A point load with a value of F \(= \) 100 kN is applied 18 cm below the top edge of the plate on the opposite free end. The thickness of the plate is taken as 1 cm for size optimization. The numerical results are discussed later in Section \ref{sec:applications}. At this juncture, the focus is rendered on the parameters that can affect the initial geometry needed for skeletonization. 

As shown in Figure \ref{fig:TopOpt_CantPlate_DifVol_b} to \ref{fig:TopOpt_CantPlate_DifVol_d}, different volume fractions ($V_f$) will lead to distinct topologically optimized outputs. Similarly, different mesh sizes can vary these outputs. Filter radius will also affect the binary conversion of these optimized models. One can conduct comparative analysis to find the optimal volume fraction, discretization, and filter radius for the problem as conducted by \cite{9906146,bendsoe2003topology}. However, since the goal of this study is to provide an approach to extract frame models, generate CAD models and perform structural design validations, less emphasis is placed on achieving very precise topology-optimized solutions. 

For a better input to the skeleton extraction, a topology-optimized output without noise (such as short branches stemming out) is desirable. Therefore, it is recommended to select suitable problem parameters in order to obtain an accurate frame output after skeletonization. However, we show that an initial topology-optimized structure with a similar critical load path (see Figure \ref{fig:TopOpt_CantPlate_DifVol_c} and \ref{fig:TopOpt_CantPlate_DifVol_d}) will produce the same skeleton and optimized frame. 

\subsection{Skeletonization}\label{subsec:Skeletonization}

Skeletonization is the process of finding the topological skeleton of a structure in any spatial dimension. In 2D skeletonization, it reduces a 2D binary object of a continuum to a one-pixel-wide representation with the same connectivity as the original \cite{deng2000fast,zhang1984fast}. In this study, the 2D binary image data obtained from the topology optimization is passed through the skeletonization process to extract the medial axis of the continuum structure. 

\subsubsection{Review of 2D Binary Image Skeletonization Techniques}\label{subsec:review_2DSkel}

2D binary image skeletonization techniques have become increasingly influential in image-processing applications such as character recognition, military route finding, medical imaging, computer graphics and chromosome analysis \cite{deng2000fast, tagliasacchi20163d}. Binary image skeletonization techniques are broadly divided into two main categories: i) distance transform methods \cite{xia1989skeletonization} and ii) parallel thinning methods \cite{deng2000fast, chen1990new}. The latter can be used to capture the structural information of objects, making it the preferred method for civil engineering applications where the preservation of such information is vital. Hence in this study, the parallel thinning method is used for topological skeleton extraction.

\begin{figure}[H] 
	\centering
    \captionsetup{justification=centering}
	\subfloat[Tracing border pixels]{
		\includegraphics[width=0.49\textwidth]{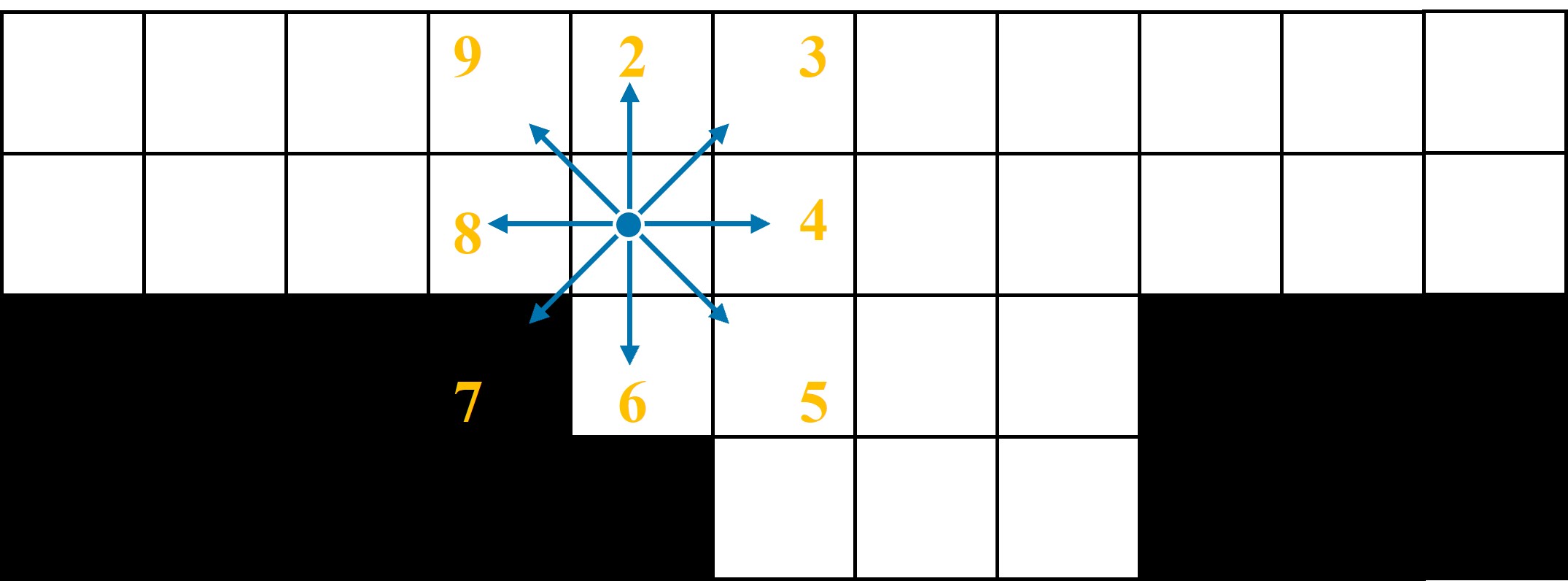}
		\label{fig:pix_cat_a}
	}
	\centering
	\subfloat[Detected pixels for removal]{
		\includegraphics[width=0.49\textwidth]{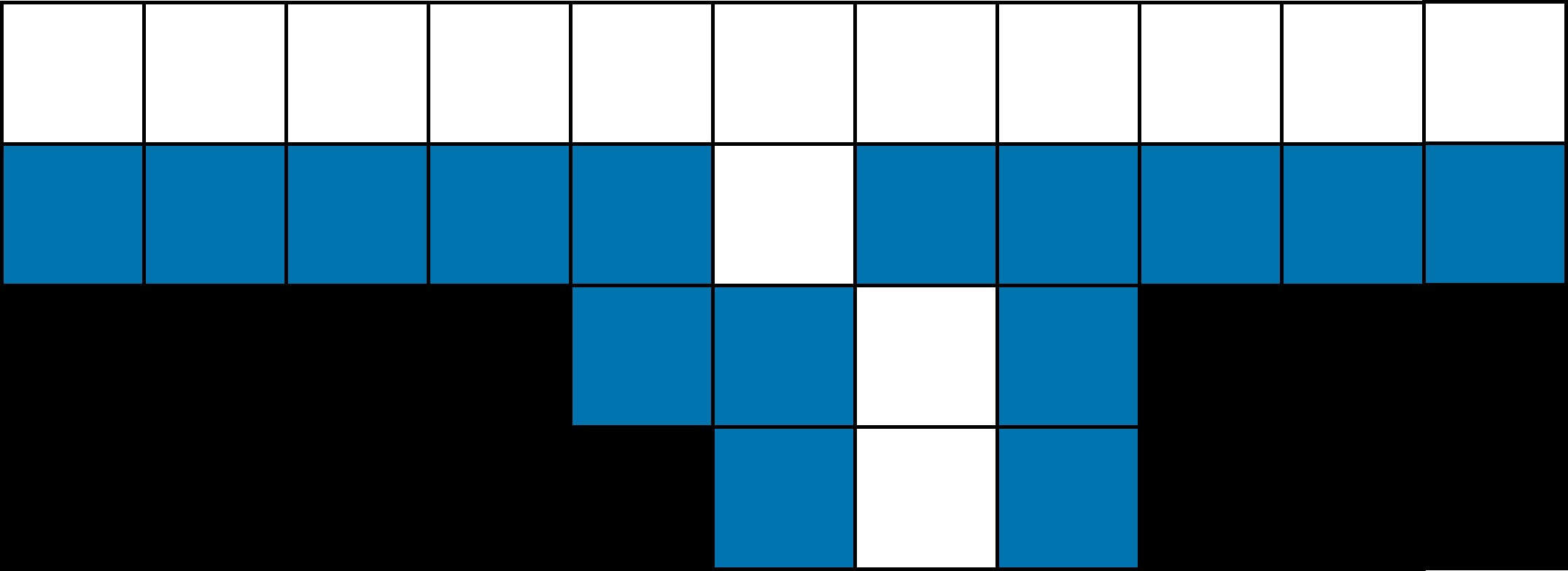}
		\label{fig:pix_cat_b}
	} \\[-3pt]
        \centering
	\subfloat[White pixel chain – Skeleton]{
		\includegraphics[width=0.49\textwidth]{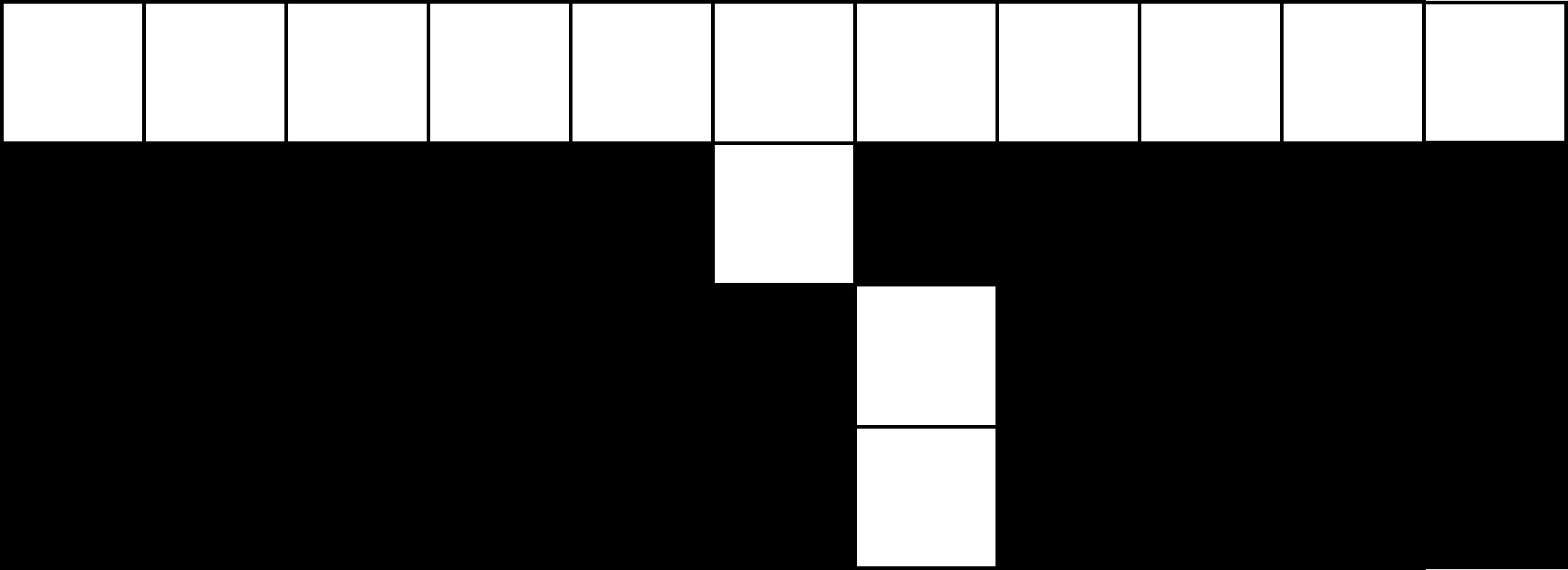}
		\label{fig:pix_cat_c}
	}
	\centering
	\subfloat[Pixel Types (Red–End, Green–Joint, White–Regular)]{
		\includegraphics[width=0.49\textwidth]{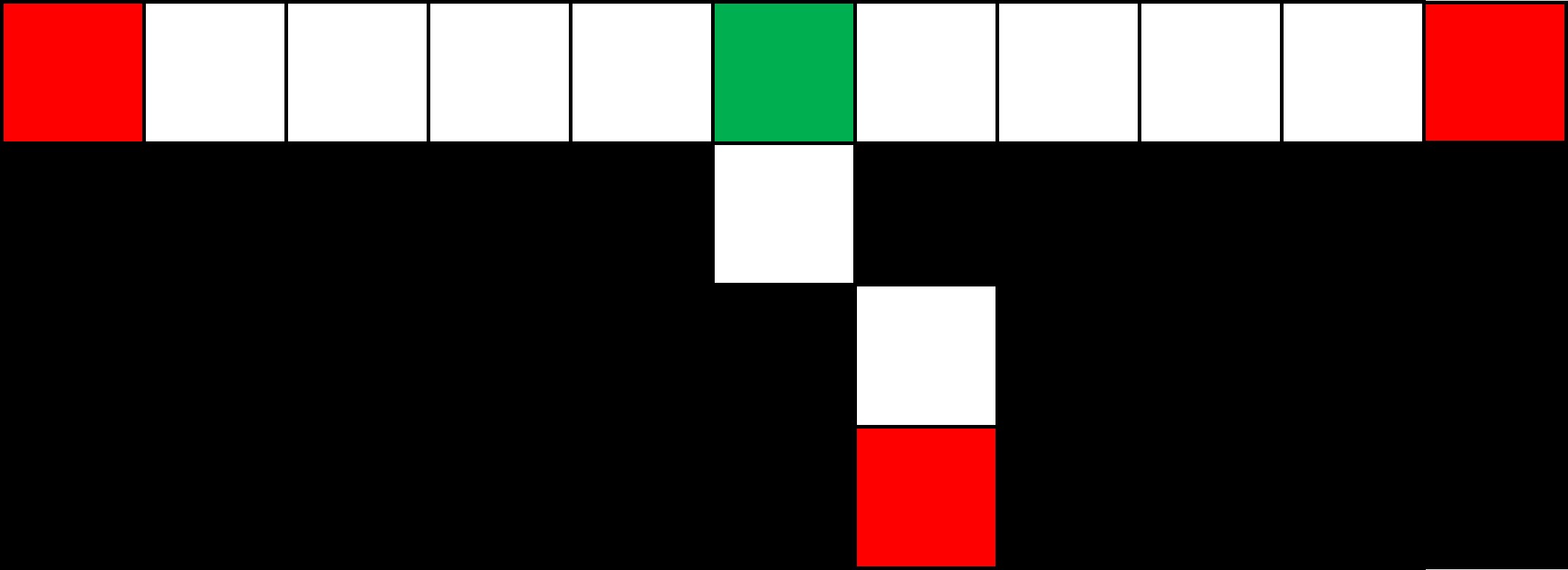}
		\label{fig:pix_cat_d}
	} \\[-3pt] 
	\centering
	\caption{Illustration of skeletonization and pixels categorization}
	\label{fig:pix_cat}
\end{figure}

Thinning is a technique used to extract the skeleton of an object from an image, while preserving its topology \cite{bertrand1995parallel,cheriet2007character, rosenfeld1975characterization}. Thinning algorithms remove as many pixels as possible without affecting the overall topology of the object. This is achieved by removing simple pixels that are not end/boundary pixels, while preserving important information about the geometry of the object. Illustrations of the thinning process are depicted in Figures \ref{fig:pix_cat_a} through \ref{fig:pix_cat_d}.

Parallel thinning algorithms take into account another condition, as removing all simple pixels that are not end pixels simultaneously may alter the topology of the image. A common approach in 2D is to use a directional strategy to delete pixels in parallel  \cite{bertrand1995parallel}, such as removing all border pixels of a particular type that are simple and non-ends at each iteration. Thinning algorithms work by repeatedly converting removable white pixels to black until no more can be done. For this, the Zhang-Suen thinning algorithm \cite{zhang1984fast} is the by far the most widely used technique because of its robustness and ease of implementation in image processing workflows.

\subsubsection{Zhang-Suen Thinning Algorithm}\label{subsec:zhang_thin_alg}

The Zhang-Suen Thinning algorithm is a widely used parallel thinning algorithm that involves removing pixels from an image based on local neighborhood patterns, which is devised in 1984 \cite{zhang1984fast}. For its versatility, the Zhang-Suen algorithm has found recent applications in various fields, including medical imaging \cite{zhang2015retinal}, pattern recognition \cite{chen2018fingerprint,zhao2020handwritten}, robotics \cite{fani2020obstacle}, and structural analysis \cite{bakhshi2019application}. It is a two-pass algorithm that uses two sets of checks to remove pixels from the image, and the flagged pixels are deleted after all pixels in the image have been visited. The resulting image is the skeleton, which meets the criterion of being as thin as possible, connected, and centered. 

\subsubsection{Thresholding and Otsu's Technique}\label{subsec:thresholding}

The input for the skeletonization algorithm is a binary image obtained through thresholding the optimized geometry. This thresholding process binarizes the optimized pixel model based on its relative density or pixel densities, resulting in a binary image where all pixels represent either \(1\) and void \(0\) pixels. The threshold value $\eta$ should be determined so that the given material volume constraint is met. Binary image shown in Figure \ref{fig:skel_cant_ex_0.5_b} is obtained by thresholding the optimized output with only the pixels above a relative density of \( \eta=0.5 \). While it is possible to iteratively determine the threshold value for geometry with the volume even closer to \( V_f \bar{V} \), it is not recommended because of the high computational costs.

Another way to automatically determine a threshold value is by using Otsu’s technique, which is widely used in various applications in computer vision and image processing \cite{talab2016detection, meiburger2016skeletonization}. In the simplest form, Otsu’s thresholding technique processes the input pixel data and returns a single threshold value that separates the pixels into two classes such as 1 and 0. This threshold is determined by minimizing intra-class intensity variance, or equivalently, by maximizing inter-class variance. In the topology-optimized outputs, pixel intensities are equal to the relative densities of each element. Otsu’s thresholding technique can be used to extract the binary image from the topology-optimized model. However, this algorithm performs poorly in the case of high filter radius usage for topology optimization.

\subsubsection{Modified Zhang-Suen Thinning Algorithm}\label{subsec:mod_zhang_thin_alg}
After obtaining the binary image, the skeleton extraction process begins, aiming to preserve the critical load path revealed by topology optimization. To accomplish this, a modified version of the Zhang-Suen thinning algorithm is employed in this work, which includes two key modifications.
\begin{figure}[H] 
	\centering
    \captionsetup{justification=centering}
	\subfloat[Chosen cantilever example]{
		\includegraphics[width=0.49\textwidth]{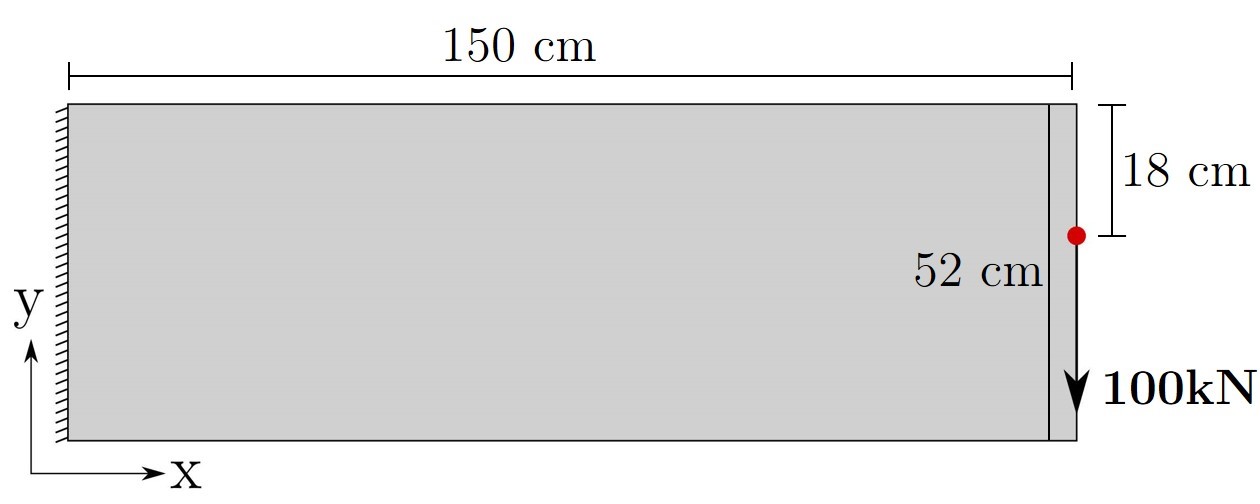}
		\label{fig:skel_cant_ex_0.5_a}
	}
	\centering
	\subfloat[Topology optimized model for $V_f=0.5$]{
		\includegraphics[width=0.49\textwidth]{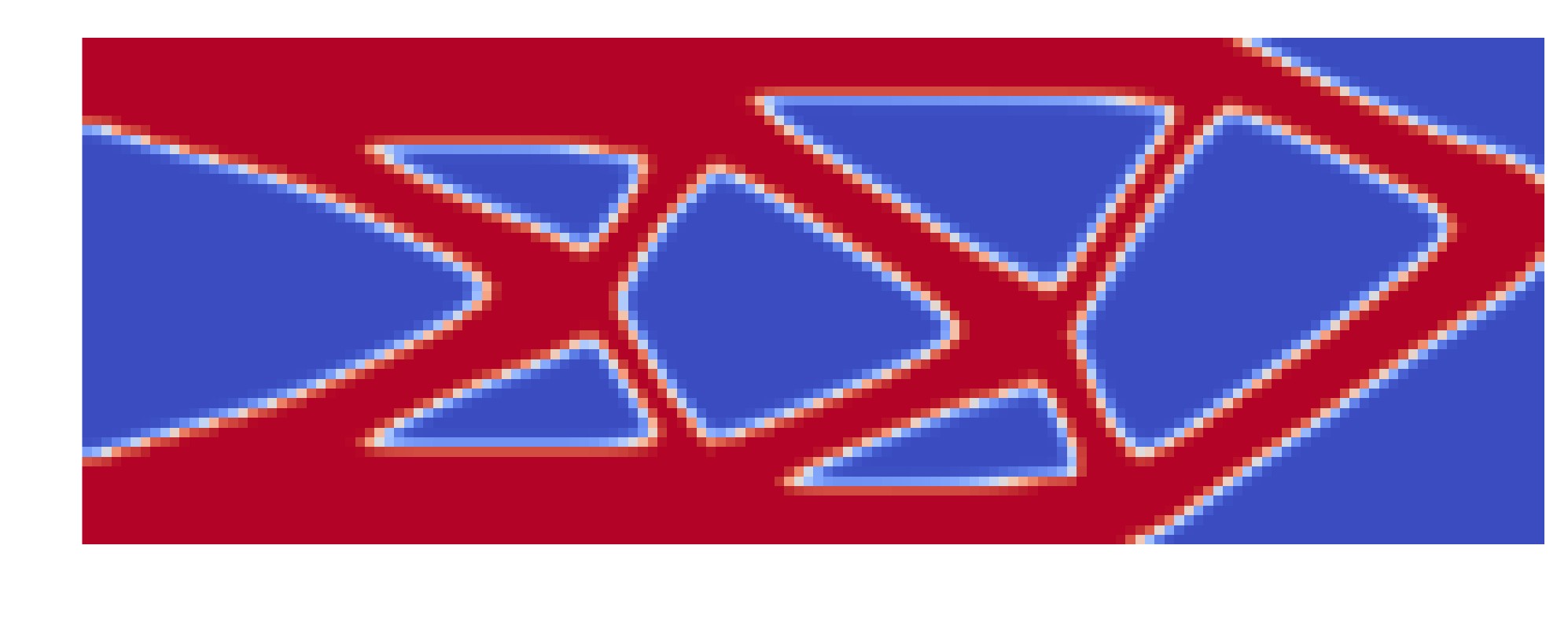}
		\label{fig:skel_cant_ex_0.5_b}
	} \\[-3pt]
        \centering
	\subfloat[Binary image generated for corresponding model in (b)]{
		\includegraphics[width=0.49\textwidth]{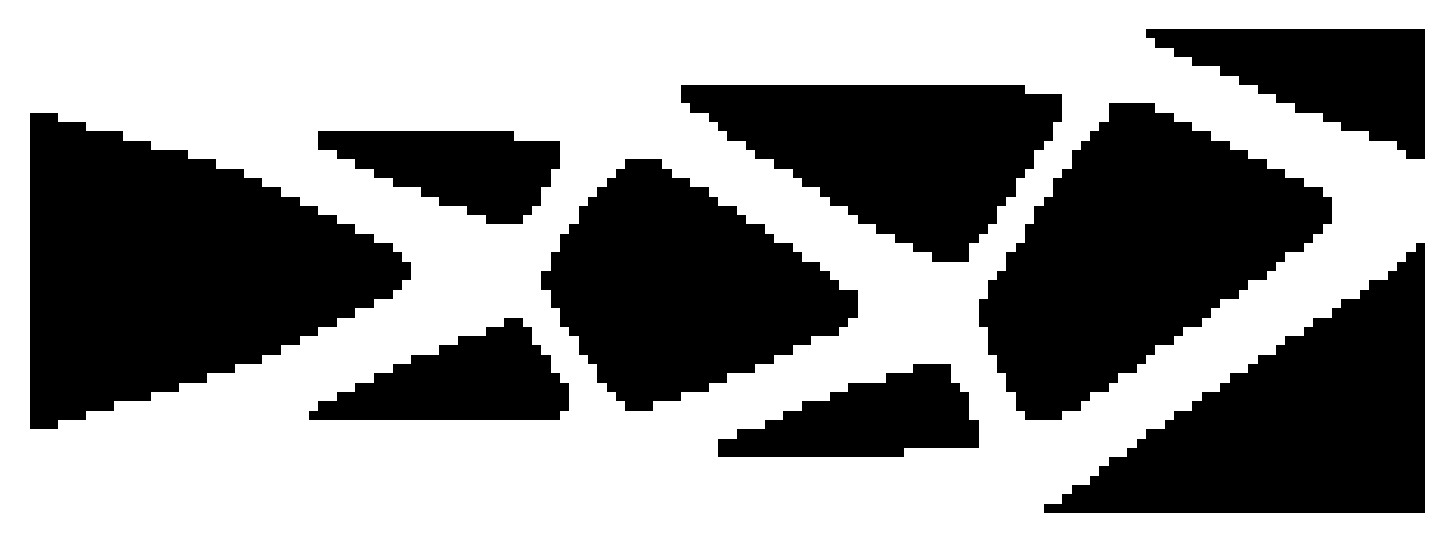}
		\label{fig:skel_cant_ex_0.5_c}
	} 
	\centering
	\subfloat[Direct skeleton obtained for topology optimized binary data in (c)]{
		\includegraphics[width=0.49\textwidth]{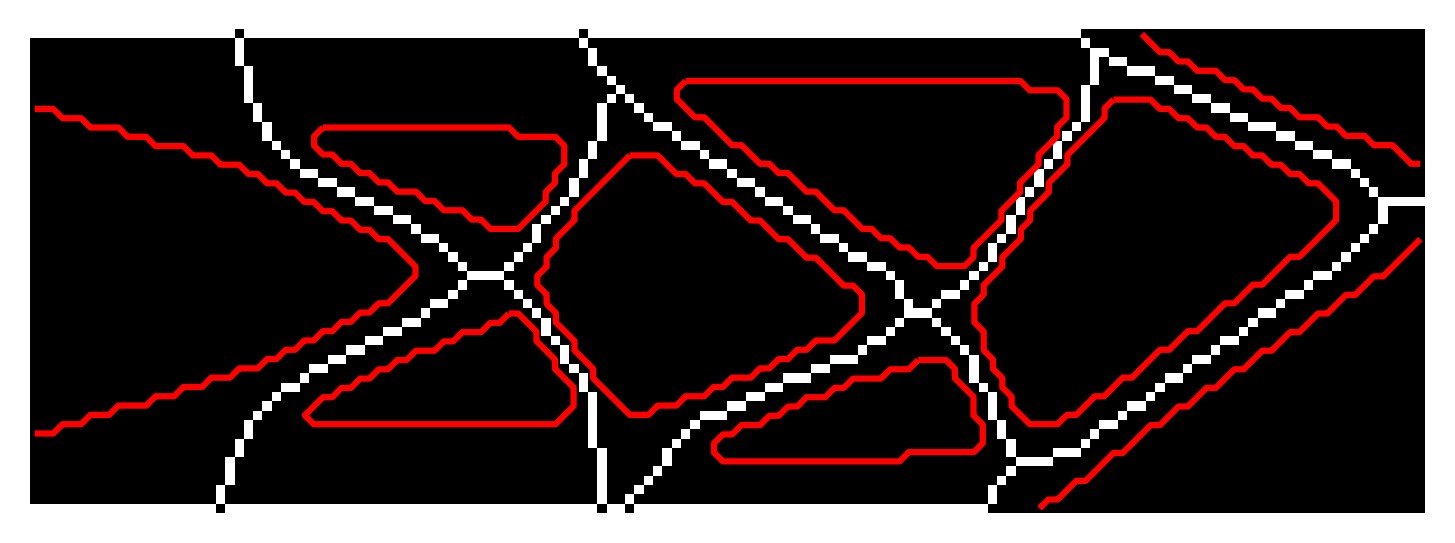}
		\label{fig:skel_cant_ex_0.5_d}
        } \\[-3pt]
  	\centering
	\subfloat[Skeleton obtained after adding extra row and column of pixels around the image]{
		\includegraphics[width=0.49\textwidth]{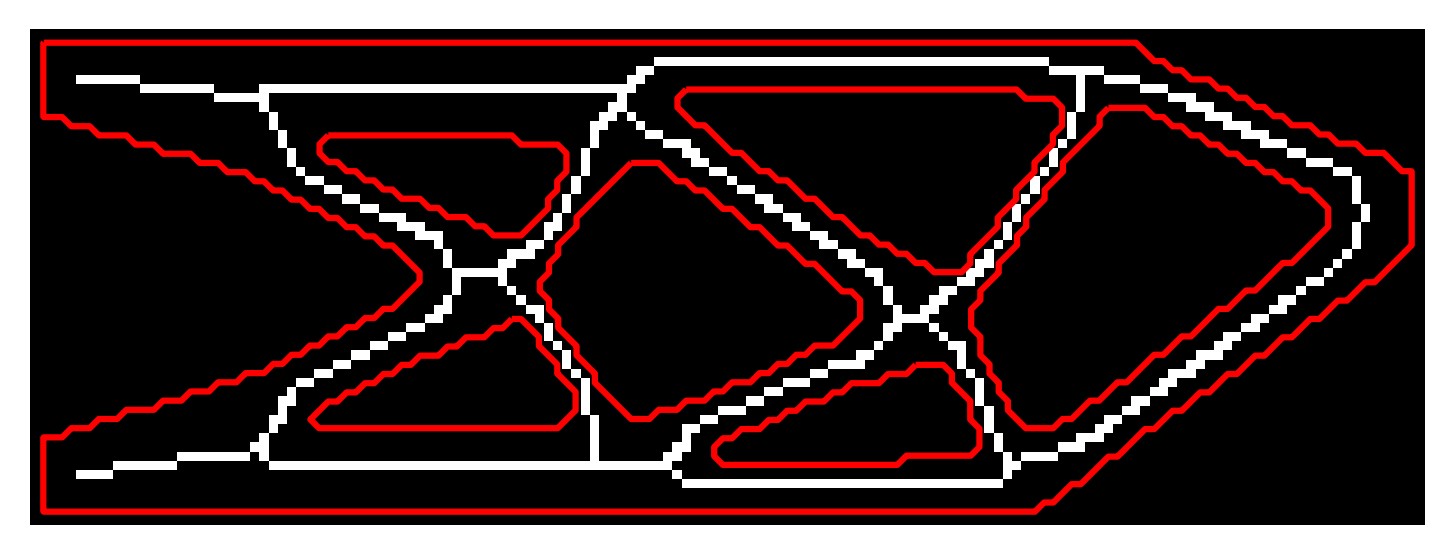}
		\label{fig:skel_cant_ex_0.5t_e}
	} 
        \centering
	\subfloat[Skeleton in (e) after removing added extra row and column of pixels]{
		\includegraphics[width=0.49\textwidth]{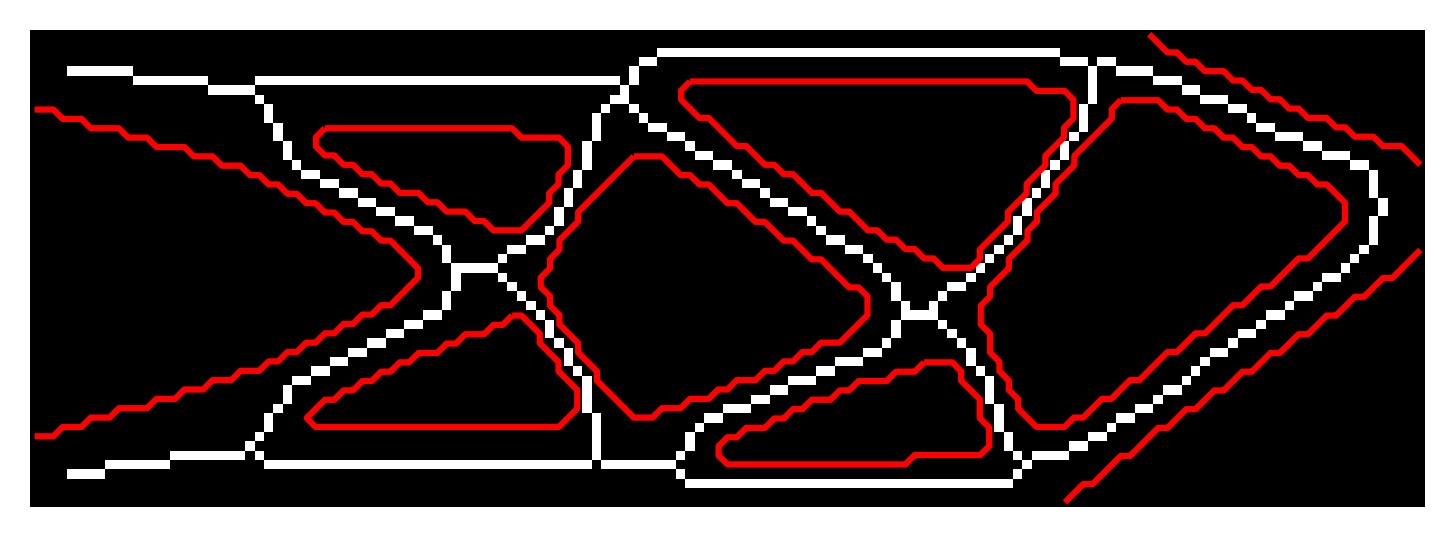}
		\label{fig:skel_cant_ex_0.5_f}
	} \\[-3pt]
	\centering
	\subfloat[Skeleton obtained after tagging boundary pixels (as left mid pixels) and loaded element pixel]{
		\includegraphics[width=0.49\textwidth]{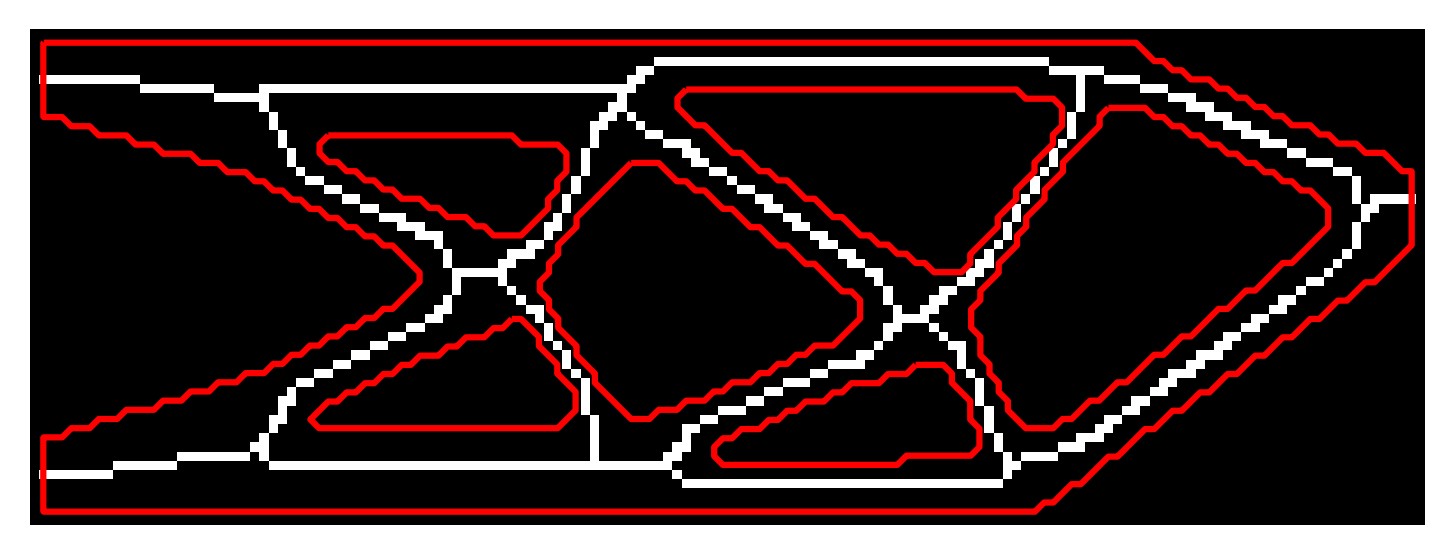}
		\label{fig:skel_cant_ex_0.5_g}
        }
  	\centering
	\subfloat[Skeleton in (g) after removing added extra row and column of pixels]{
		\includegraphics[width=0.49\textwidth]{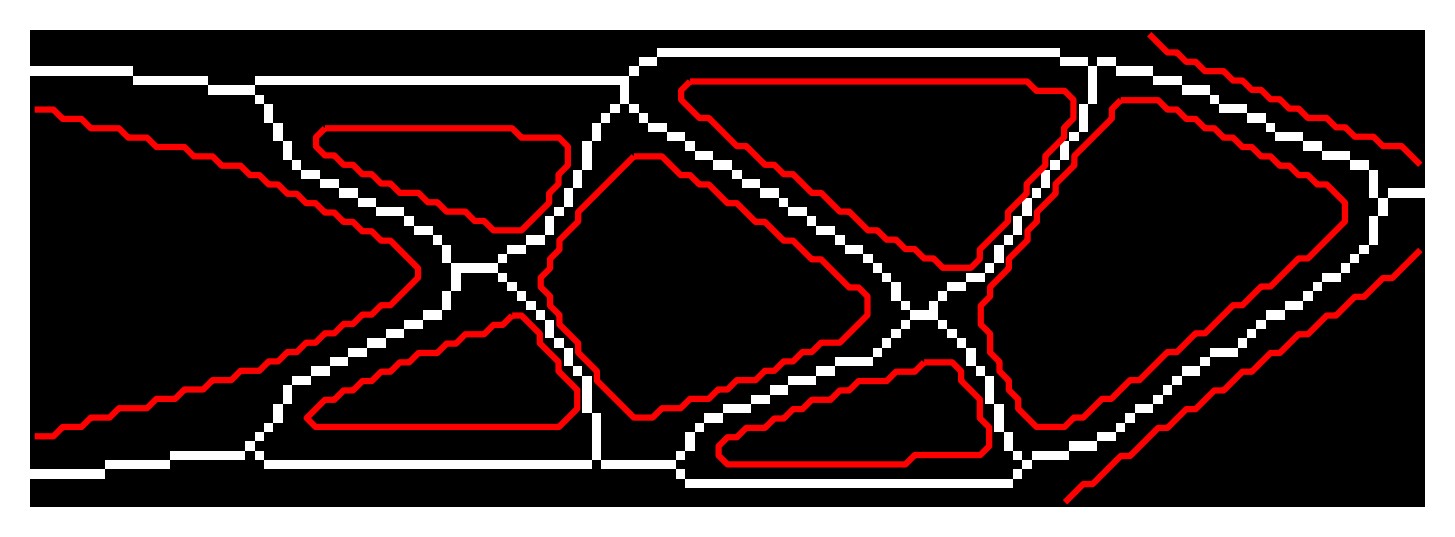}
		\label{fig:skel_cant_ex_0.5_h}
	} \\[-3pt]
        \centering
	\subfloat[Skeleton obtained after tagging boundary pixels (as left top and bottom pixels) and loaded element pixel]{
		\includegraphics[width=0.49\textwidth]{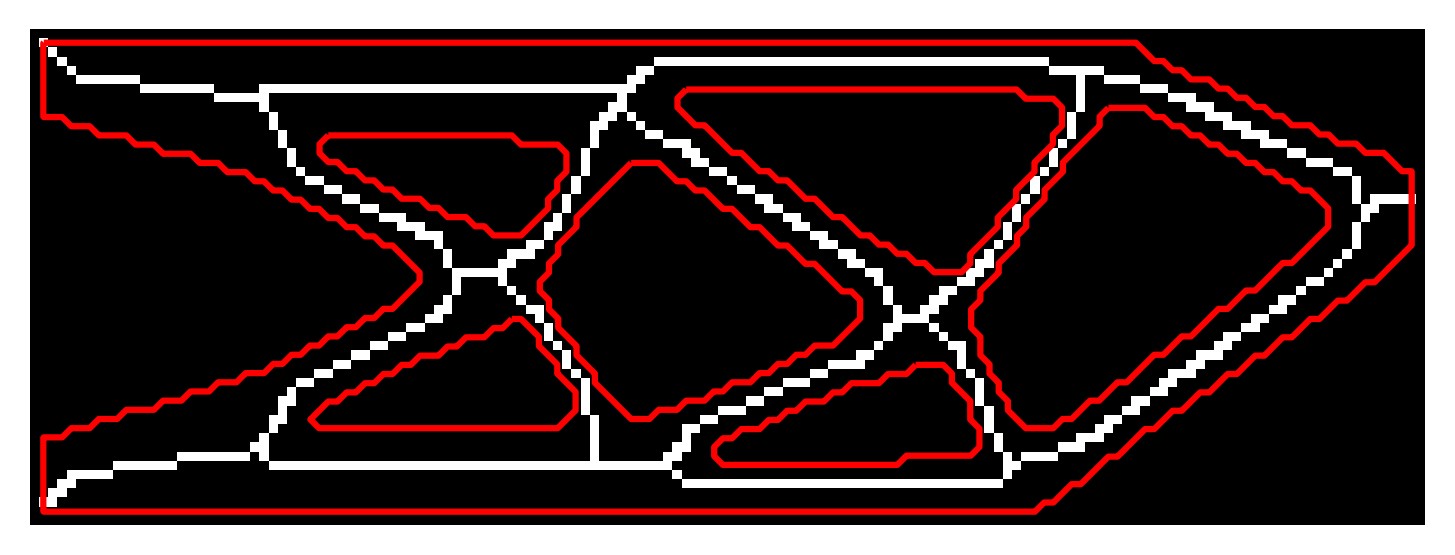}
		\label{fig:skel_cant_ex_0.5_i}
	}
	\centering
	\subfloat[Skeleton in (i) after removing added extra row and column of pixels]{
		\includegraphics[width=0.49\textwidth]{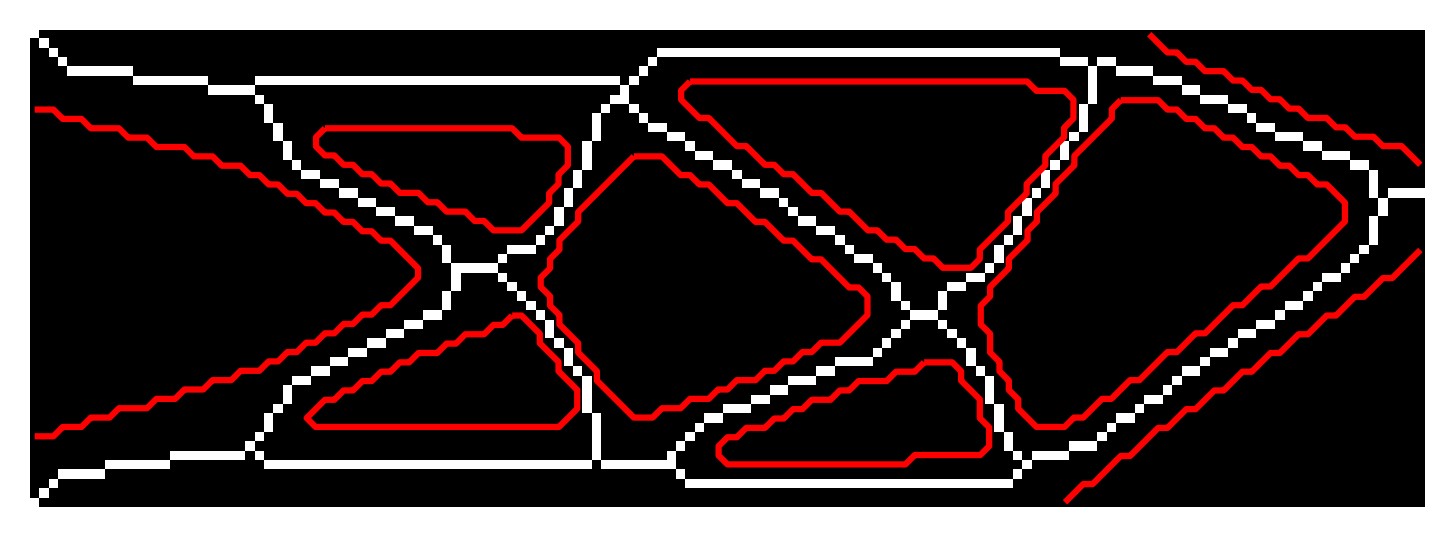}
		\label{fig:skel_cant_ex_0.5_j}
	} \\[-3pt] 
	\centering
	\caption{Skeletonization of a cantilever plate example ($V_f$ = $0.5$)}
	\label{fig:skel_cant_ex_0.5}
\end{figure}
The first modification involves the introduction of a tagging mechanism to mark certain features in the structure as non-removable. This tagging process ensures that important structural elements are preserved during the thinning process. The algorithm initially detects and flags all simple untagged non-end edge pixels. A rechecking process is then performed to maintain the connectivity of the remaining pixels. The flagged pixels are iteratively eliminated until no further removal is possible, resulting in a 1-wide pixel chain that represents the skeleton of the structure.

The second modification addresses the issue of pixel connectivity near the image borders. To ensure accurate skeleton generation for pixels located near the borders, an extra row and column of pixels are added around the binary image before applying the thinning algorithm. This additional padding guarantees that the medial axis is correctly generated. Once the skeletonization process is completed, the extra rows and columns are removed, restoring the binary image to its original size. In the case of the binary image shown in Figure \ref{fig:skel_cant_ex_0.5_c}, the modified Zhang-Suen algorithm is applied after adding the extra borders, resulting in the skeletonized image depicted in Figure \ref{fig:skel_cant_ex_0.5_f}. This skeleton accurately represents the underlying structure while preserving the important features. The proposed modifications to the Zhang-Suen algorithm are employed in the skeletonization of the topology-optimal model. The resulting skeleton is processed further to extract the frame model for subsequent analysis.

\subsection{Frame Extraction}\label{subsec:frame_skel}
This section focuses on the process of converting topology-preserved pixel chains (skeleton) into a simplified 2D frame model using graph algorithms. By characterizing the pixel chain as a graph model with nodes and edges, it is possible to obtain simpler geometric representations. Therefore, the skeleton is transformed into an unweighted undirected graph, leveraging its inherent properties \cite{babin2018skeletonization}. While it is feasible to derive element cross-sections directly from the pixel chain widths in the topology-optimized model, this study opts to determine cross-sections through size optimizations of the frames. This approach is more practical as it avoids the burden of retaining unnecessary chain widths that will be optimized and modified later on \cite{YIN2020113102}.

\subsubsection{Graph Construction from Pixel Chain}\label{subsec:graph_cons}

\begin{figure}[H] 
	\centering
    \captionsetup{justification=centering}
	\subfloat[Skeleton (medial axis of Figure~\ref{fig:pix_cat_a})]{
		\includegraphics[width=0.49\textwidth]{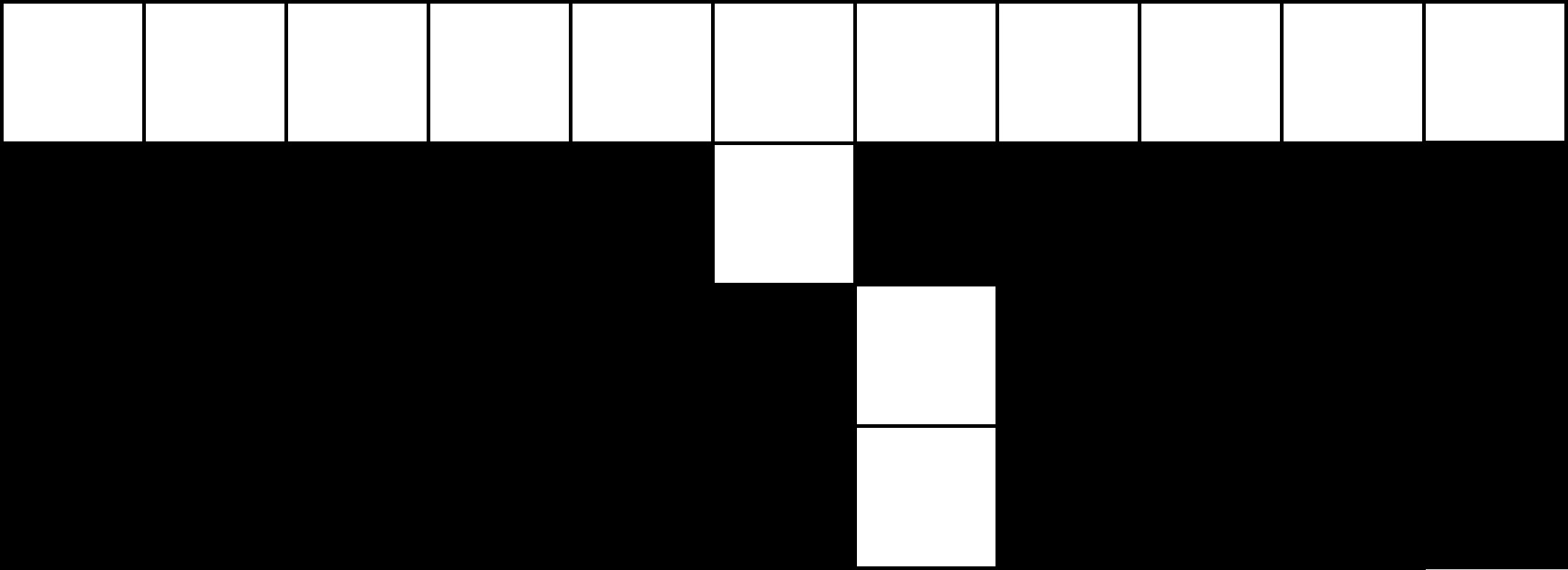}
		\label{fig:compact_graph_a}
	}
	\centering
	\subfloat[Directly converted graph model from skeleton]{
		\includegraphics[width=0.49\textwidth]{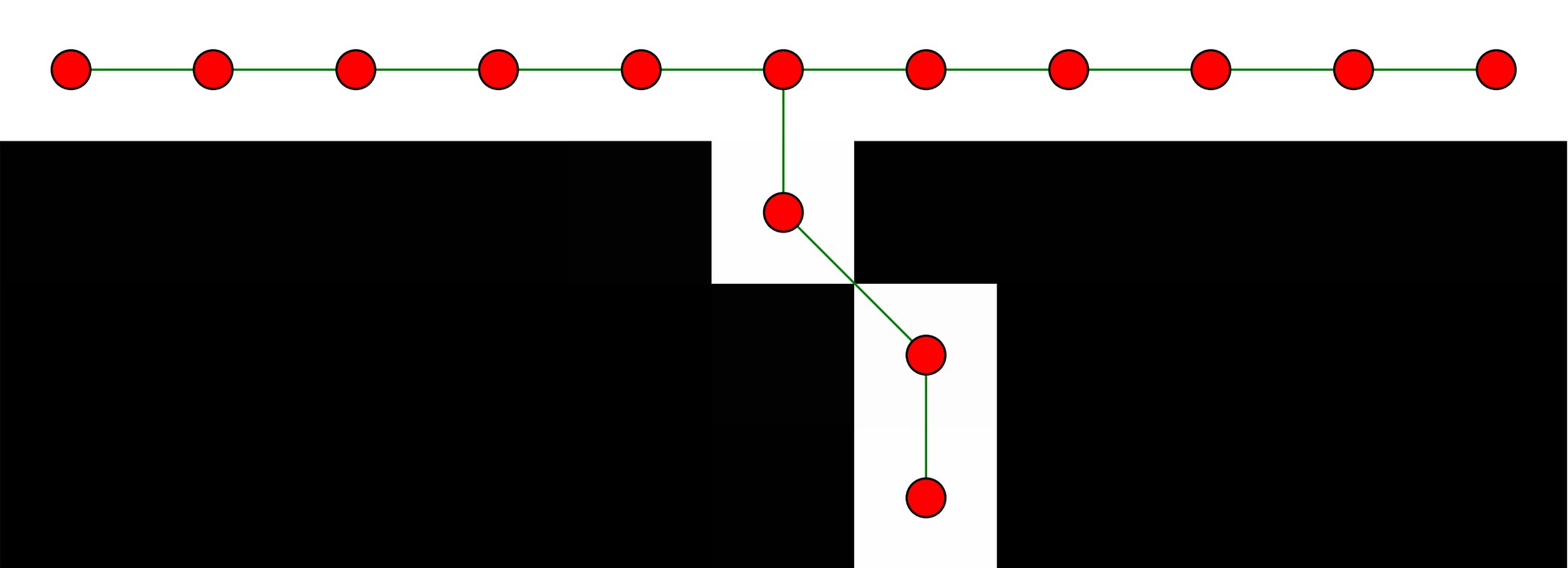}
		\label{fig:compact_graph_b}
	} \\[-3pt]
        \centering
	\subfloat[Identified important pixels in skeleton (Red–End, Green–Joint, White–Regular)]{
		\includegraphics[width=0.49\textwidth]{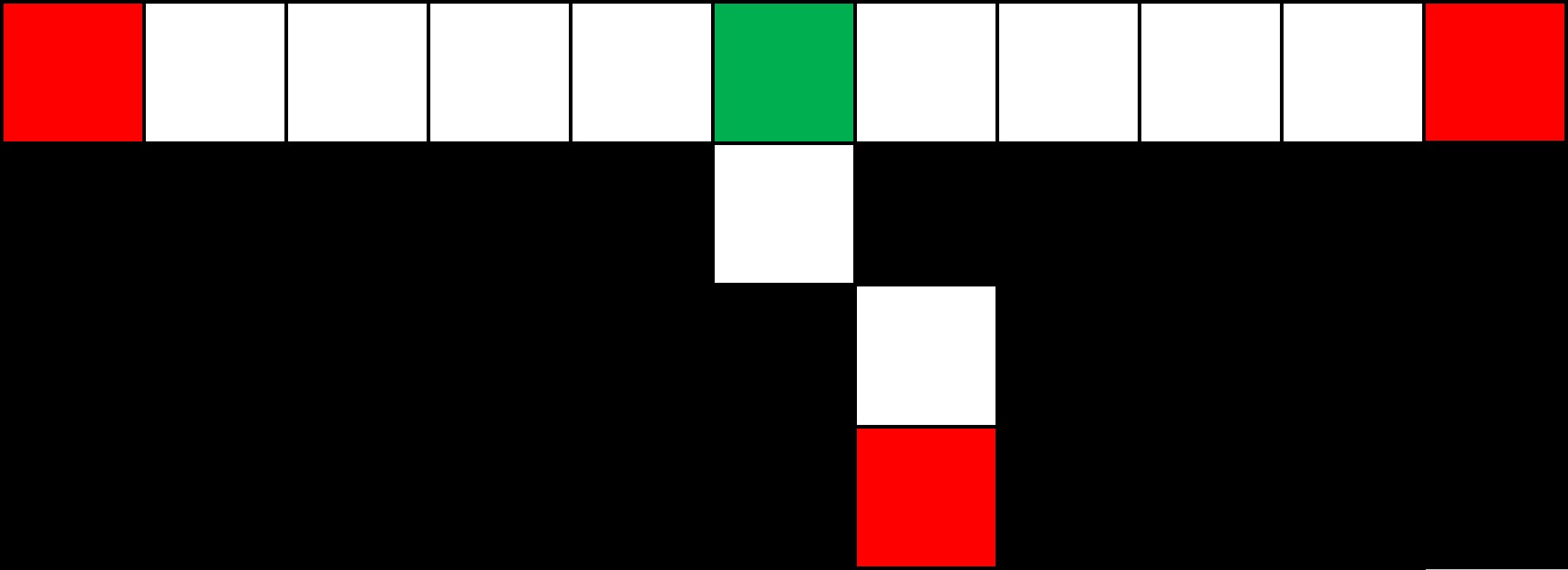}
		\label{fig:compact_graph_c}
	}
	\centering
	\subfloat[Compact graph model generated using important pixels]{
		\includegraphics[width=0.49\textwidth]{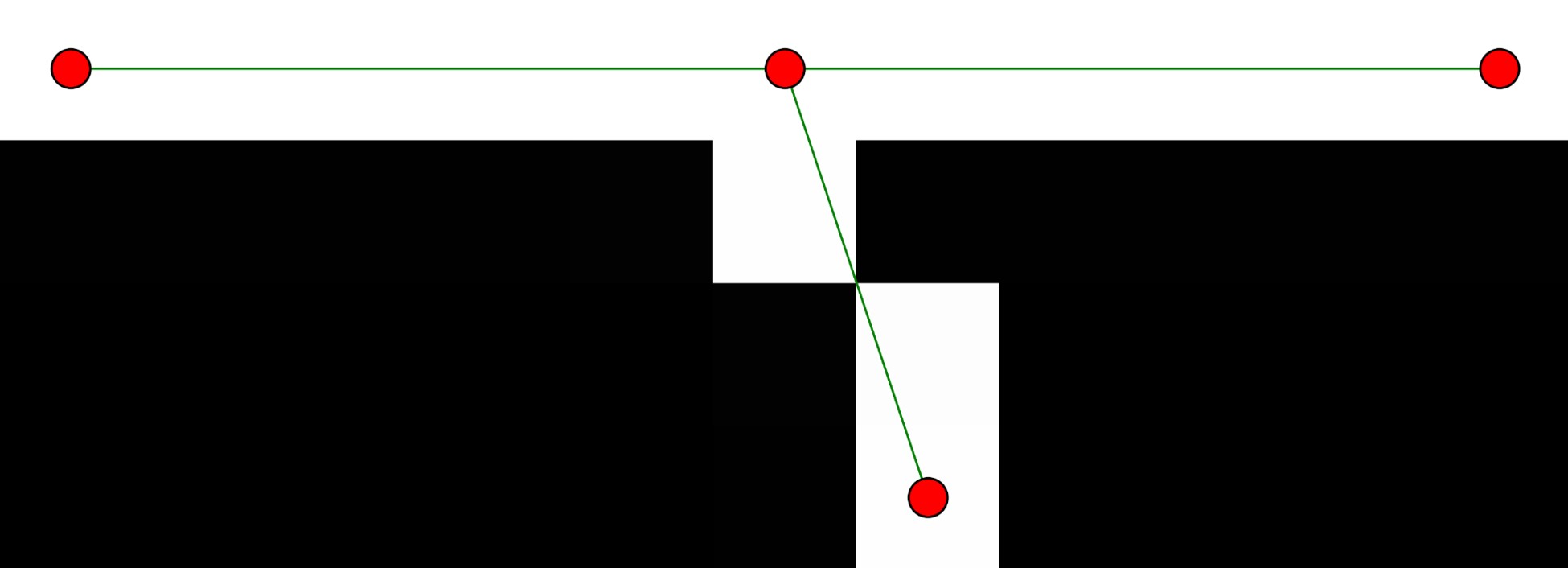}
		\label{fig:compact_graph_d}
	} \\[-3pt] 
	\centering
	\caption{Compact graph model generation using important pixels}
	\label{fig:compact_graph}
\end{figure}

A simple graph $G$ is represented as $G = (V, E)$, where $V$ denotes the set of nodes and $E$ represents the set of edges \cite{wilson1979introduction, west2001introduction}. The graph model is constructed by characterizing the pixel chain using a graph, consisting of nodes and edges. Initially, each node is assigned coordinates corresponding to the centroid of the corresponding pixel. However, directly converting the edges of the graph into frame members, as illustrated in Figure \ref{fig:compact_graph_b}, is neither feasible nor practical. This is due to the creation of short and ineffective members in the structure, leading to increased complexity. Therefore, it is necessary to obtain a more compact graph, depicted in Figure \ref{fig:compact_graph_d}, by extracting only the edges that connect the important pixels.

In this context, pixels are categorized into three types based on their topology: regular, end, and joint pixels. These types are distinct from the concept of tagged pixels, which serve the purpose of representing load and boundary conditions and are treated separately. Similarly, nodes can also be classified based on the number of edges they are shared by, into regular, end, and joint nodes, with their degree representing the number of edges close to the node. The information on important pixels and nodes is integrated into the graph model for operations on the compact graph model.
The end pixel, joint pixel, and tagged pixels are considered featured pixels, which provide the basic compact layout of the graph. The process of identifying edges involves extracting two pixels in featured pixels that are connected through a path consisting of untagged regular pixels. Once all pixels have been classified, edges are identified by starting from the joint pixel and moving along its 8-neighbors until a featured pixel is reached \cite{YIN2020113102, babin2018skeletonization}. A graph model made with these edges and nodes will produce a compact graph needed for the further process as shown in Figure \ref{fig:compact_graph}. For interested readers, a similar work of compact graph generation without the implementation of tagging is presented by Babin et al. \cite{babin2018skeletonization}. 

The compact graph model obtained is still not efficient for structural applications because of the short branches, as depicted in Figure \ref{fig:compact_graph_b}. To address this, the model is modified by pruning redundant branches and contracting short edges using a merge ratio defined by the user \cite{YIN2020113102}. The modification results in a more efficient graph model for structural applications.
\begin{figure}[H] 
	\centering
    \captionsetup{justification=centering}
	\subfloat[Direct conversion of skeleton to graph model (complex and inefficient)]{
		\includegraphics[width=0.49\textwidth]{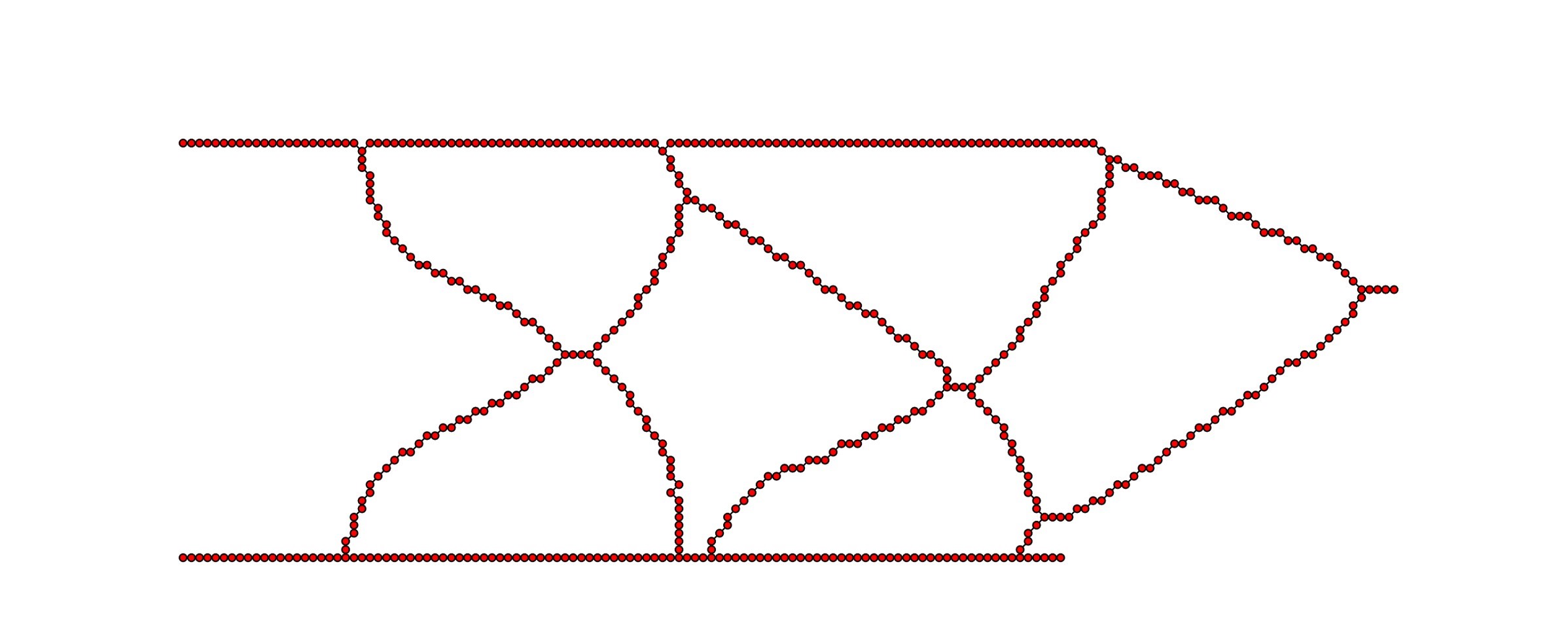}
		\label{fig:skel_to_frame_a}
	}
	\centering
	\subfloat[Initial compact graph model extracted from skeleton]{
		\includegraphics[width=0.49\textwidth]{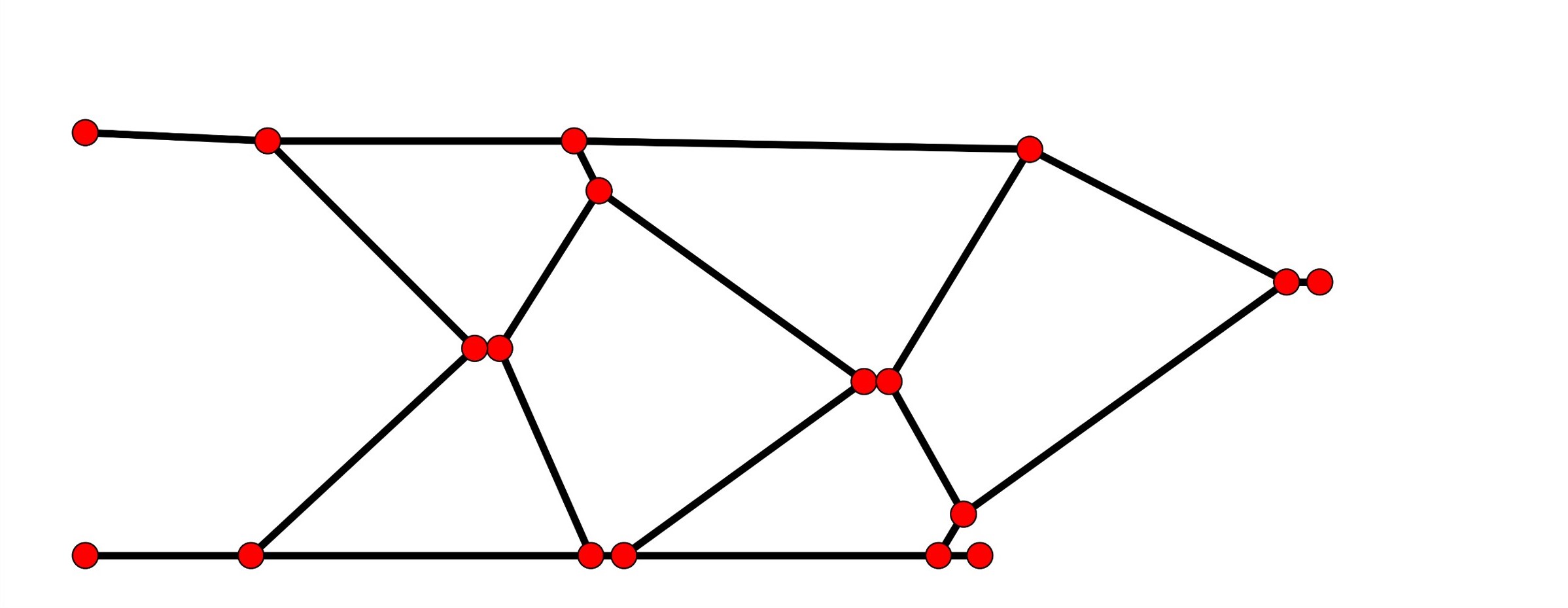}
		\label{fig:skel_to_frame_b}
	} \\[-3pt]
    \centering
	\subfloat[Identified redundant branch for pruning]{
		\includegraphics[width=0.49\textwidth]{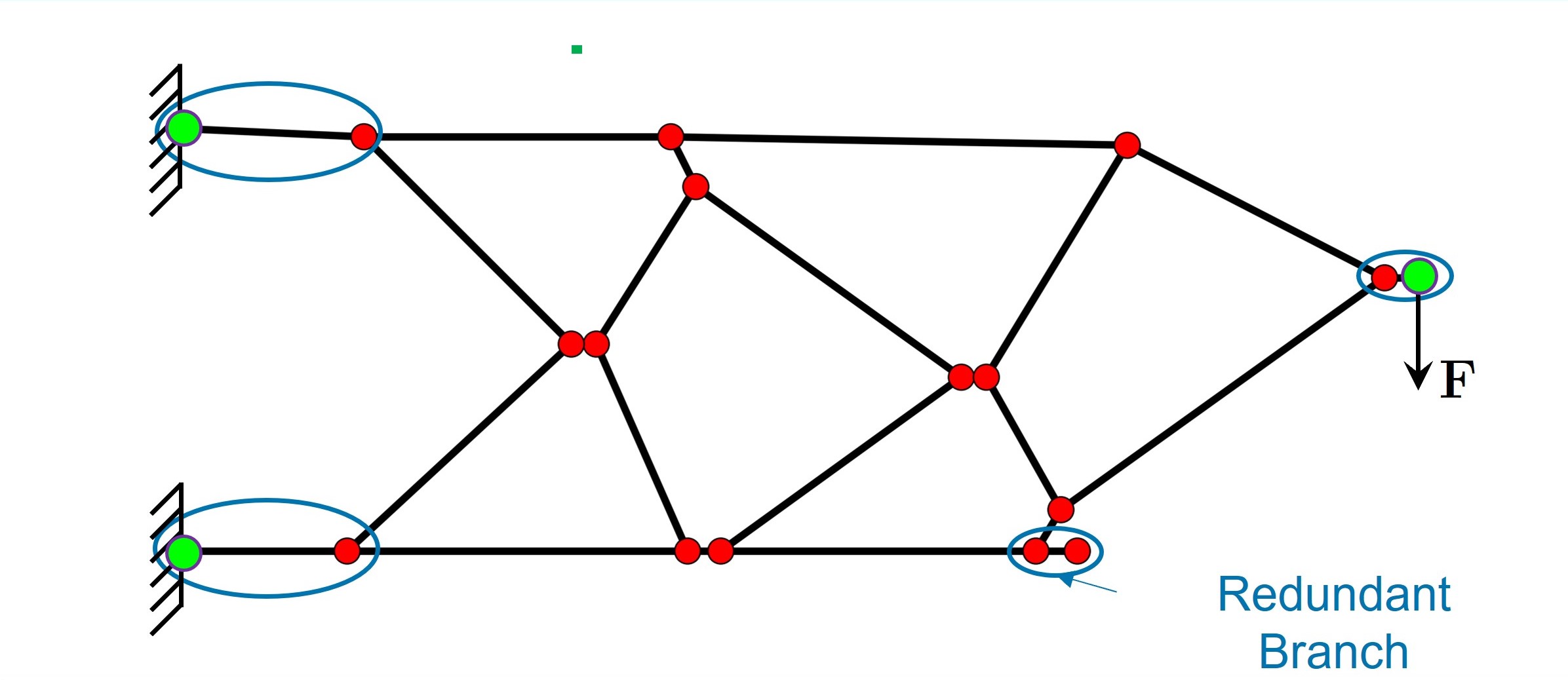}
		\label{fig:skel_to_frame_c}
	}
	\centering
	\subfloat[Graph model after pruning redundant members]{
		\includegraphics[width=0.49\textwidth]{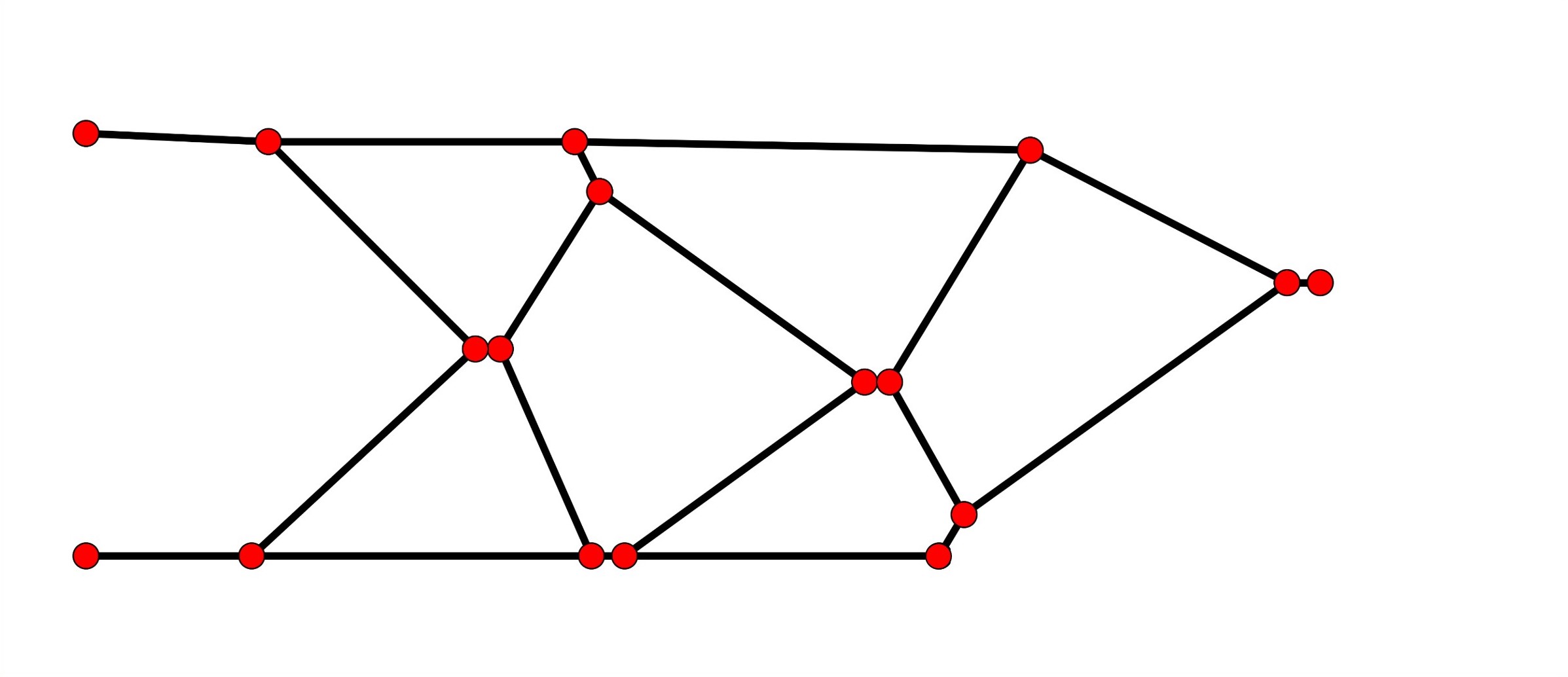}
		\label{fig:skel_to_frame_d}
	} \\[-3pt] 
    \centering
	\subfloat[Identified short members for contraction]{
		\includegraphics[width=0.49\textwidth]{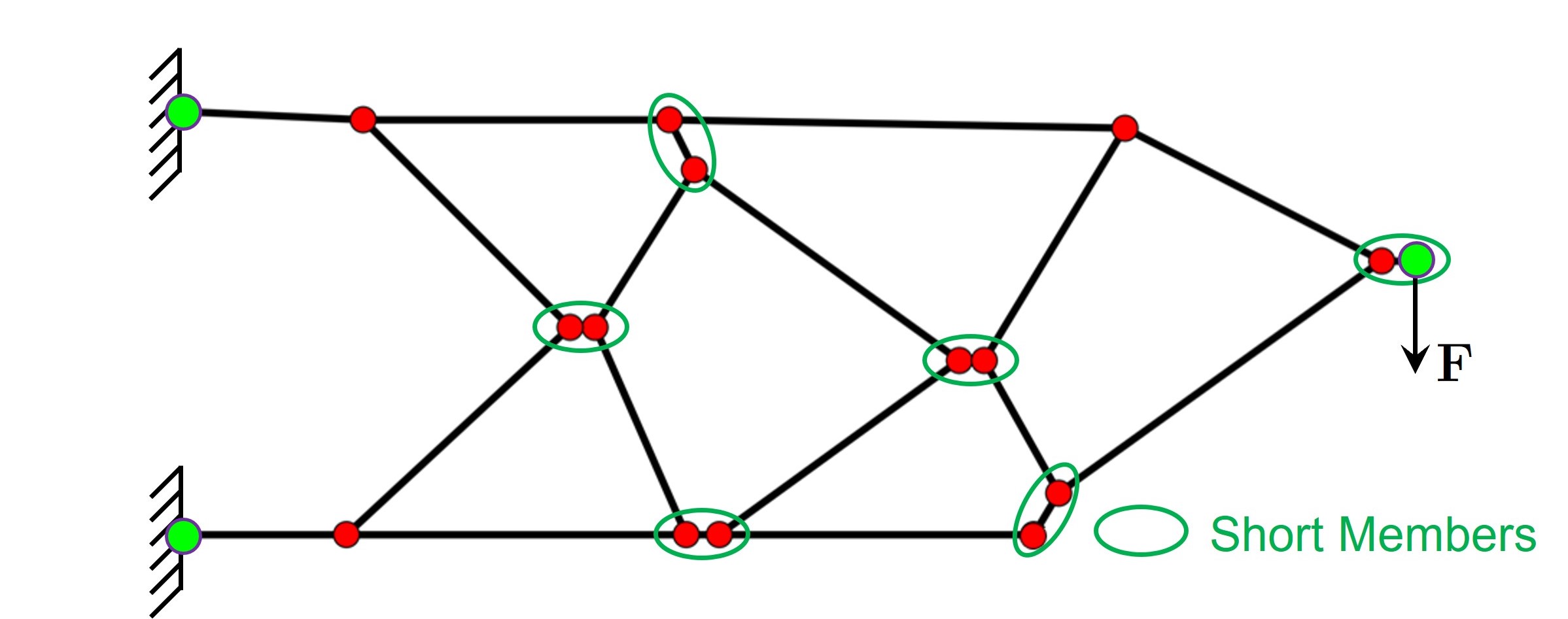}
		\label{fig:skel_to_frame_e}
	}
	\centering
	\subfloat[Graph model after contracting short members]{
		\includegraphics[width=0.49\textwidth]{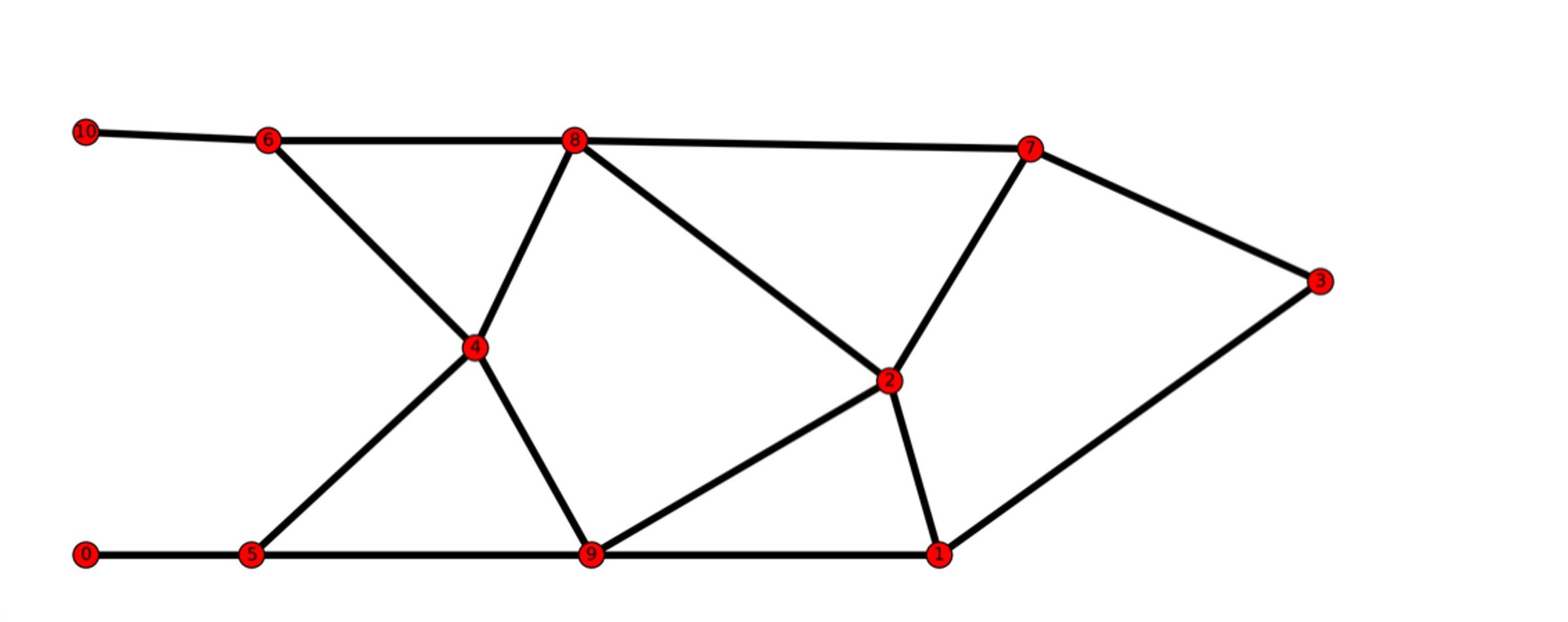}
		\label{fig:skel_to_frame_f}
	} \\[-3pt] 
	\centering
	\caption{Skeleton to frame generation using graph model}
	\label{fig:skel_to_frame}
\end{figure}
\subsubsection{Pruning Redundant Edges and Contraction of Short Edges}\label{subsec:pruning&contraction}

Skeletons from the medial axis transformations may have extra branches generated. To identify useful structural frame members using a graph model, redundant edges can be pruned instead of editing the pixel chain. This pruning is done by removing edges with a degree of one that are not connected to any tagged nodes. An illustration of the pruning of a redundant branch is shown in Figure \ref{fig:skel_to_frame_c} and \ref{fig:skel_to_frame_d} where tagged nodes are used to identify only the redundant branches without detecting the edges adjacent to loading and boundary regions.

\begin{figure}[H] 
	\centering
    \captionsetup{justification=centering}
	\subfloat[Graph model before constraining the angle]{
		\includegraphics[width=0.49\textwidth]{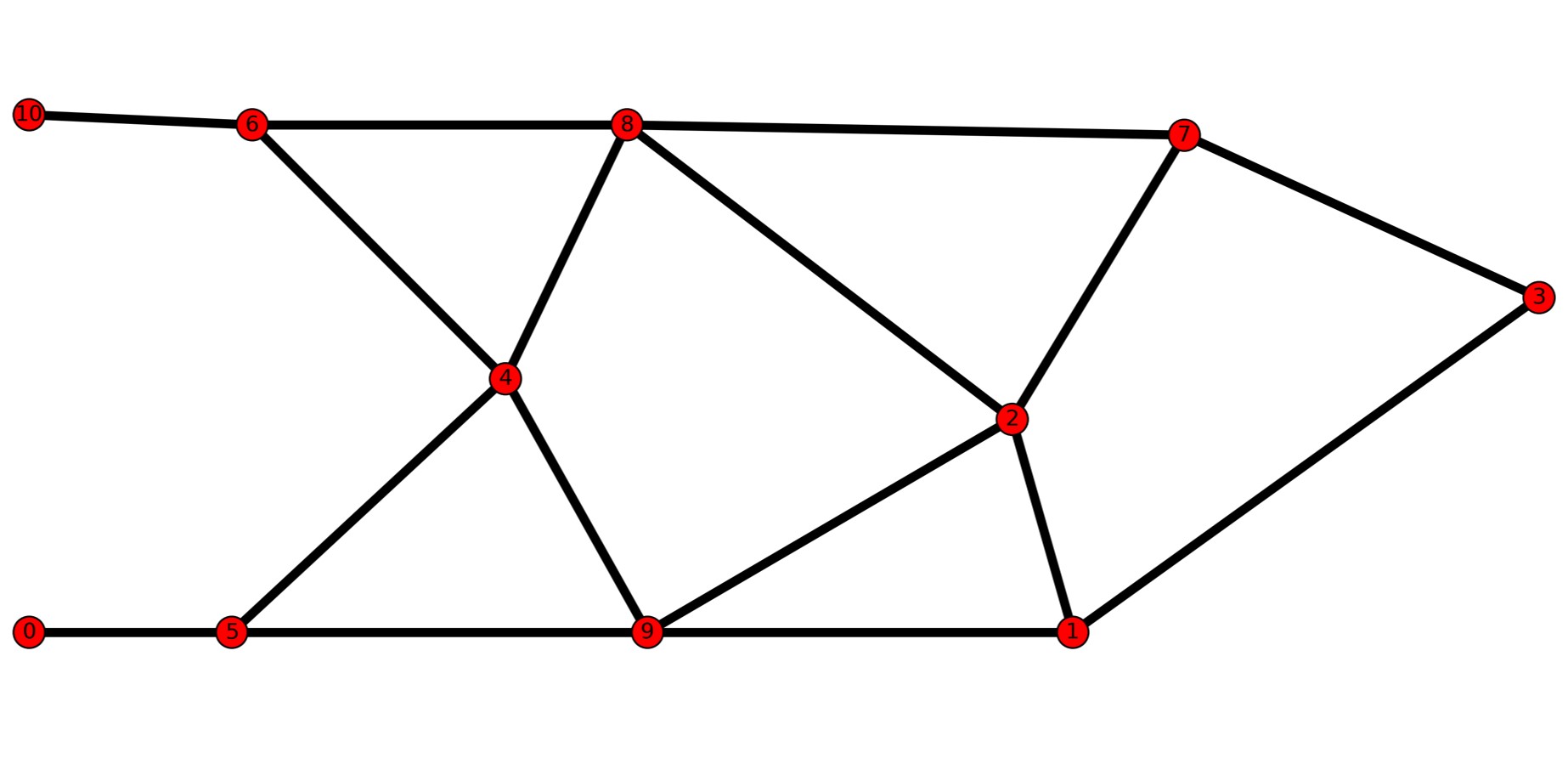}
		\label{fig:cons_angle_a}
	}
	\centering
	\subfloat[Graph model after constraining the angle]{
		\includegraphics[width=0.49\textwidth]{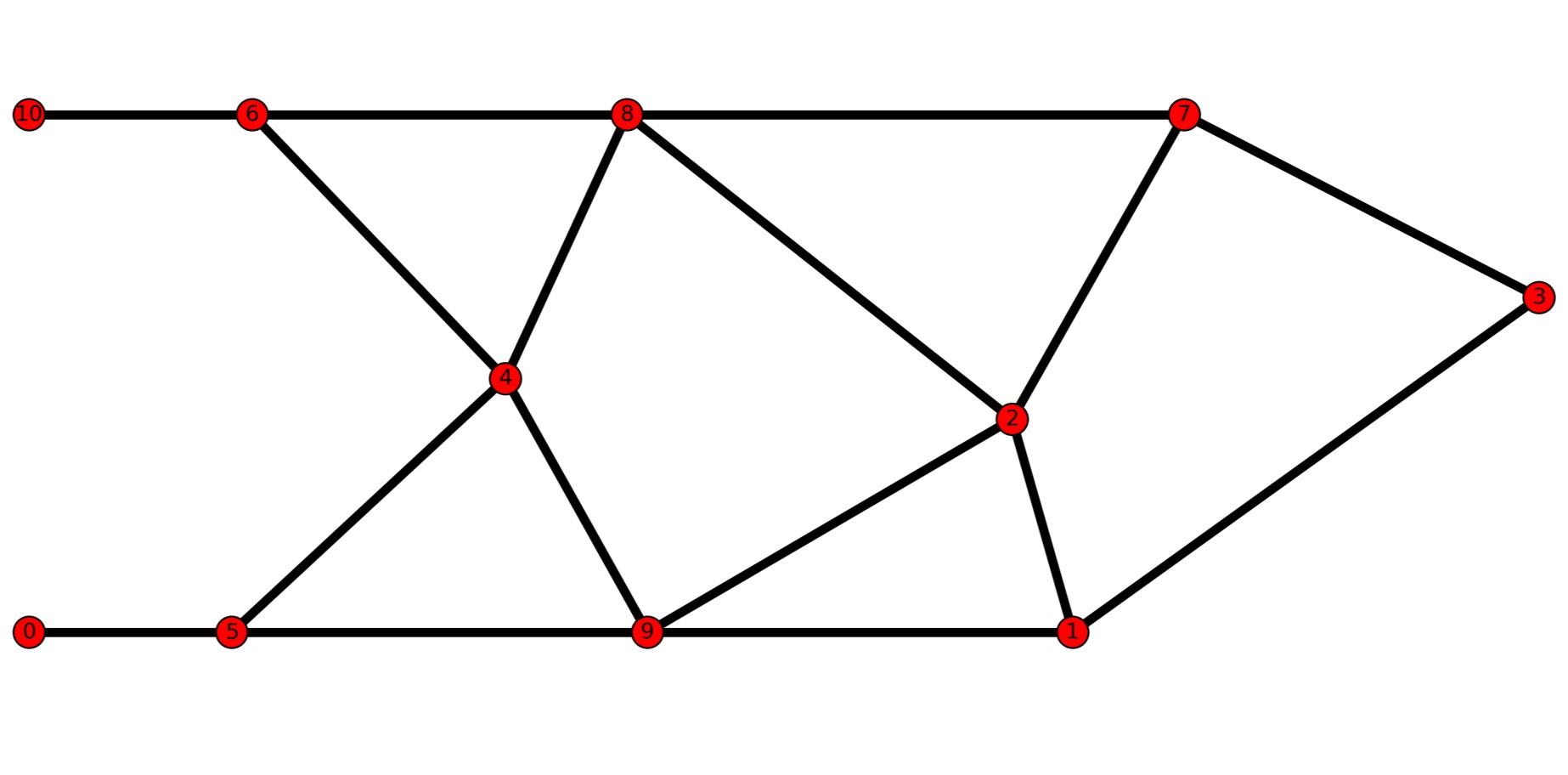}
		\label{fig:cons_angle_b} 
        }
	\centering
	\caption{Constraining the angle of the members}
	\label{fig:cons_angle}
\end{figure}

To address short edges, a merge ratio ($\zeta$) must be defined to identify and contract them. Short edges are those that are shorter than the merge ratio of the total length of all members sharing the same node \cite{YIN2020113102}. The end nodes of the short edges are then contracted until none remain in the graph model. This study uses a merge ratio of 0.05 to 0.1. After pruning and edge contraction, graph models can be optimized for structural analysis applications using layout optimization. This further optimizes bending stresses in the members \cite{zhang2022nodal}. Angle constraining can also be used to obtain a more efficient frame model by constraining small-angled frames. Figure \ref{fig:cons_angle} shows an example where top and bottom frames are horizontally aligned using a small constrain angle of 10°.
%

\subsection{Sequential Size and Layout Optimization}\label{subsec:seq_size&layout_opt}

The extracted structural frame model from the skeleton typically represents only the initial layout of the optimized load path. The optimality achieved through topology optimization is partially compromised in the initial frame model. This compromise is inevitable since skeletonization and frame extraction are unaware of the  structural design considerations. Therefore, sequential size and layout optimizations are carried out to recover the loss in optimality.

The boundary and loading conditions applied to the frame are identical to those used in topology optimization. Subsequently, the size and layout of the members are optimized sequentially. In size optimization, the cross-sectional area of the circular members is updated, while in layout optimization, the nodal coordinates of the members are updated. The optimization problems are solved using both Method of Moving Asymptotes (MMA) and Sequential Quadratic Programming (SQP) algorithms, and their results are compared to reveal computational performances in generating the optimal designs.

In size optimization, solid circular cross-sections are employed in this study. Initially, a constant cross-sectional area is assumed for all members. This area is determined based on the prescribed material volume of the frame structure. The material volume is calculated by considering the volume fraction \(V_f\) used in topology optimization and the assigned thickness value. Since this is a 2D problem, the user needs to define the thickness of the domain to obtain the volume \(\bar{V} \) used for size optimization. Throughout the size and layout optimization, the volume of the frame structure is constrained to be equal to the volume defined by \( V_f \bar{V} \).

\begin{figure}[H]
\centering
\includegraphics[width=0.99\textwidth]{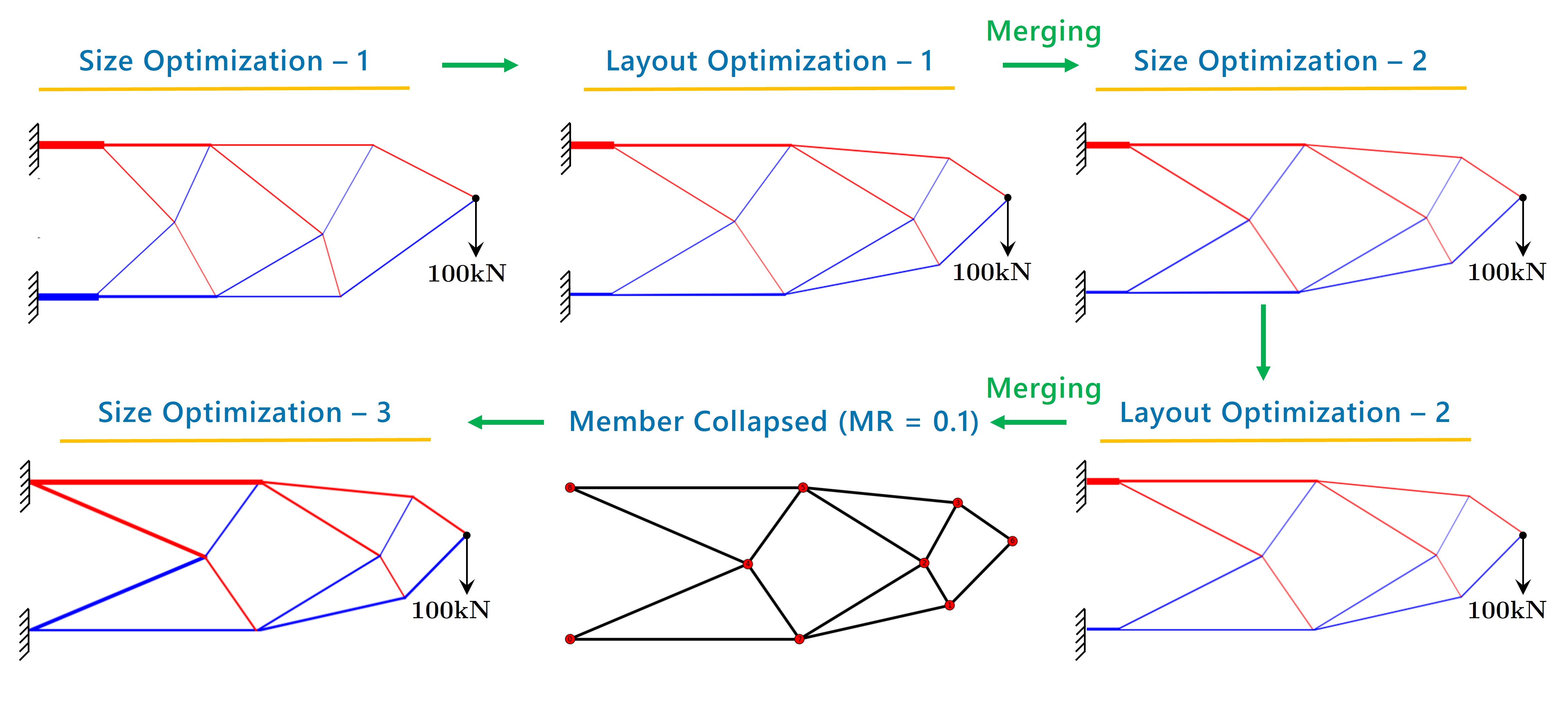}
\\[-12pt]
\caption{Sequential size and layout optimization of a cantilever example. Members in tension and compression are represented by red and blue color respectively}
\label{fig:seq_size_lay}
\end{figure}
The sizes of the members are bounded by defining upper and lower bounds according to the user requirements. However, these bounds can be omitted, allowing the sizes to be determined by the algorithm. Additionally, any other constraints related to the member cross-sections can be enforced as constraints. For example, the allowable member cross-sections can be constrained based on the allowable stress in tension and compression, as specified by the relevant standards and codes of practice.
\begin{figure}[h] 
	\centering
	\subfloat[]{
		\includegraphics[width=0.475\textwidth]{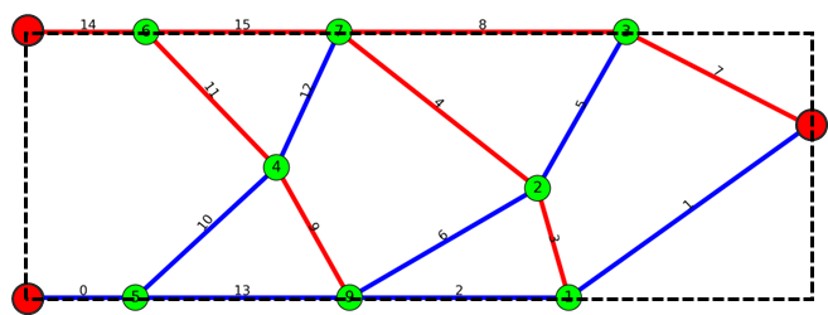}
		\label{fig:lay_constr_a}
	}
	\centering
	\subfloat[]{
		\includegraphics[width=0.475\textwidth]{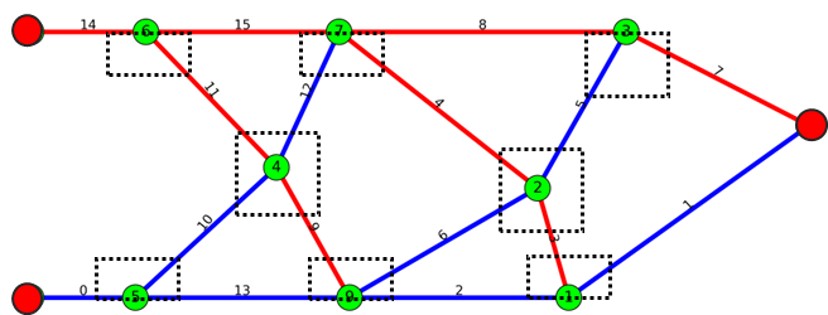}
		\label{fig:lay_constr_b} 
        }
        \\[-6pt]
	\centering
	\caption{Illustration of a Box-bound constraint for layout optimization. (a) Box-bound constraint as design domain for all varying nodes, (b) Box-bound constraint separately for each varying node}
	\label{fig:lay_constr}
\end{figure}
In layout optimization, the upper and lower bounds (simple box-bound constraints) need to be prescribed. It can be given as either the whole design domain or small allowable regions for each node. In the first case, all the varying nodes except the tagged boundary and loaded node are allowed to vary within the whole design domain but in the second case, it is allowed to vary within a small rectangular region.
\begin{figure}[h] 
	\centering
	\subfloat[]{
		\includegraphics[width=0.475\textwidth]{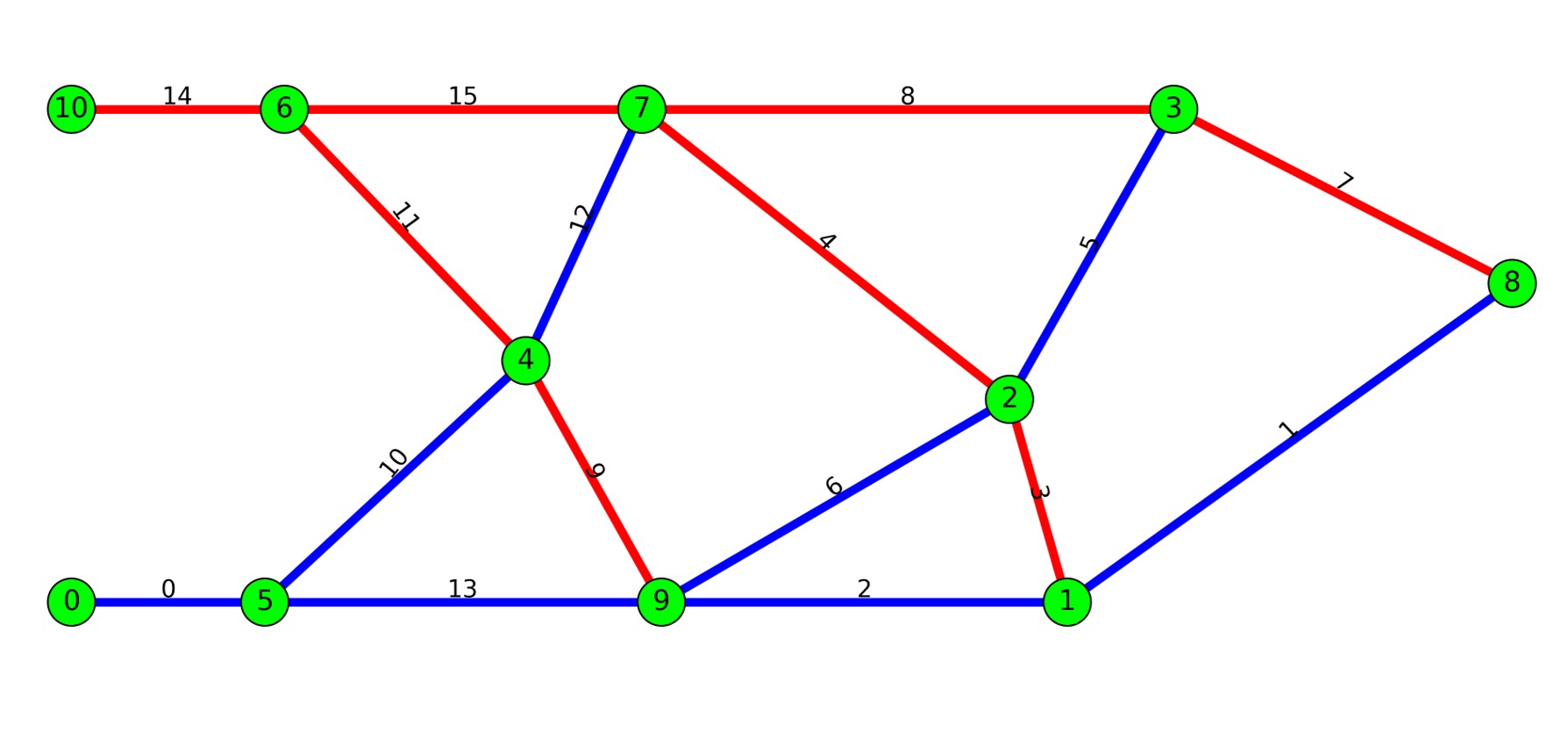}
		\label{fig:lay_init_final_a}
	}
	\centering
	\subfloat[]{
		\includegraphics[width=0.475\textwidth]{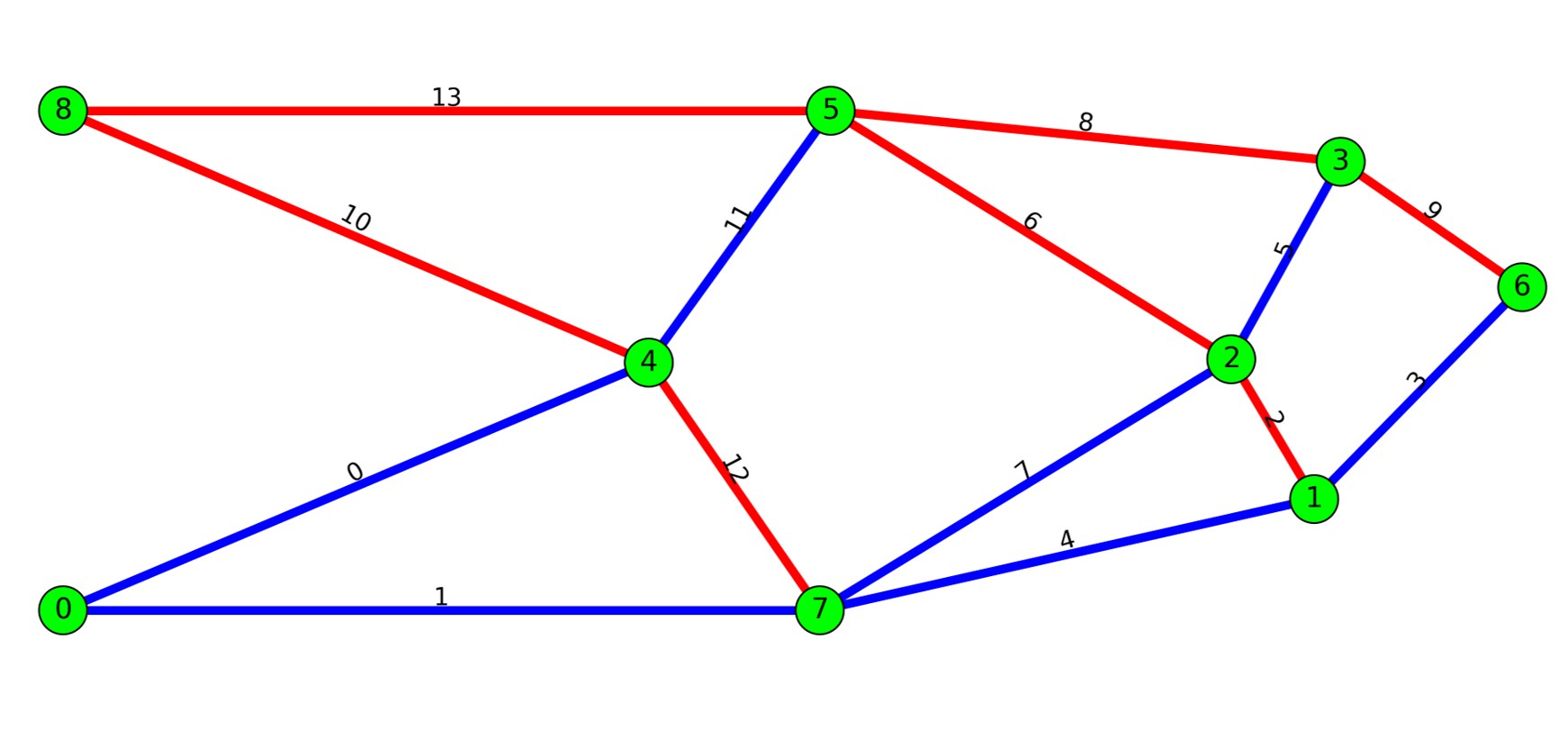}
		\label{fig:lay_init_final_b} 
        }
        \\[-6pt]
	\centering
	\caption{Final layout (b) obtained for sequential size and layout optimization of initial frame model (a).}
	\label{fig:lay_init_final}
\end{figure}
During layout optimization, the length of certain members may reduce and become as short members. These short members are removed, and their end nodes are merged based on the defined merge ratio as explained in edge contraction (see \ref{subsec:pruning&contraction}). Both size and layout optimization are performed sequentially until the optimum design is obtained. The optimal design is identified when there is minimal change in the objective function within two consecutive iterations. Next, the CAD model generation of the optimal design is explored.

\subsection{CAD Model Generation}\label{subsec:CAD_modelgen}
The graph model with nodes, edges and optimized member cross-sections is sufficient to create a structural frame model. From the structural frame model details, the CAD model can be easily generated. There are several methods to perform this generation \cite{YIN2020113102, smith2016application}, and in this study, solid models are generated using the Constructive Solid Geometry (CSG) tree technique.
\begin{figure}[H] 
	\centering
	\subfloat[]{
		\includegraphics[width=0.49\textwidth]{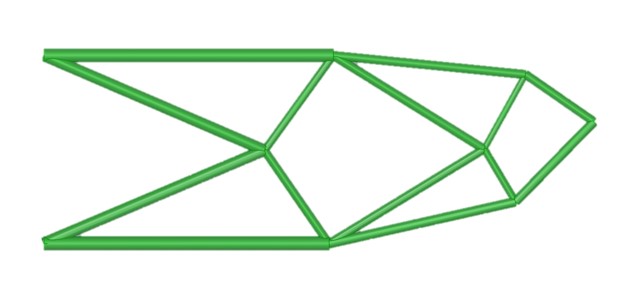}
		\label{fig:Illus_CAD_a}
	}
	\centering
	\subfloat[]{
		\includegraphics[width=0.49\textwidth]{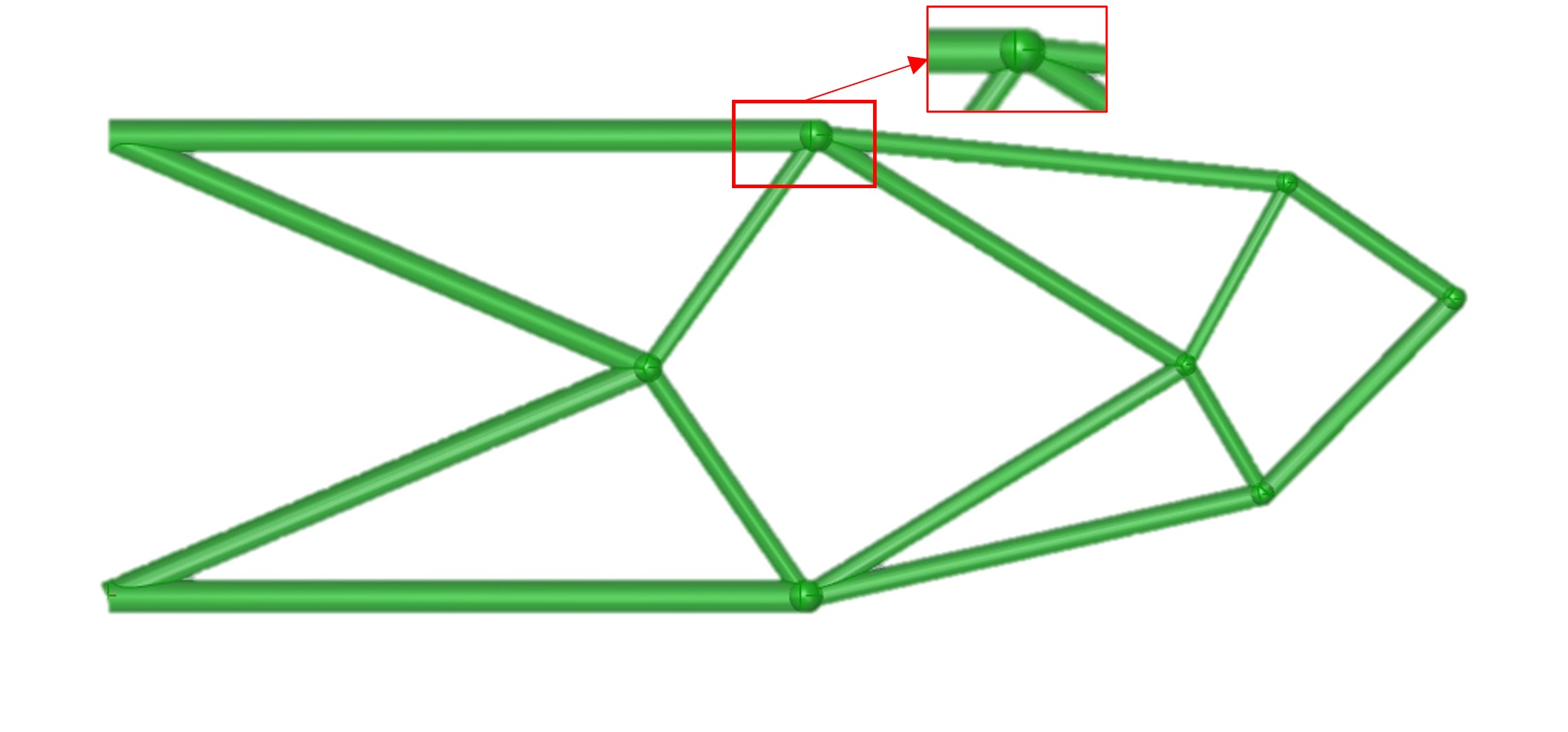}
		\label{fig:Illus_CAD_b}
	} \\[-3pt]
        \centering
	\subfloat[]{
		\includegraphics[width=0.49\textwidth]{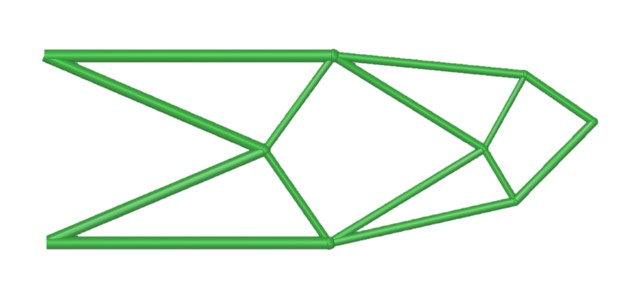}
		\label{fig:Illus_CAD_c}
	}
	\centering
	\subfloat[]{
		\includegraphics[width=0.49\textwidth]{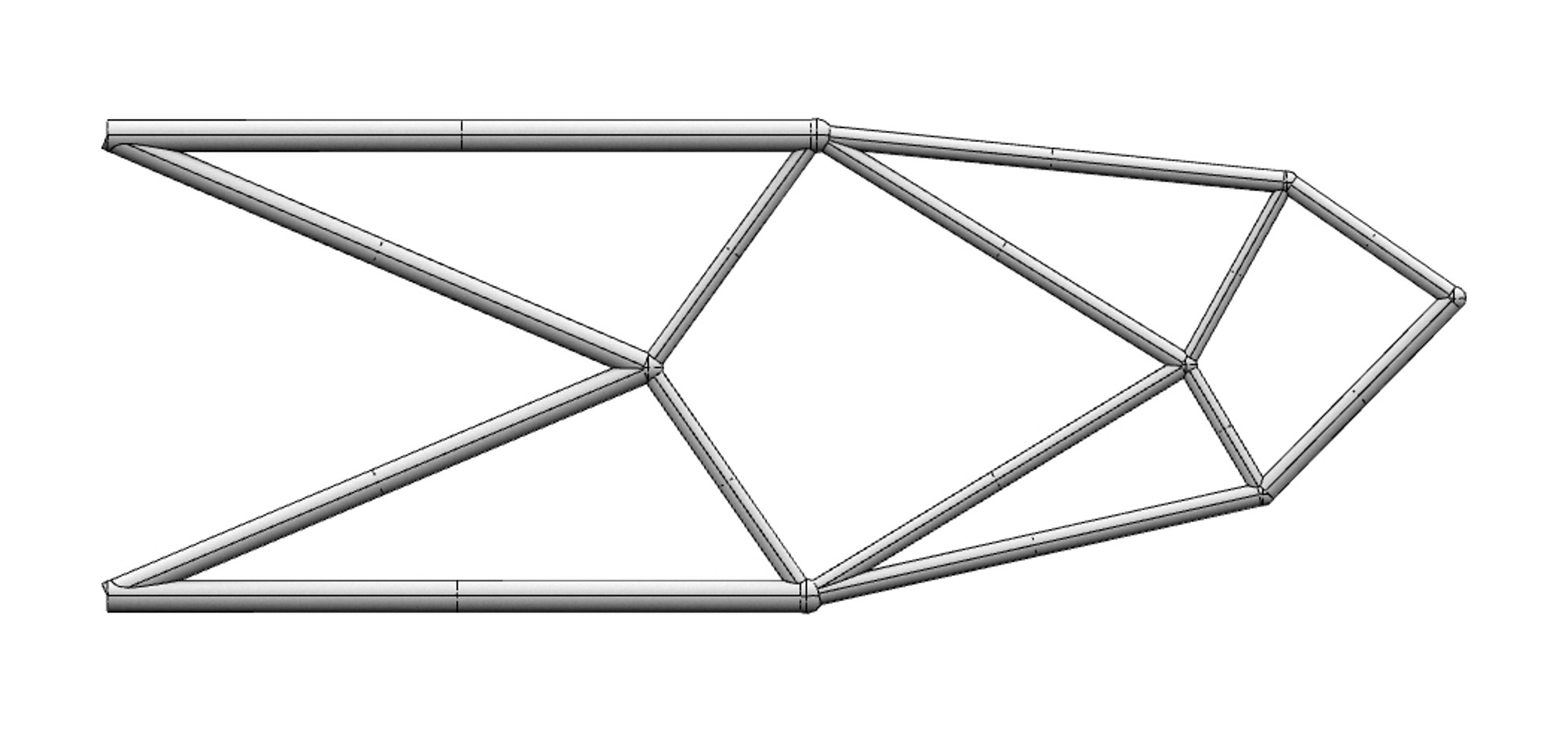}
		\label{fig:Illus_CAD_d}
	} \\[-3pt]
	\centering
	\caption{Illustration of CAD model generation method. (a) Extruded cross sections along the center line of members, (b) Spheres generation at the joints of the structure, (c) Solid model created using Boolean operations, (d) CAD model generated using CSG}
	\label{fig:Illus_CAD}
\end{figure}
First, each member is created by extruding a selected cross-section along the length of the line element, with the cross-sectional area determined from size optimization. In the case of additive manufacturing products, at each joint, a sphere can be created with a maximum radius of all cylinders connected to the joint. After this, a single solid model can be created using Boolean operations and a single watertight model is created. 

\newpage

\section{Structural Analysis and Design of the Topology-, Size- and Layout-optimal Model}
\label{sec:design_informed_stropt}
In this Section, the structural integrity of the final topology-, size- and layout-optimal model is assessed to verify the strength, stiffness and stability requirements. Herein, the guidelines given in EN 1993-1-1: Eurocode 3: Design of steel structures \cite{standard2006eurocode} are followed. By following a similar approach it is straightforward to extend the proposed design process to other standard codes of practice.
\subsection{Section Assignment}\label{subsec:assign_sect}
Engineers would greatly benefit from the ability to utilize a wide range of industry-standard sections for manufacturing optimized structures. By determining the optimal cross-sectional areas of the members through size optimization, we can leverage these findings to filter and select the appropriate standardized sections. However, when the cross-sectional shape changes (such as channel, angle, universal beam, universal column, circular hollow and rectangular hollow etc.,), the slenderness of the members varies, leading to additional checks on buckling. To avoid buckling-induced instability effects, it is vital to check the buckling of the members at the design stage. 

\subsection{Structural Analysis}\label{subsec:str_analysis}

The cross-section assigned optimal models are next subjected to a structural analysis to determine the internal forces. The structural idealization of the optimal design yields bending moments, shear forces and axial forces as member internal forces. At this juncture, a Python algorithm is developed using the \textit{anastruct} module for various load combinations. The guidelines given in EN 1990: Eurocode - Basis of structural design \cite{en20021} and EN 1991-1-1: Eurocode 1: Actions on structures \cite{en20022} are followed in calculating critical load combinations for ultimate and serviceability limit states. The analysis results of examples used in this study are presented in Section \ref{sec:applications}.

\subsection{Structural Design}\label{subsec:str_design}

In this study, structural steel design is focused to produce structurally robust designs. However, it is noted that the extension to other materials is straightforward using the proposed workflow. Herein, the EN 1993-1-1 \cite{en20051} is utilized for the structural design of frame elements. All possible failure modes of steel frame members are considered; such as axial tension, axial compression, flexural buckling, bending, shear and excessive deflection. In the proposed scope of this study, steel connections are idealized as welded connections and are not considered for designing purposes. Additionally, local buckling and lateral torsional buckling are not critical for the nature of applications considered in this study. The governing equations of various failure modes are outlined below.

\subsubsection{Tension Members}\label{subsec:tens_mem}

Tension resistance $(N_{t,Rd})$ according to Eurocode 3 is the smaller value between the design plastic resistance of the gross section $(N_{pl,Rd})$ and design the ultimate resistance of the net cross-section $(N_{u,Rd})$ at holes for fasteners. In this study, the produced optimized are considered with rigid connections (such as arc welding) and hence no holes for fasteners are provided.  Thus, the member tensile resistance is the design plastic resistance of the gross section as expressed in equation \eqref{eq: N_t,Rd},

\begin{equation}\label{eq: N_t,Rd}
    N_{t,Rd}=N_{pl,Rd}=\frac{Af_y}{{\gamma }_{mo}\ }=Af_y
\end{equation}

where $A$ is the gross cross-sectional area, $f_y$ is the yield strength. The partial safety ${\gamma }_{mo}$ is considered unity as per the relevant UK National Annexe \cite{en19931}. For axial tension of the members, the minimum required cross-sectional area is constrained as,

\begin{equation}\label{eq: A - N_t,Rd}
    A\ge \frac{N_{t,Ed}}{f_y\ } \ .
\end{equation}
where $N_{t,Ed}$ is the design axial tensile force of the member.
\subsubsection{Compression Members}\label{subsec:comp_mem}
The compression resistance and minimum required cross-sectional area for the non-slender sections are found using equations \eqref{eq: N_c,Rd}, \eqref{eq: A - N_c,Rd}. For convenience, in this study, the slender sections (Class 04) have been omitted in the design of compression and flexural members.

\begin{equation}\label{eq: N_c,Rd}
    N_{c,Rd}=\frac{Af_y}{{\gamma }_{mo}\ }=Af_y
\end{equation}

\begin{equation}\label{eq: A - N_c,Rd}
    A\ge \frac{N_{c,Ed}}{f_y\ } 
\end{equation}

\subsubsection{Flexural Buckling}\label{subsec:buckling}

The buckling check and constraint formulation of tubular members are presented here. Similarly, for the other sections, similar equations can be easily derived. Members of Circular Hollow Sections (CHS) were considered first for the optimized structures because they have higher radii of gyration for a given volume of material than solid cross-sections and consequently higher flexural buckling resistances \cite{yeetal2021}. The buckling capacity of a CHS (Class 1-3) is given in equation \eqref{eq: N_b,Rd} by taking ${\gamma }_{M1}=1.0,$

\begin{equation}\label{eq: N_b,Rd}
    N_{b,Rd}=\chi Af_y \ .
\end{equation}

Here, the reduction factor $\chi $ for flexural buckling is given as
\begin{equation}\label{eq: buck_red_fac}
    \mathrm{\chi = \frac{1}{\phi + \sqrt{{\phi }^{2} - {\overline{\lambda }}^{2}}} \le 1.0}
\end{equation}
with
\begin{equation}\label{eq: buck_phi}
    \mathrm{ \phi =0.5 \left[ 1 + \alpha \left(\overline{\lambda } - 0.2 \right) + {\overline{\lambda }}^{2}\right]}
\end{equation}
where $\alpha $ is an imperfection factor and $\overline{\lambda}$ is the non-dimensional member slenderness defined as,
\begin{equation}\label{eq: slenderness_ratio}
    \overline{\lambda }\mathrm{=}\sqrt{\frac{f_y}{{\sigma }_F}}
\end{equation}
where ${\sigma }_F$ is the elastic flexural buckling stress given by,
\begin{equation}\label{eq: flex_buck_stress}
    {\sigma }_F\mathrm{=}\frac{{\pi }^{\mathrm{2}}EI}{{A\mathrm{(}KL\mathrm{)}}^{\mathrm{2}}}
\end{equation}
where $E$ is the elastic modulus, $I$ is the second moment of area, $A$ is the cross-sectional area and $L$ is the length of the member. The effective length factor $K$ for fixed-ended members is considered $\mathrm{0.7}$ since the rigid-frame structure is considered in this study. The value of the imperfection factor $\alpha $ can be obtained from Table 6.1 for the different buckling curves identified in Table 6.2 of the Eurocode 3 \cite{en20051} for the tubular members. 

The critical area can be found by equating the average applied axial stress in the member to the local and global buckling stresses. The design variables representing the cross-section namely the diameter $d$ and wall thickness $t$ can be found from this criterion.

The elastic local buckling Stress ${\sigma }_L$ of a cross-section can be considered as,
\begin{equation}\label{eq: loc_buck_stress}
    {\sigma }_L\mathrm{=}\frac{\mathrm{2}E}{\sqrt{\mathrm{3}\left(\mathrm{1-}{\nu }^{\mathrm{2}}\right)}}\frac{t}{d}\mathrm{=1.21}E\frac{t}{d}
\end{equation}

where Poisson's ratio $\nu \mathrm{=}0.\mathrm{3}$. Similarly, the average axial stress ${\sigma }_s$ of a CHS from an applied axial force, $N$ can be approximately represented by,

\begin{equation}\label{eq: avg_axial_stress}
    {\sigma }_s\mathrm{=}\frac{N}{\pi dt}
\end{equation}

The critical values of the geometry can be found by equating the applied axial stress ${\sigma }_s$, the local buckling stress ${\sigma }_L$ and the flexural buckling stress ${\sigma }_F$,

\begin{equation}\label{eq: stress_eqn_relation}
    \frac{N}{\pi dt} =1.21E\frac{t}{d} = \frac{{\pi }^{2}Ed^{2}}{{8(0.7L)}^{2}}
\end{equation}

Therefore, the critical design variables, diameter d and wall thickness t can be determined from,
\begin{equation}\label{eq: thickness_cons}
    t \ge 0.513{\left(\frac{N}{E}\right)}^{1/2} 
\end{equation}
and
\begin{equation}\label{eq: diameter_stress}
    d \ge 0.627 {\left(\frac{NL^{4}}{E}\right)}^{1/6}
\end{equation}
Also, to avoid excessively slender members for practical considerations and serviceability purposes, the following slenderness limitation is also imposed as per Ye et al. \cite{yeetal2021}.
\begin{equation}\label{eq: slenderness_ratio}
    \lambda = \frac{L}{i}\ \ge 150
\end{equation}
\subsubsection{Bending}\label{subsec:towXsec}
The bending resistance $(M_{c,Rd})$ according to Eurocode 3 varies with the structural class of the member. Since the slender members (class 4) are avoided in this workflow, only the first three classes are considered. The bending resistance $(M_{c,Rd})$ for class 1 and 2 cross-sections is given by,
\begin{equation}\label{eq: M_c,Rd_class-1,2}
    M_{c,Rd}=M_{pl,Rd}=\frac{w_{pl}f_y}{{\gamma }_{mo}}
\end{equation}
and for class 3 cross-sections,
\begin{equation}\label{eq: M_c,Rd_class-3}
    M_{c,Rd}=M_{el,Rd}=\frac{w_{el,min}f_y}{{\gamma }_{mo}}
\end{equation}
where  $w_{pl}$ is the plastic section modulus of the section and $w_{el,min}$ is the elastic section modulus of the section. 
\subsubsection{Shear}\label{subsec:shear}
The shear resistance $\mathrm{(}V_{c,Rd}\mathrm{)}$ is discussed in clause 6.2.6 of Eurocode 3. The shear resistance can be taken as $V_{pl,Rd}$ based on a plastic distribution or $V_{el,Rd}$ based on an elastic distribution. In this study, unusual sections are avoided and hence torsion effects are neglected. Hence the shear resistance is considered based on the plastic distribution. The shear resistance $(V_{c,Rd})$ in the absence of torsion is given by,
\begin{equation}\label{eq: shear_res}
    V_{c,Rd} = V_{pl,Rd} = \frac{A_v\left({f_y} /{\sqrt{3}}\right)}{{\gamma }_{mo}}
\end{equation}
where $A_v$ is the shear area. For CHS, the shear area $A_v$ is given as $2A/\pi$ \cite{en20051}.
\subsubsection{Bending and Shear}\label{subsec:bending&shear}
According to clause 6.2.8 of Eurocode 3 \cite{en20051}, when the design shear force exceeds 50\% of the shear resistance, a reduced moment resistance is obtained by reducing the yield strength. This reduction takes the form,
\begin{equation}\label{eq: yield_str_red}
    f_{y,new} = (1-\rho )f_y \ \ \text{where} \ \ \rho = {\left(2\frac{V_{Ed}}{V_{pl,Rd}}-1\right)}^2
\end{equation}
where $V_{Ed}$ is the design shear force.
\subsubsection{Combined Axial and Bending}\label{subsec:comp_axial&bend}
Frame action (provision of rigid joints) of a structure can induce bending combined with axial forces. Hence it is necessary to check for combined effects in these structures. Here, the members are assumed to be subjected to bending about one axis only, so that equations are further simplified. For combined effects, the linear summation of the utilization ratios for each resistance should be less than unity.
\subsubsection{Excessive Deflections}\label{subsec:excess_def}
The deflection check is carried out under serviceability conditions considering the maximum deflection of the whole structure. Maximum deflection due to unfactored loading is usually compared with the limits specified by the designer. For these limits, since there are no specific guidelines given in EN 1993-1 \cite{en20051}, the UK national Annex (NA 2.23) \cite{en19931} recommended values are used. Generally, a limit of span/180 is suggested for cantilevers. For beams, the limit is span/200 and it is reduced to span/360 when it is carrying brittle finishes. More stringent limits can be used based on the requirement of the structure or structural element.

\begin{equation}\label{eq: sum_utilization_ratio - 1}
    \frac{N_{Ed}}{N_{Rd}}+\ \frac{M_{Ed}}{M_{c,Rd}} \le 1
\end{equation}

\begin{equation}\label{eq: sum_utilization_ratio - 2}
    \frac{N_{Ed}}{N_{b,Rd}}+\ k_{yy}\frac{M_{Ed}}{M_{b,Rd}} \le 1
\end{equation}

\begin{equation}\label{eq: sum_utilization_ratio - 3}
    \frac{N_{Ed}}{N_{b,Rd}}+\ k_{zy}\frac{M_{Ed}}{M_{b,Rd}} \le 1
\end{equation}
where $N_{Ed}$ is the design axial force, $M_{Ed}$ is the design bending moment, $N_{Rd}$ is axial resistance, $M_{c,Rd}$ is bending resistance, $N_{b,Rd}$ is the buckling resistance and $M_{b,Rd}$ is the design buckling resistance moment. $k_{yy}$ and $k_{zy}$ are interaction factors that are given in EN 1993-1-1 Annex A and B. It is noted that the check using equation \eqref{eq: sum_utilization_ratio - 1} is sufficient for a member when the axial action is tension. But equations \eqref{eq: sum_utilization_ratio - 2} and \eqref{eq: sum_utilization_ratio - 3} should be considered when the member is subjected to axial compression.

\newpage
\section{Applications}
\label{sec:applications}
This section provides numerical results of the proposed approach for selected examples. First in Section \ref{sec:cant_plate}, a cantilevered plate example is used to illustrate the proposed workflow in detail. Then a simply supported beam example is presented in Section \ref{sec:simp_beam}.

\subsection{A cantilever subjected to a point load}\label{sec:cant_plate}
\begin{figure}[H]
\centering
\includegraphics[width=0.99\textwidth]{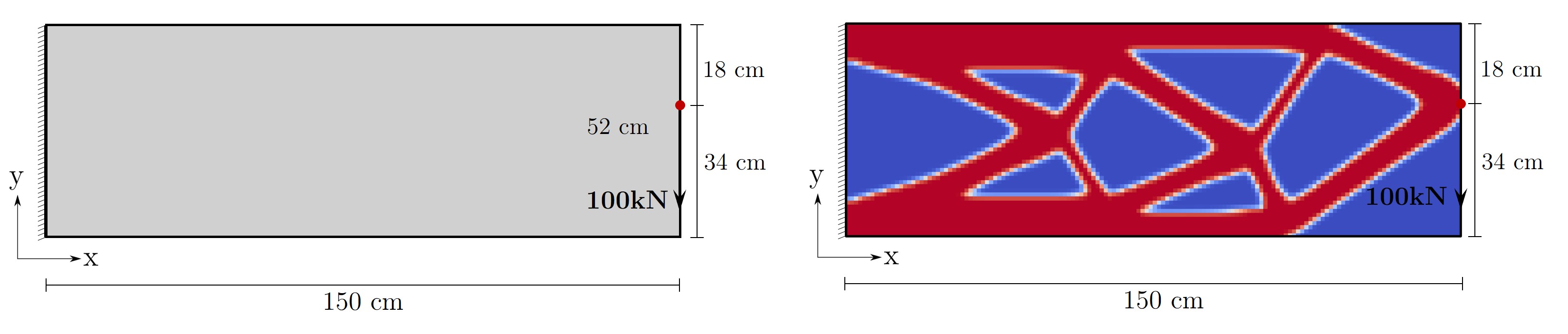}
\caption{Topology optimization of a cantilever plate example ($\overline{E}$ = $2.1 \times 10^5$ N/mm$^2$, ${E}_{min}$ = $1 \times 10^{-9}$ N/mm$^2$, $\nu$ = $0.3$, $V_f$ = $0.3$, $p$ = $3$, $R$ = $1.2$, Maximum iterations = $200$)}
\label{fig:topopt_cantplate}
\end{figure}

In the analysis, a 150 cm $\times$ 52 cm plate is chosen, as illustrated in Figure \ref{fig:topopt_cantplate}. One face of the plate is fixed, and a point load with a value of F = 100 kN is applied at $\overline{y}$ = 18 cm from the top of the plate on the opposite free end.  For the size optimization process, the plate is assumed to have a uniform thickness of 1 cm.

\begin{figure}[H]
\centering
\includegraphics[width=0.99\textwidth]{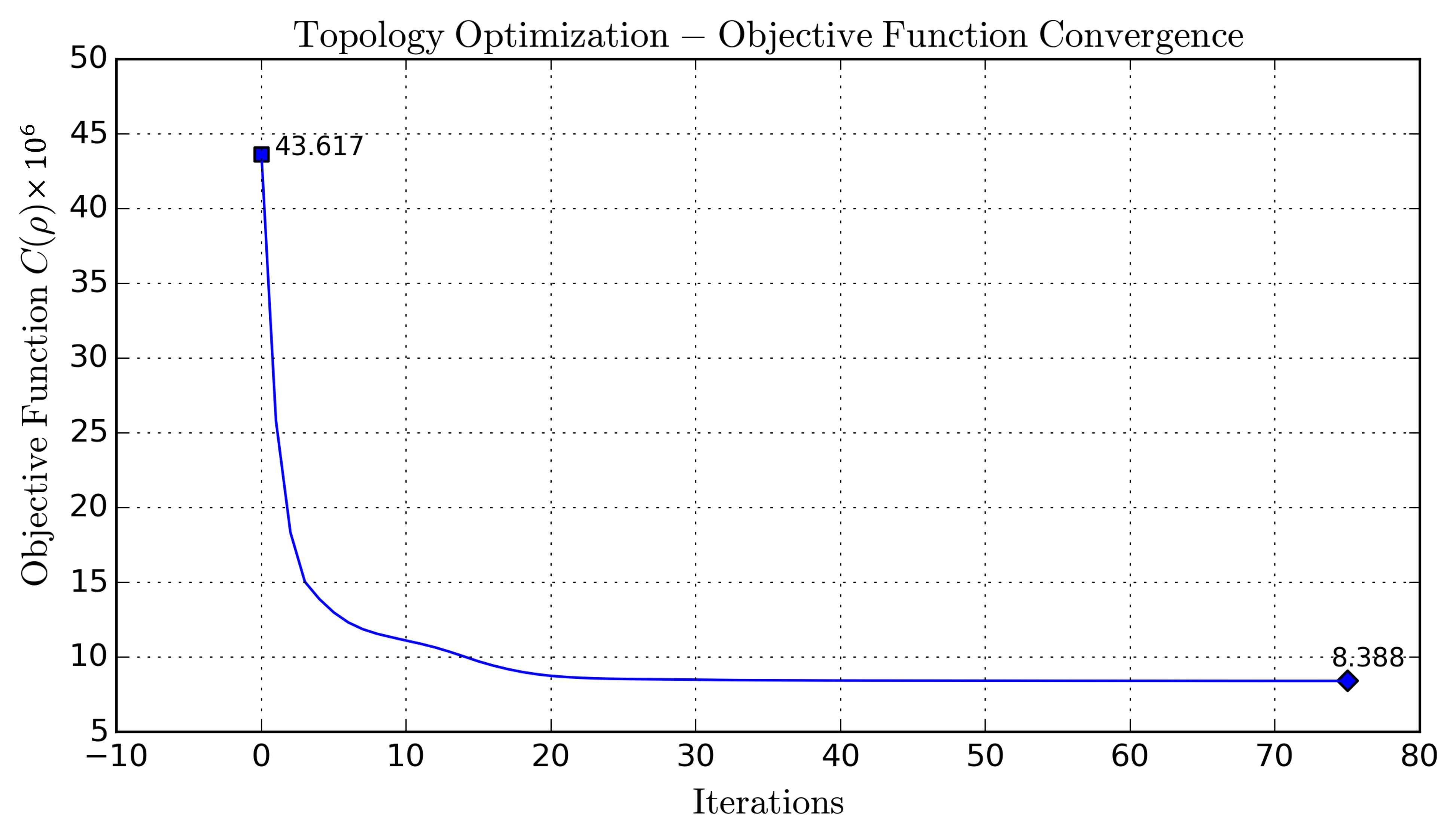}
\caption{Objective function convergence during topology optimization of cantilever example}
\label{fig:cant_obj_con}
\end{figure}

In this example, topology optimization was performed on a discretized domain consisting of 150 $\mathrm{\times}$ 52 $\mathrm{\times}$ 1 linear quadrilateral elements. The objective function $C_{(\rho)}$ was minimized following the approach detailed in Section \ref{subsec:TopOpt_Review}. The cantilever example employed specific parameters: a volume fraction ($V_f$) of 0.5, penalization power ($p$) set to 3, a filter radius ($R$) of 1.2, and a maximum of 200 iterations. The convergence of the objective function was analyzed and illustrated in Figure \ref{fig:cant_obj_con}.

The 2D binary image data obtained from the topology optimization is skeletonized and the frame model is extracted as discussed in Sections \ref{subsec:Skeletonization} and \ref{subsec:frame_skel}. Fixed-end boundary pixels and force-applied pixels were tagged so no alteration is allowed during any stages of optimization. Sequential size and layout optimization are performed for the extracted frame model by taking the thickness of the plate as 1 cm. First, sequential optimization starts with size optimization.  The lower bound of the areas is taken as circular areas with a radius of 0.5 cm and the upper bound radius is set to 10 cm. The initial sizes of the members are assumed to have the same cross-sectional area of $A=7.09\ cm^2$ which is found by dividing the volume $V=0.5\overline{V}=3900\ cm^3$ by the total length of the members. After the first size optimization iteration, the optimized cross sections for the initial layout are known. Next, the layout optimization starts to find the optimal layout. The lower and upper bounds can be given as either the whole design domain or small allowable regions for each node. During layout optimization, some members can become short, and it can be removed by the edge contraction as explained in Section \ref{subsec:pruning&contraction}. In this example, a merge ratio of 0.1 is used to remove short members after each layout optimization. The parameters used for both size and layout optimization are given in Table~\ref{table:3}.

\begin{table}[h!]
    \centering
    \caption{Parameters used for Size and Layout Optimization}
    \label{table:3}   
    \begin{tabular}{p{3.5in} c} \hline 
    \centering \textbf{Parameter} & \textbf{Value} \\ \hline 
    Volume constraint ($V_{(A,s)})$ & $\le 3900\ {cm}^3$ \\ \hline 
    Tolerance (${\epsilon }_{size}\ ,{\epsilon}_{layout}$) & ${10}^{-4}$ \\ \hline 
    Minimum Area ($A_{min}$) & $0.785\ {cm}^2$ \\ \hline 
    Maximum Area ($A_{max}$) & $\mathrm{314.160}\ {cm}^2$ \\ \hline 
    Bounds for Nodes & Whole design domain \\ \hline 
    Maximum Iterations for SQP Algorithm for each Size and Layout Optimization ($i_{max}$) & 200 \\ \hline 
    Maximum Iterations for MMA Algorithm for each Size and Layout Optimization ($i_{max}$) & 20 \\ \hline 
    Merge Ratio $(\zeta )$ & 0.1 \\ \hline 
    \end{tabular}
\end{table}

Sequential optimization stops when it converged according to the defined tolerance value below.
\[\left|\frac{C_{i-1}-C_i}{C_i}\ \right|<{\epsilon }_{frame}\ \]

In this study, tolerance ${\epsilon }_{frame}={10}^{-4}$ is used. Here $C_i$ refers to the objective function value at the end of each size or layout optimization. 

Throughout the size and layout optimization process, the volume of the frame structure is constrained to be equal to a specified value of $V_f\overline{V}=3900\ cm^3$. Both the Sequential Quadratic Programming (SQP) and Method of Moving Asymptotes (MMA) algorithms are utilized in this example, and their results are compared. It is observed that while the SQP algorithm can handle layout optimization even when the entire design domain is defined as a box-bounds constraint with a higher maximum iteration value, the MMA algorithm fails to do so. The MMA algorithm tends to destabilize the structural problem more than the SQP algorithm under these settings. To address this issue, the maximum iteration value can be reduced or box-bound constraints can be specified for individual nodes as smaller box regions.

\begin{figure}[H]
\centering
\includegraphics[width=0.99\textwidth]{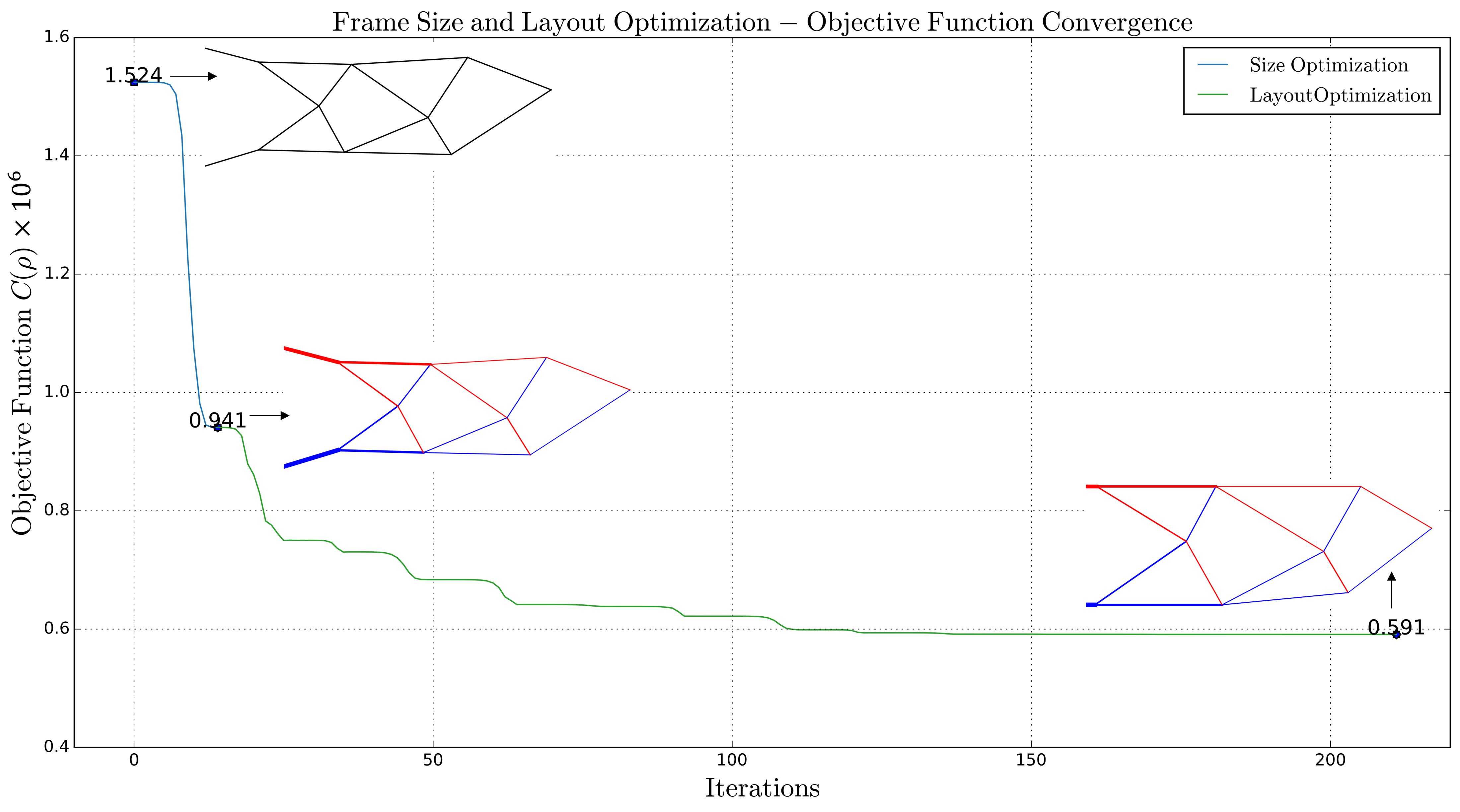}
\caption{Objective Function Convergence for Frame Size and Layout Optimization using SQP Algorithm when merging is not allowed}
\label{fig:cant_sqp}
\end{figure}

\begin{figure}[H]
\centering
\includegraphics[width=0.99\textwidth]{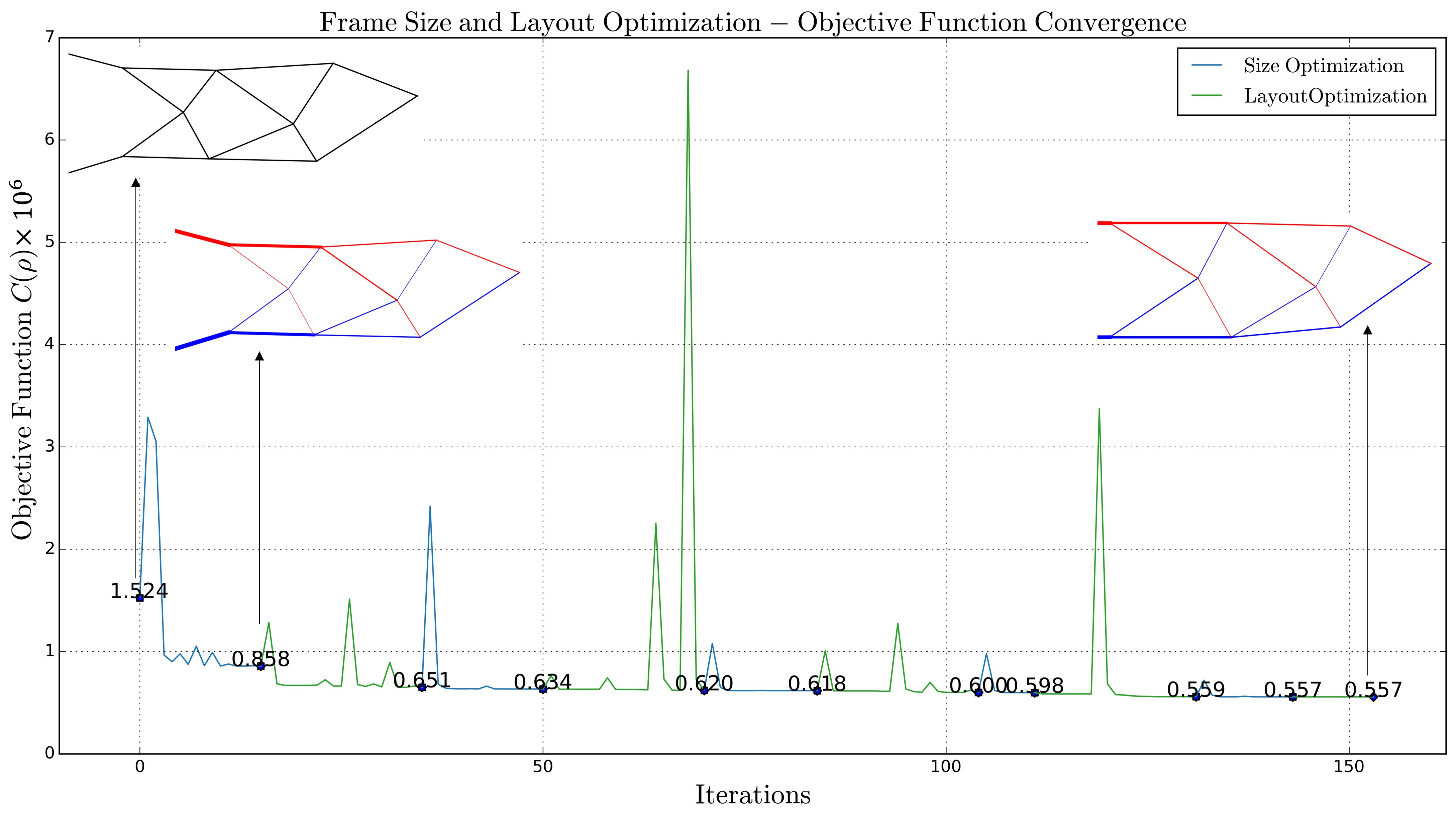}
\caption{Objective Function Convergence for Frame Size and Layout Optimization using MMA Algorithm when merging is not allowed}
\label{fig:cant_mma}
\end{figure}

\begin{figure}[H]
\centering
\includegraphics[width=0.99\textwidth]{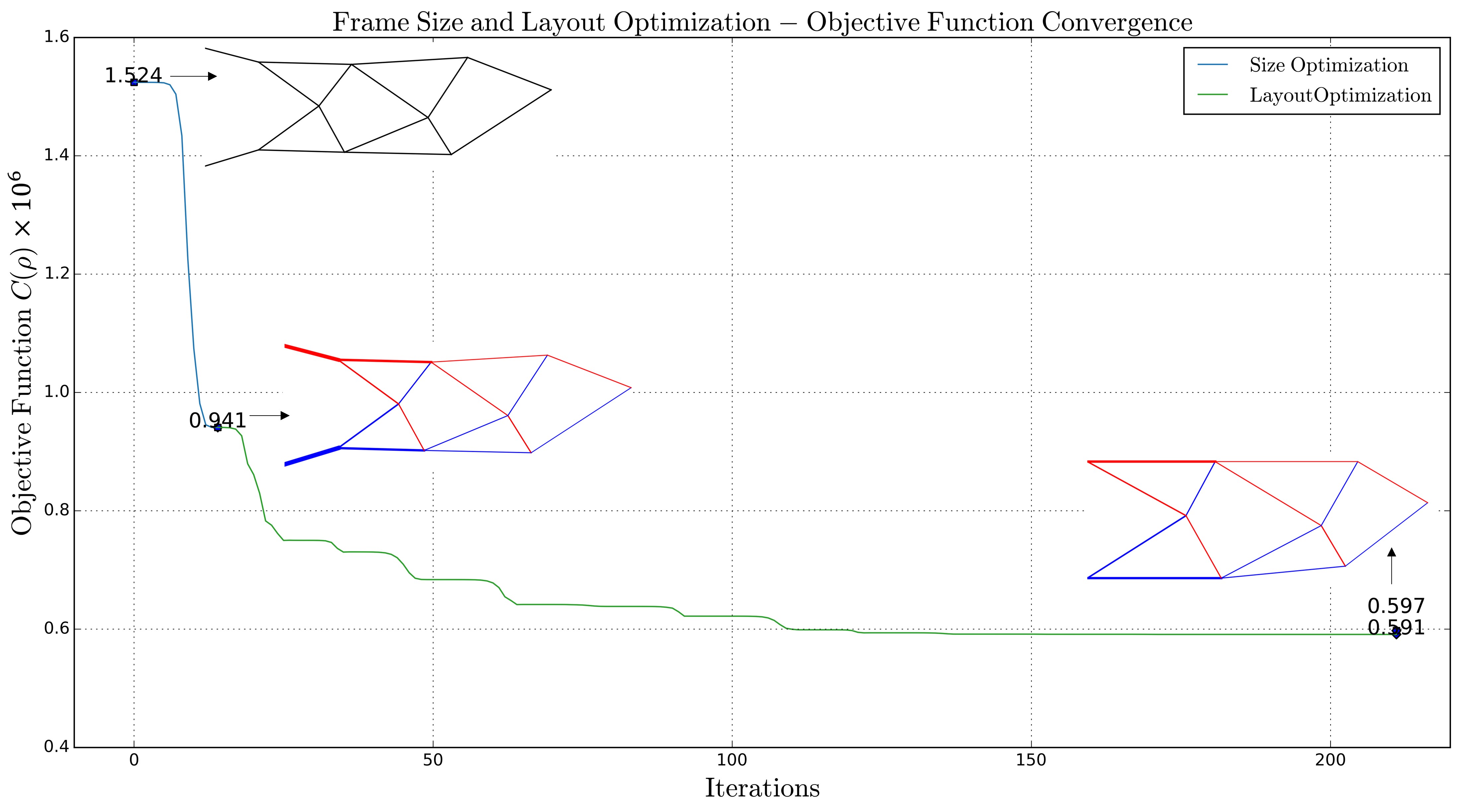}
\caption{Objective Function Convergence for Frame Size and Layout Optimization using SQP Algorithm when merging is allowed}
\label{fig:cant_sqp_merge}
\end{figure}

\begin{figure}[H]
\centering
\includegraphics[width=0.99\textwidth]{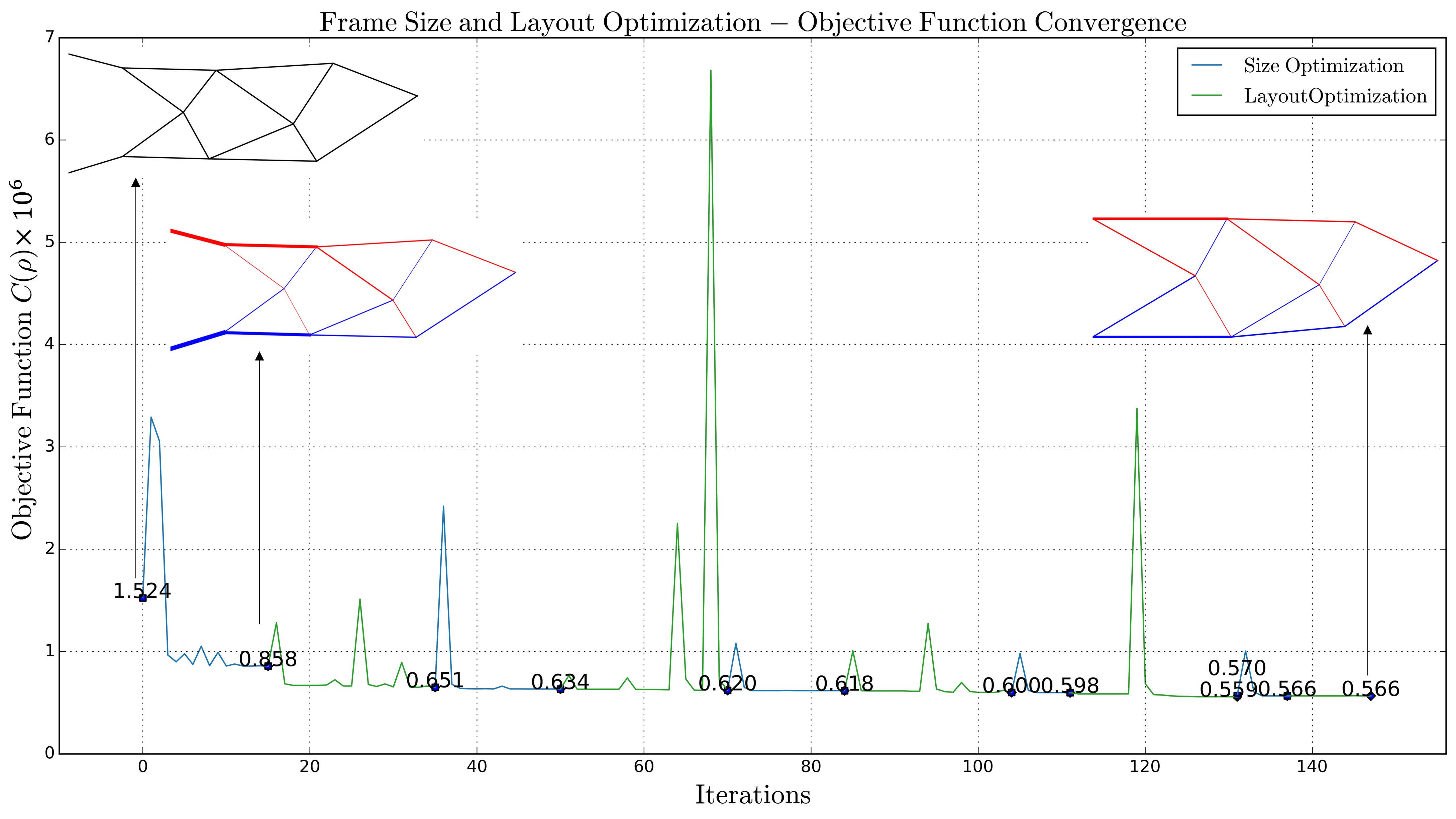}
\caption{Objective Function Convergence for Frame Size and Layout Optimization using MMA Algorithm when merging is allowed}
\label{fig:cant_mma_merge}
\end{figure}

Convergence plots, such as Figures \ref{fig:cant_sqp} and \ref{fig:cant_mma} (when merging is not allowed), as well as Figures \ref{fig:cant_sqp_merge} and \ref{fig:cant_mma_merge} (when merging is allowed), show the behavior of the objective function during optimization. These plots indicate that the MMA algorithm tends to find a more optimum output with fewer total iterations compared to the SQP algorithm. However, it should be noted that the convergence plots for the MMA algorithm exhibit some spikes, while the SQP algorithm demonstrates smoother convergence in Figures \ref{fig:cant_sqp} and \ref{fig:cant_sqp_merge}.

The MMA algorithm's ability to achieve more optimal outputs with fewer iterations suggests its effectiveness in exploring the design space and identifying promising design configurations efficiently. One possible reason for MMA's better performance could be its efficient handling of box-bound constraints, as it has been specifically developed to address such constraints.
\begin{figure}[H] 
	\centering
    \captionsetup{justification=centering}
	\subfloat[Chosen cantilever example]{
		\includegraphics[width=0.49\textwidth]{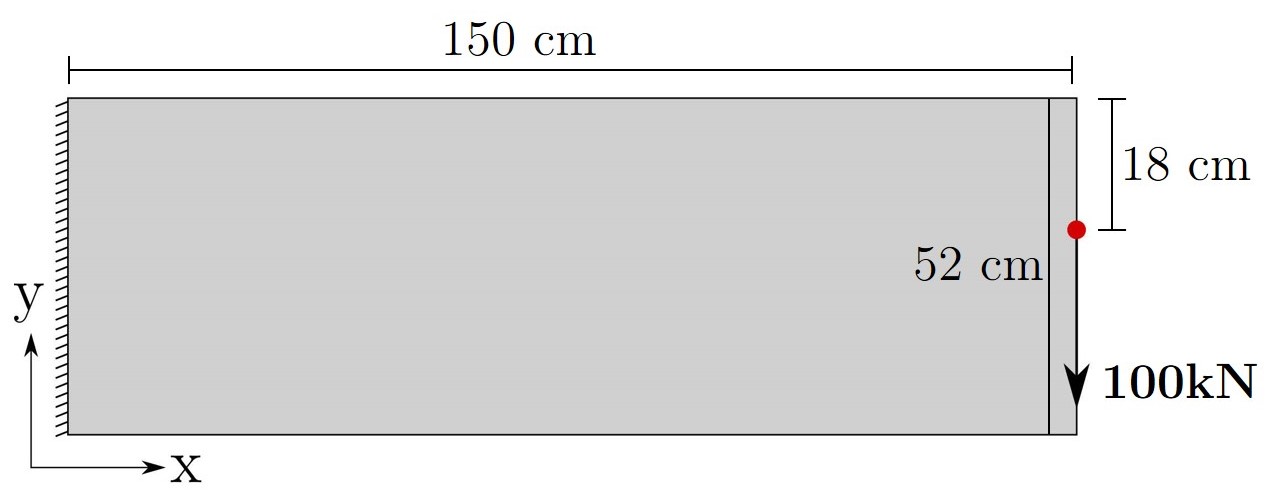}
		\label{fig:illu_cant_opt_a}
	}
	\centering
	\subfloat[Topology optimized model for $V_{f}=0.5$]{
		\includegraphics[width=0.49\textwidth]{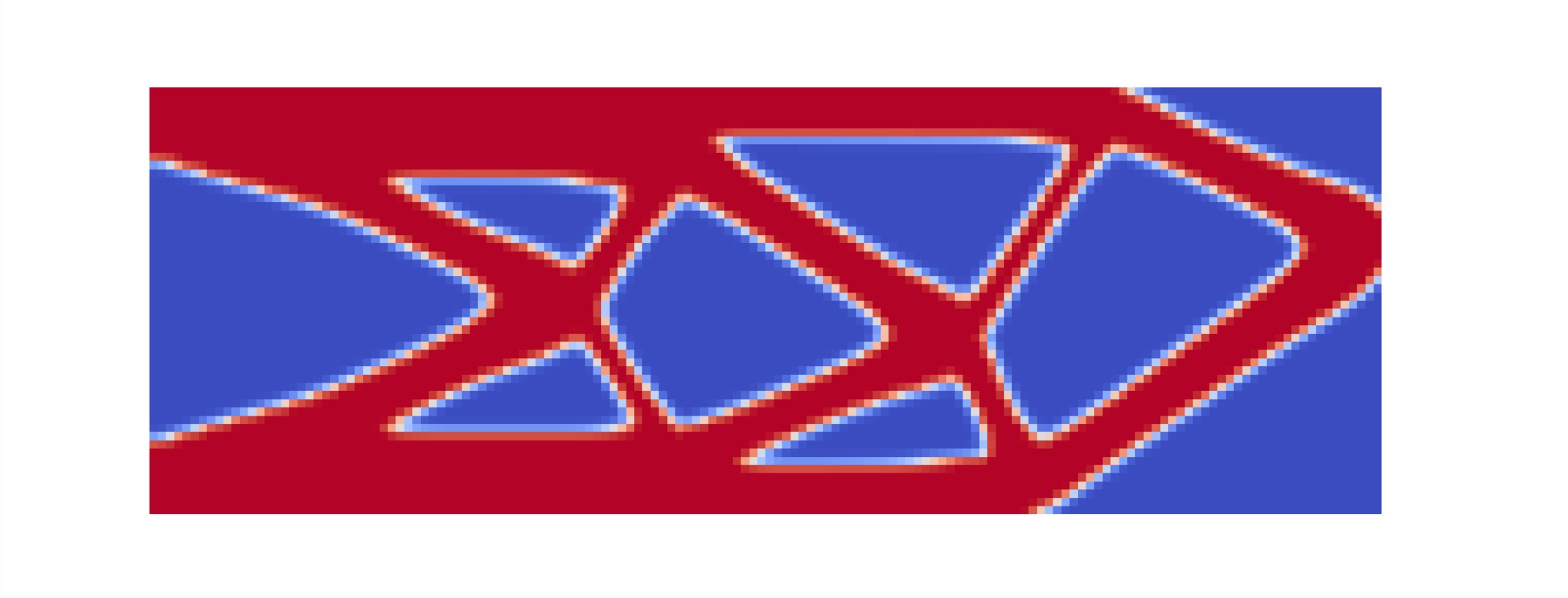}
		\label{fig:illu_cant_opt_b}
	} \\[-9pt]
        \centering
	\subfloat[Binary image generated for corresponding topology optimized model in (b)]{
		\includegraphics[width=0.49\textwidth]{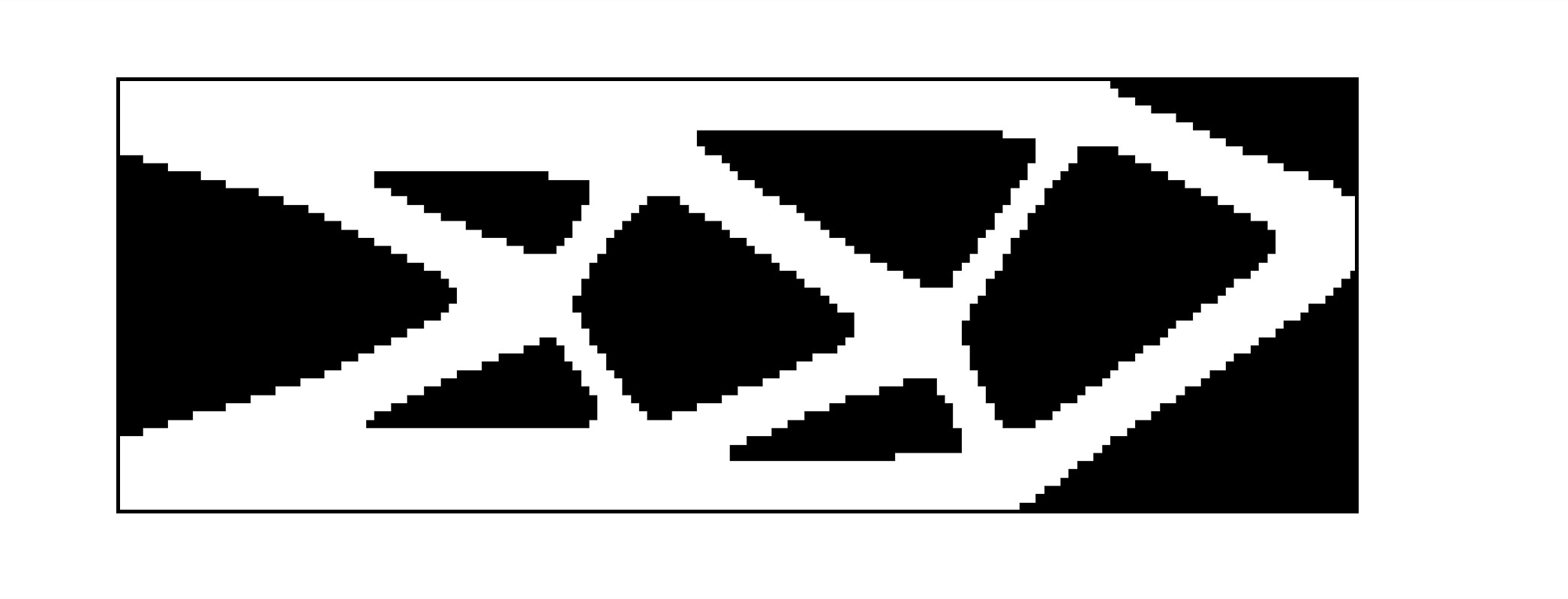}
		\label{fig:illu_cant_opt_c}
	} 
	\centering
	\subfloat[Skeleton obtained after tagging boundary pixels (as left top and bottom pixels) and loaded element pixel]{
		\includegraphics[width=0.49\textwidth]{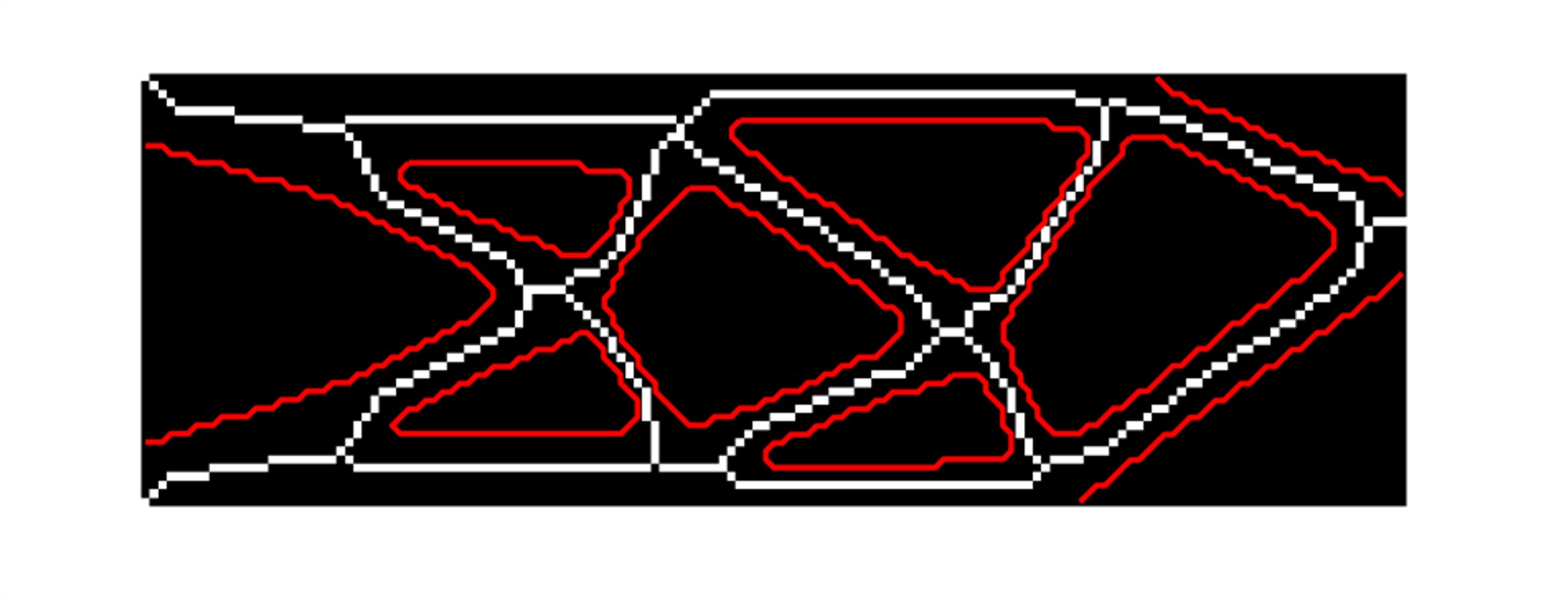}
		\label{fig:illu_cant_opt_d}
        } \\[-9pt]
  	\centering
	\subfloat[Extracted frame model from skeleton by merging short members with $\zeta$ = 0.1]{
		\includegraphics[width=0.49\textwidth]{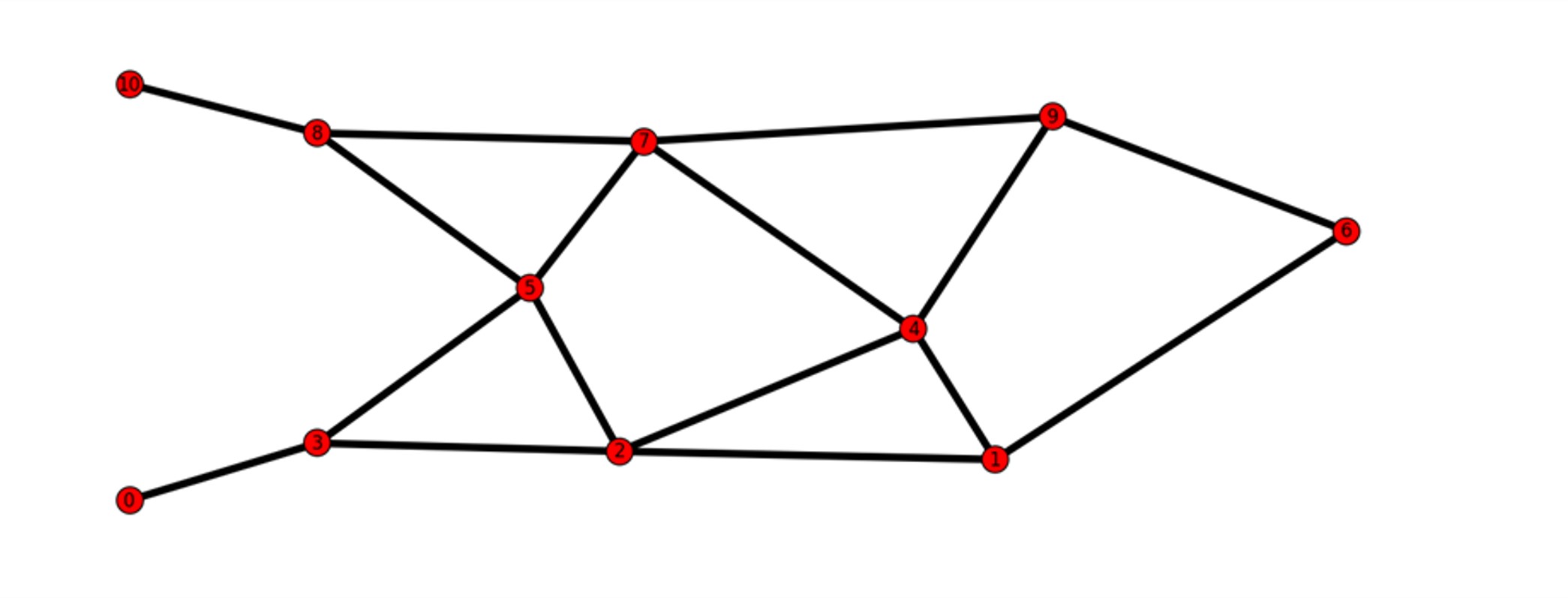}
		\label{fig:illu_cant_opt_e}
	} 
        \centering
	\subfloat[Initial frame model given for sequential size and layout optimization for both SQP and MMA algorithm]{
		\includegraphics[width=0.49\textwidth]{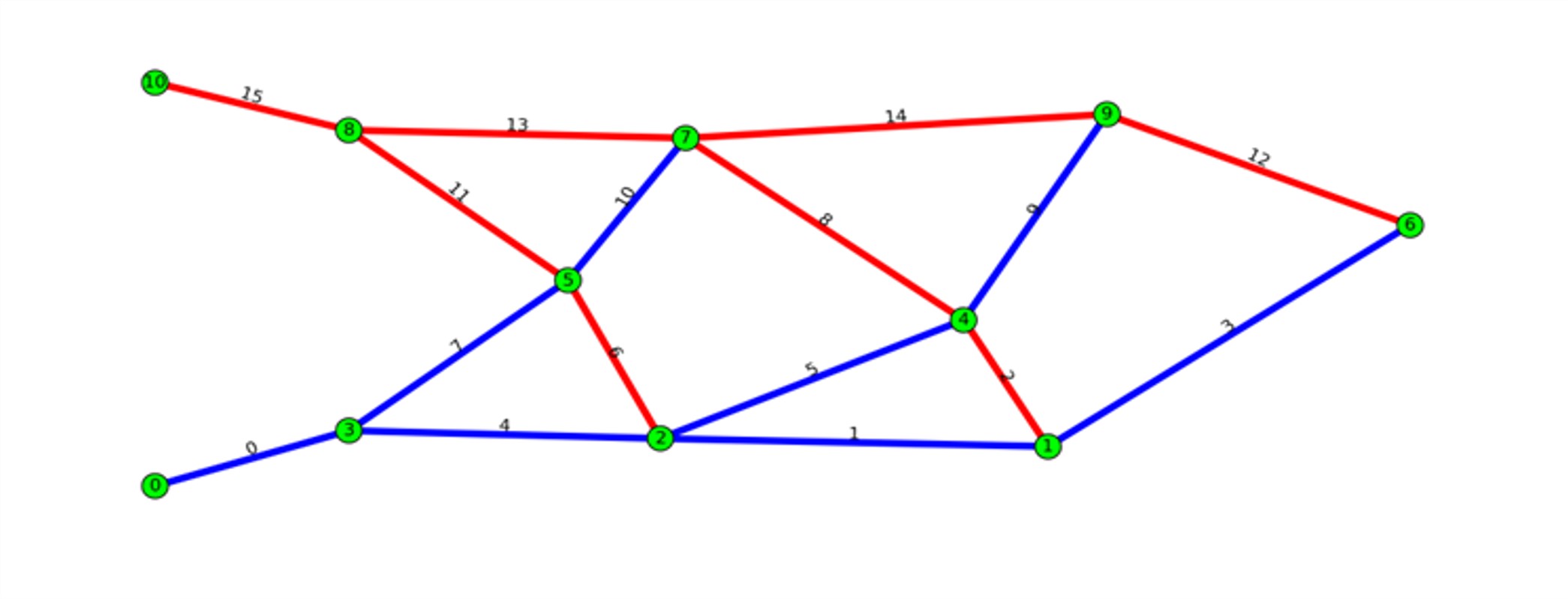}
		\label{fig:illu_cant_opt_f}
	} \\[-9pt]
	\centering
	\subfloat[Final Optimized model for MMA algorithm when merging is not allowed]{
		\includegraphics[width=0.49\textwidth]{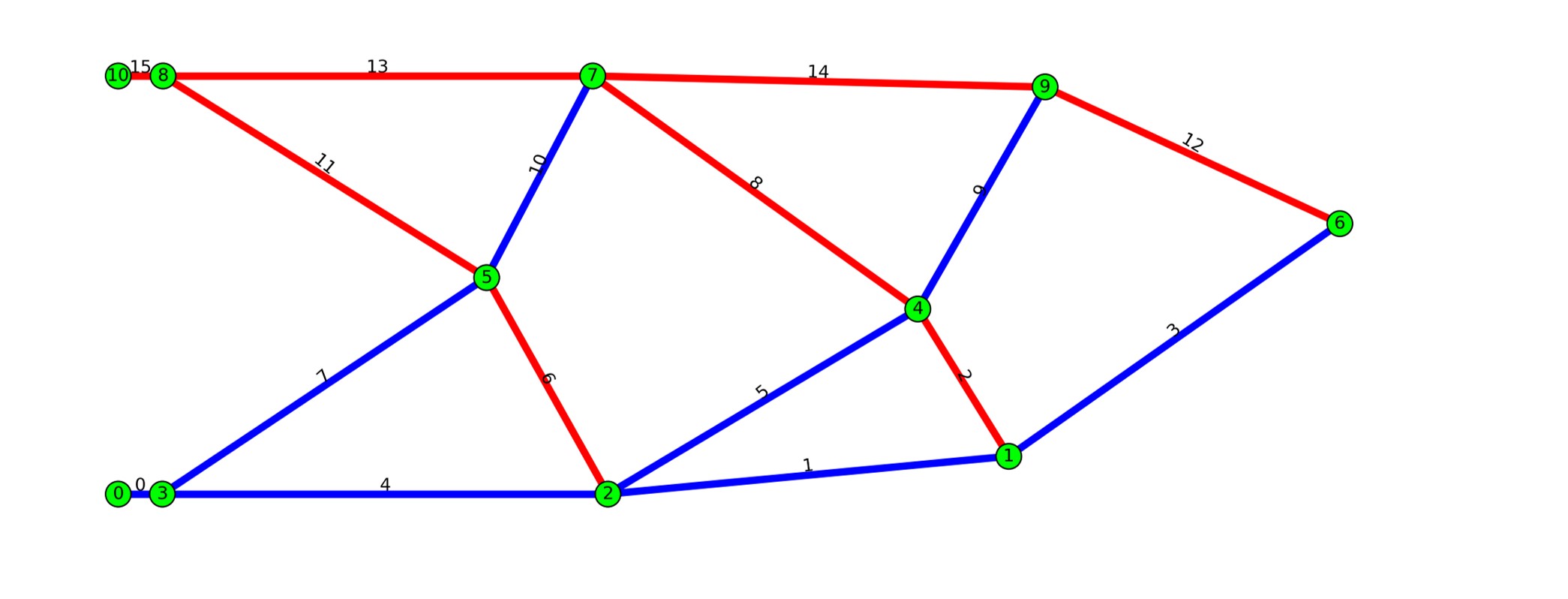}
		\label{fig:illu_cant_opt_g}
        }
  	\centering
	\subfloat[Final Optimized model for MMA algorithm when merging is allowed with $\zeta$ = 0.1]{
		\includegraphics[width=0.49\textwidth]{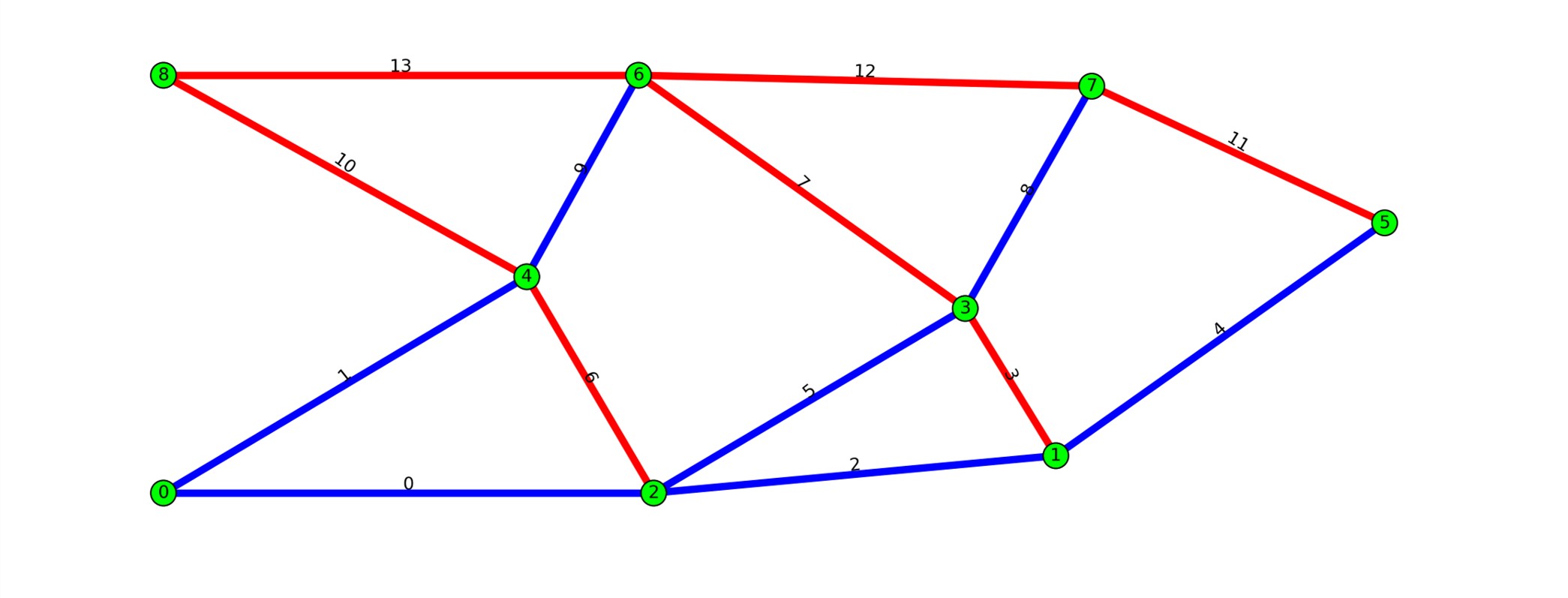}
		\label{fig:illu_cant_opt_h}
	} \\[-9pt]
        \centering
	\subfloat[Final Optimized model for SQP algorithm when merging is not allowed]{
		\includegraphics[width=0.49\textwidth]{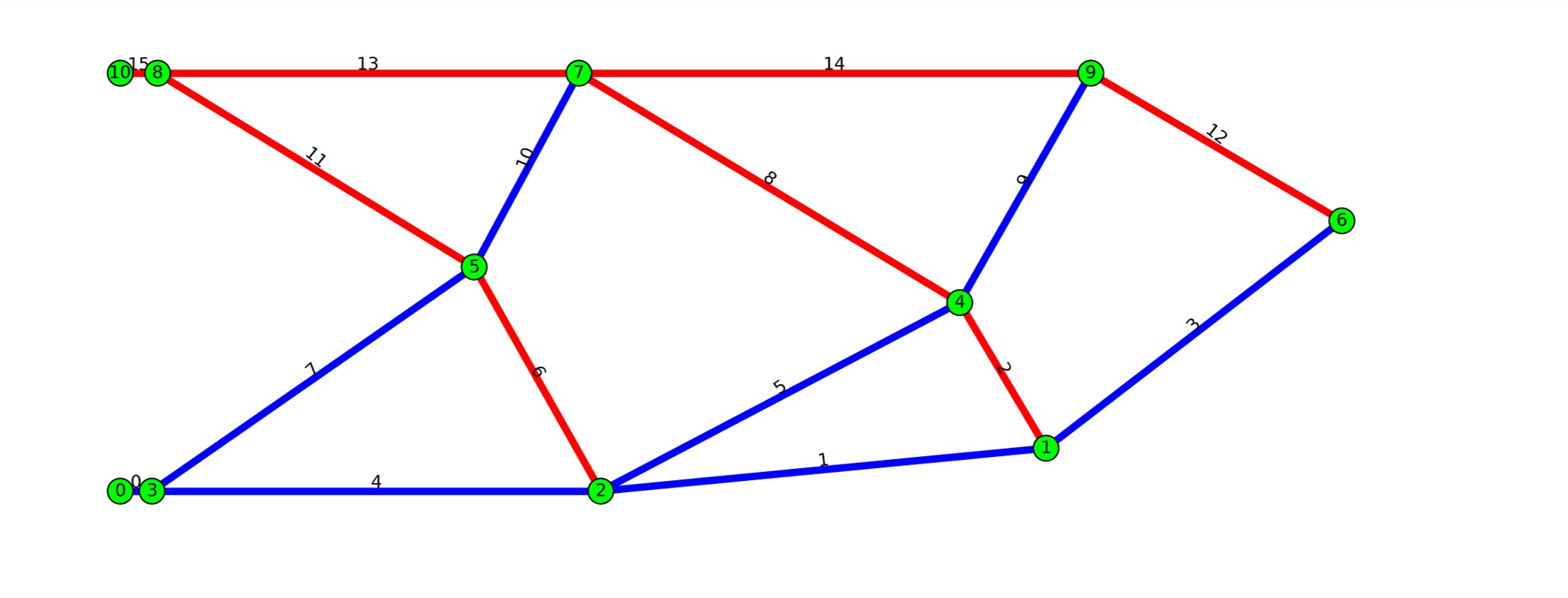}
		\label{fig:illu_cant_opt_i}
	}
	\centering
	\subfloat[Final Optimized model for SQP algorithm when merging is allowed with $\zeta$ = 0.1]{
		\includegraphics[width=0.49\textwidth]{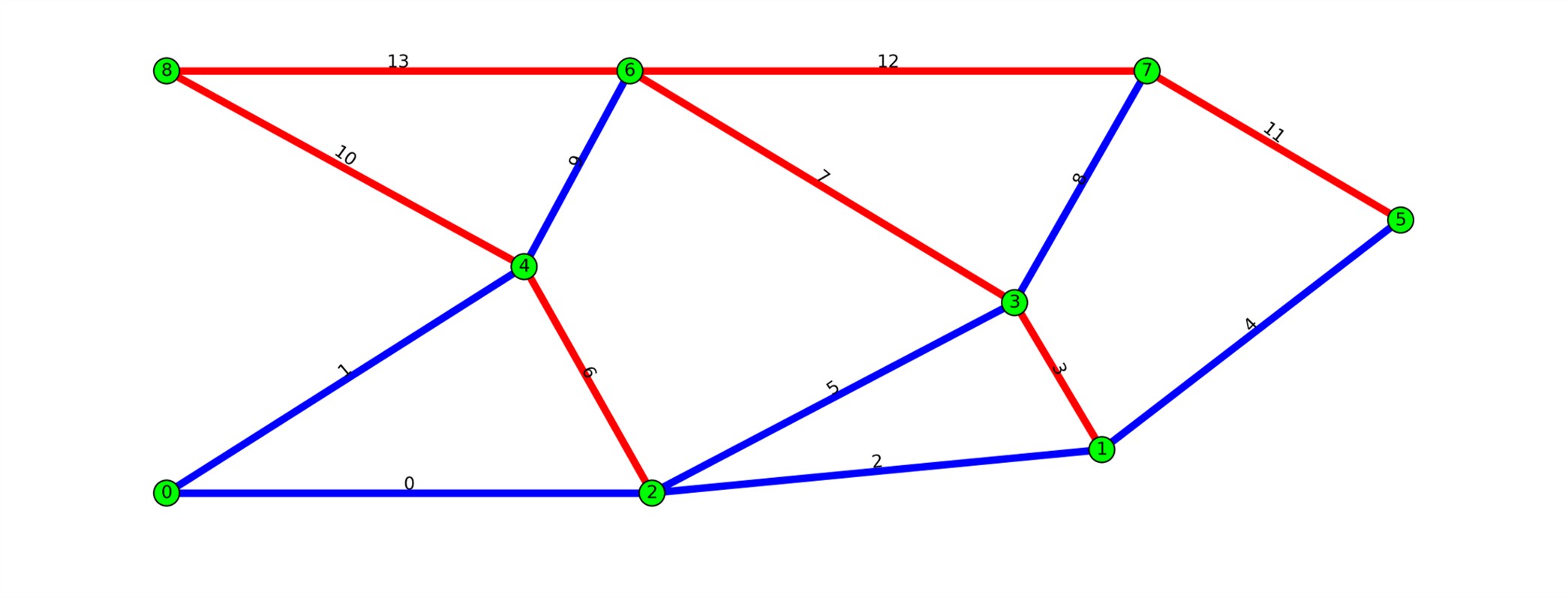}
		\label{fig:illu_cant_opt_j}
	} \\[-6pt] 
	\centering
	\caption{Illustration of Cantilever Optimization}
	\label{fig:illu_cant_opt_0.5}
\end{figure}

For this example, the initial output is given without using the angle constraint for small-angled members and observed the same optimized structure as the output. The initial and final structures of size and layout optimization for both algorithms are shown in Figures \ref{fig:illu_cant_opt_f} and \ref{fig:illu_cant_opt_j} respectively. Further, the final objective function value of both SQP and MMA algorithms is higher for the final structures when the merging is allowed with a merge ratio of 0.1. However, considering the potential manufacturing challenges, the elimination of small members is preferred over a slight compromise in the objective function. Therefore, for subsequent analysis and design, the merging is allowed (at a merge ratio of 0.1) and MMA algorithm is favoured.

After extracting the frame which is topology-, size- and layout- optimal, the next tasks include the structural analysis and design of the same. Initially, the possibility of using Square Hollow Sections (SHS) of thickness 6.3 mm and Circular Hollow Sections (CHS) of thickness 4 mm as alternatives were analyzed. The notation CHS mxn signifies that m represents the outer diameter of the circular section, while n indicates the thickness of the walls.

The in-plane dimensions of the cross-section are selected for each member based on the size optimization results and commercial availability of the section sizes. In some cases, the size optimization may yield a specific value for the required area of the section, but it might not be possible to find an exact section size with that exact area in commercial options. To address this issue, an algorithm can be employed to search for available section sizes that are equal to or larger than the desired area value obtained from the size optimization. The algorithm aims to find the closest section size with an area equal to or greater than the desired value. This ensures that a suitable section is selected for the given optimization result, even if an exact match is not available in commercial options. This approach ensures that the selected section is feasible and can be readily sourced from commercial suppliers as well. For example, member 1 requires an area of 587 $mm^2$ from the size optimization, but no section with that size is commercially available. Instead, the closest suitable option, CHS 60.3x4 with an area of 707 $mm^2$, is chosen considering the 4mm thickness constraint for all sections. Similarly, one can utilize any type of section in the structural analysis and design based on the standard codes of practice.

Herein, first, the optimal frame with assigned standard sections is structurally analyzed using the \textit{anastruct} module, and the critical member forces are found for each member. Then, the frame is checked according to Eurocode 3 guidelines for those critical forces as mentioned in Section \ref{sec:design_informed_stropt}. The analysis results and the utilization ratios according to Eurocode 3 are presented below in Figure \ref{fig:cant_strana} and Table~\ref{table:5}. The sign convention used here is as follows; axial tension -- positive, axial compression -- negative, hogging bending moment -- positive, sagging bending moment -- negative, anticlockwise shear -- positive, clockwise shear -- negative and vertical downward deflection -- positive.

\begin{figure}[h!] 
	\centering
	\subfloat[Structure]{
		\includegraphics[width=0.49\textwidth]{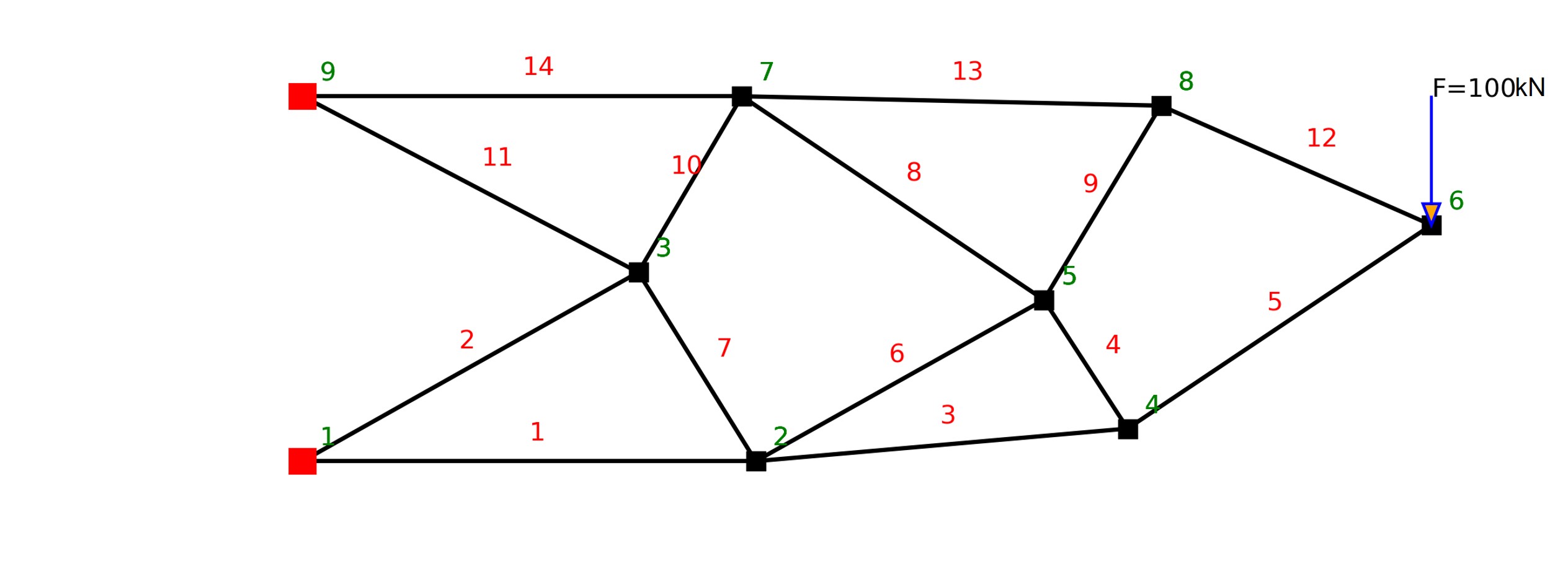}
		\label{fig:cant_strana_a}
	}
	\centering
	\subfloat[Reaction Forces]{
		\includegraphics[width=0.49\textwidth]{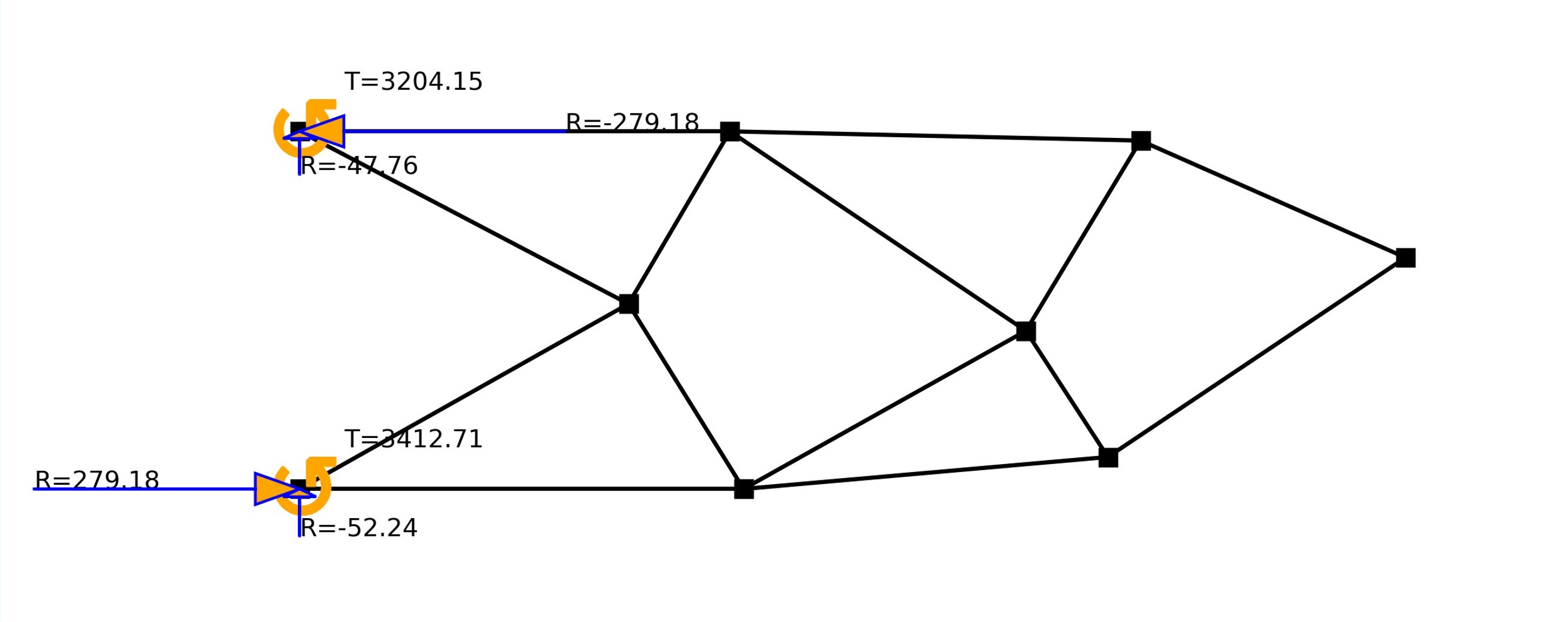}
		\label{fig:cant_strana_b}
	} \\[-3pt]
    \centering
	\subfloat[Axial Forces]{
		\includegraphics[width=0.49\textwidth]{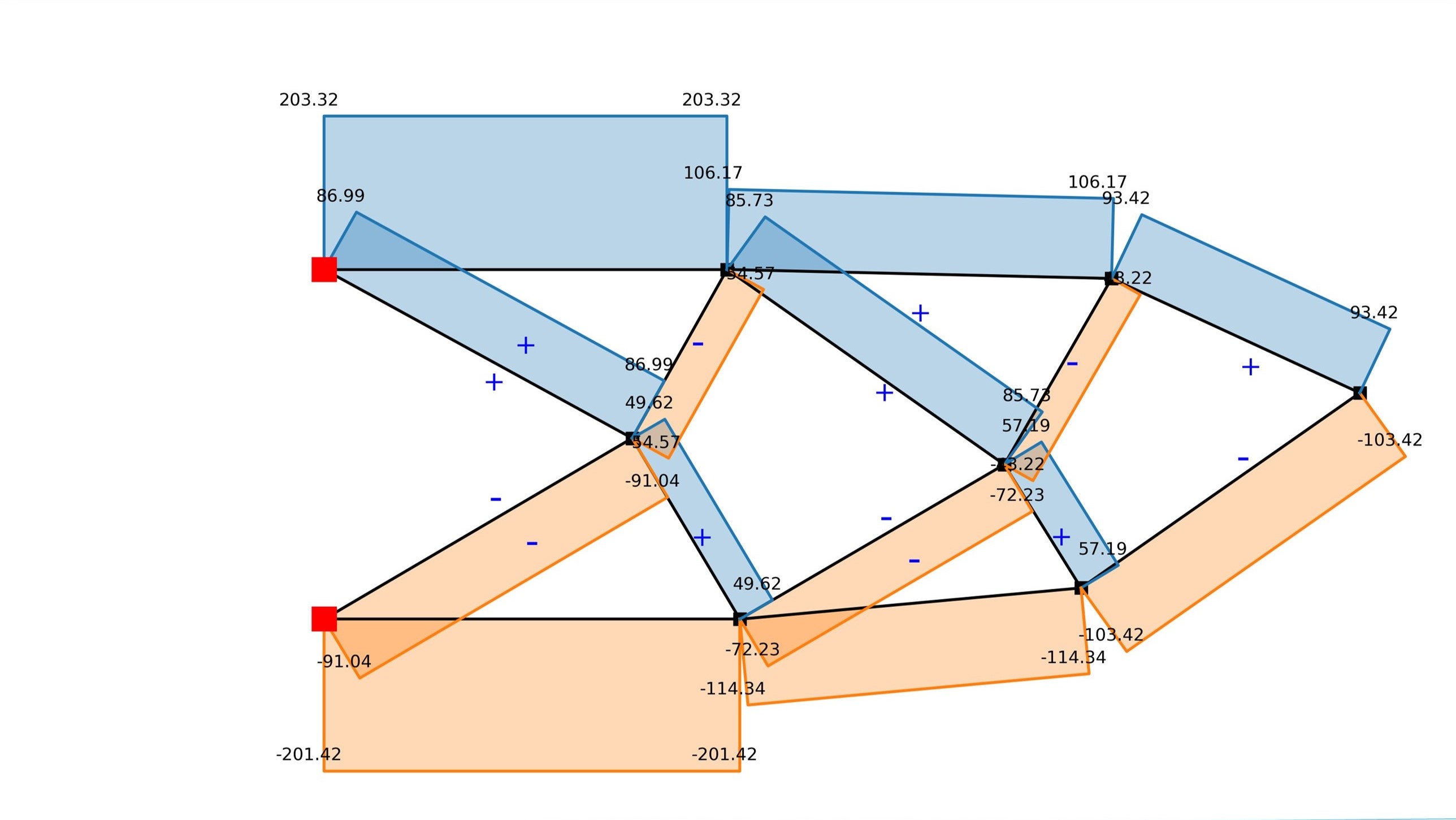}
		\label{fig:cant_strana_c}
	}
	\centering
	\subfloat[Bending Moments]{
		\includegraphics[width=0.49\textwidth]{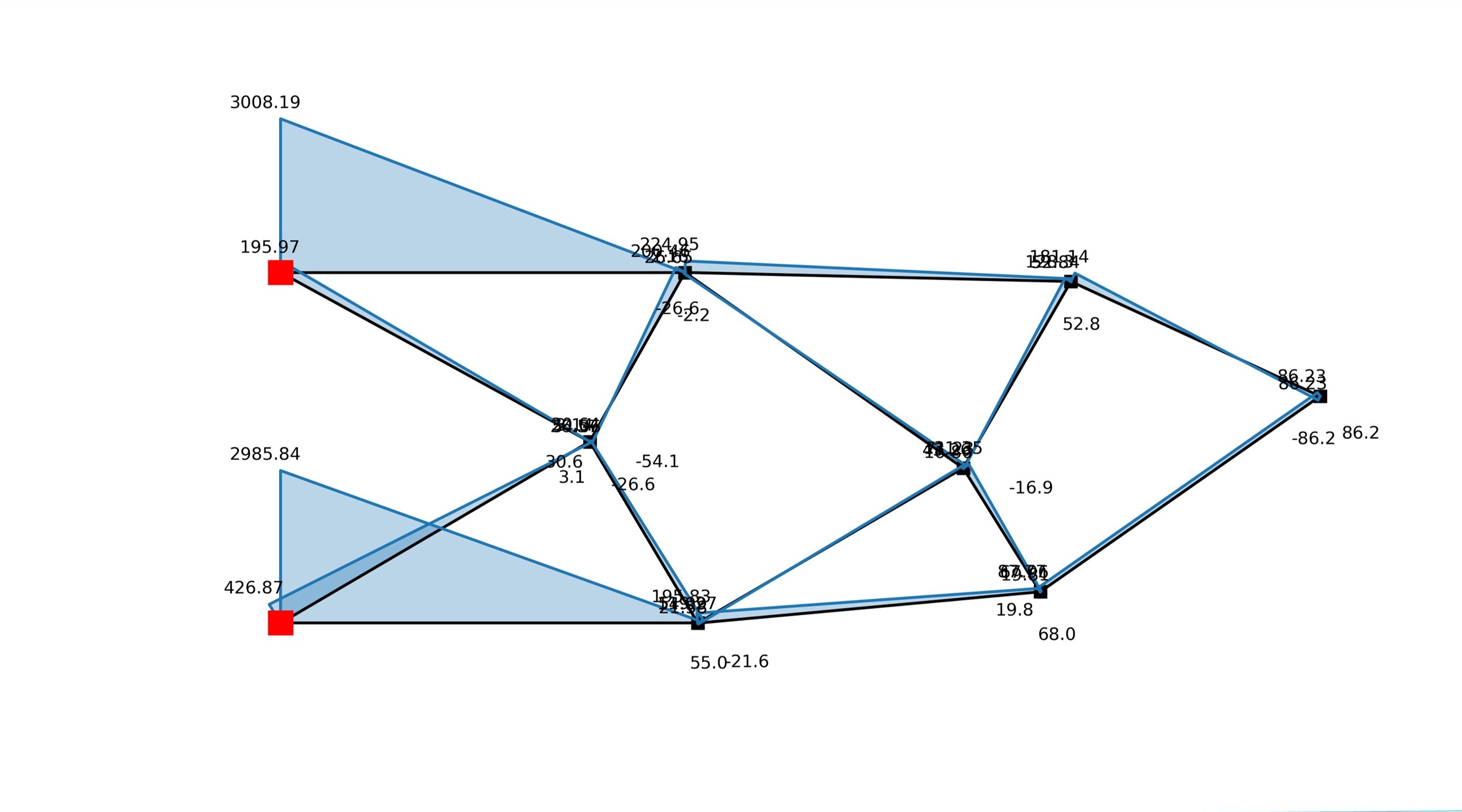}
		\label{fig:cant_strana_d}
	} \\[-3pt] 
    \centering
	\subfloat[Shear Forces]{
		\includegraphics[width=0.49\textwidth]{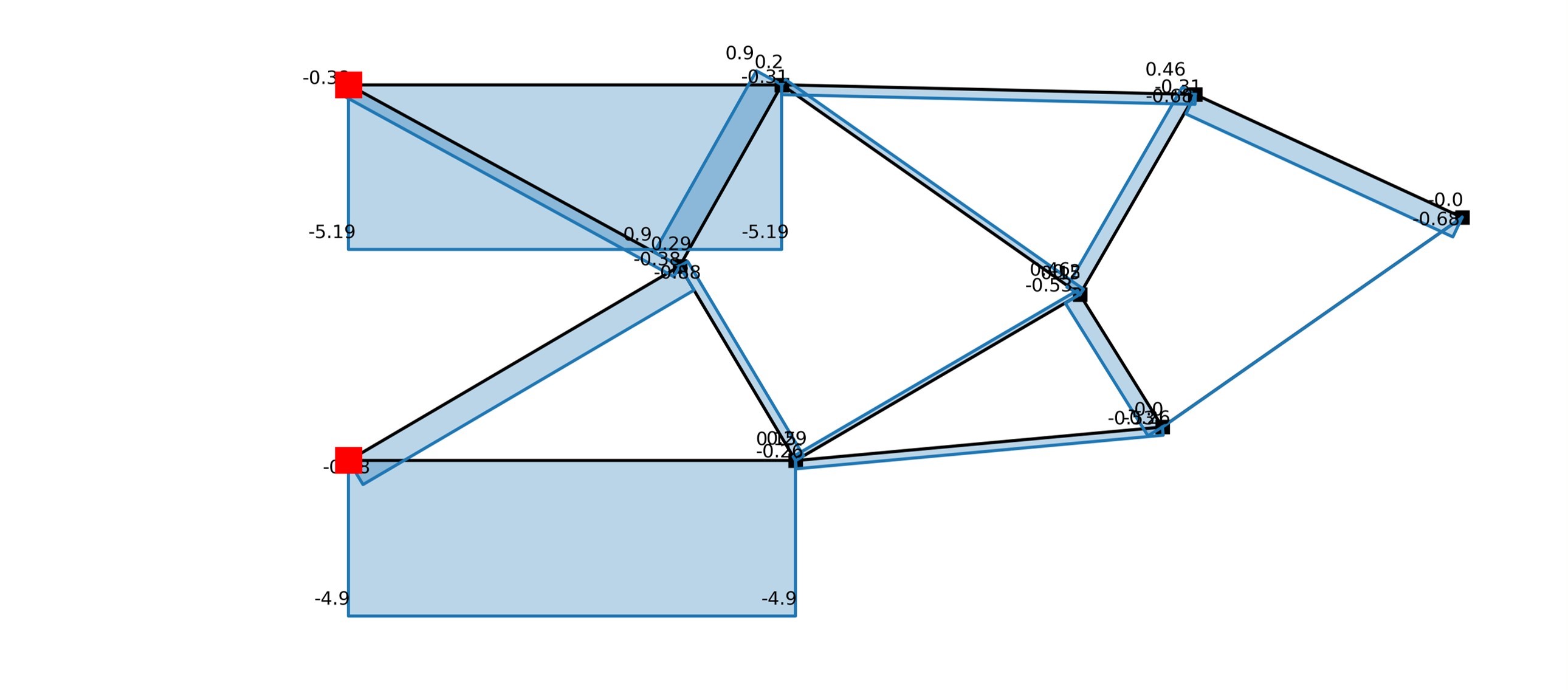}
		\label{fig:cant_strana_e}
	}
	\centering
	\subfloat[Deflection]{
		\includegraphics[width=0.49\textwidth]{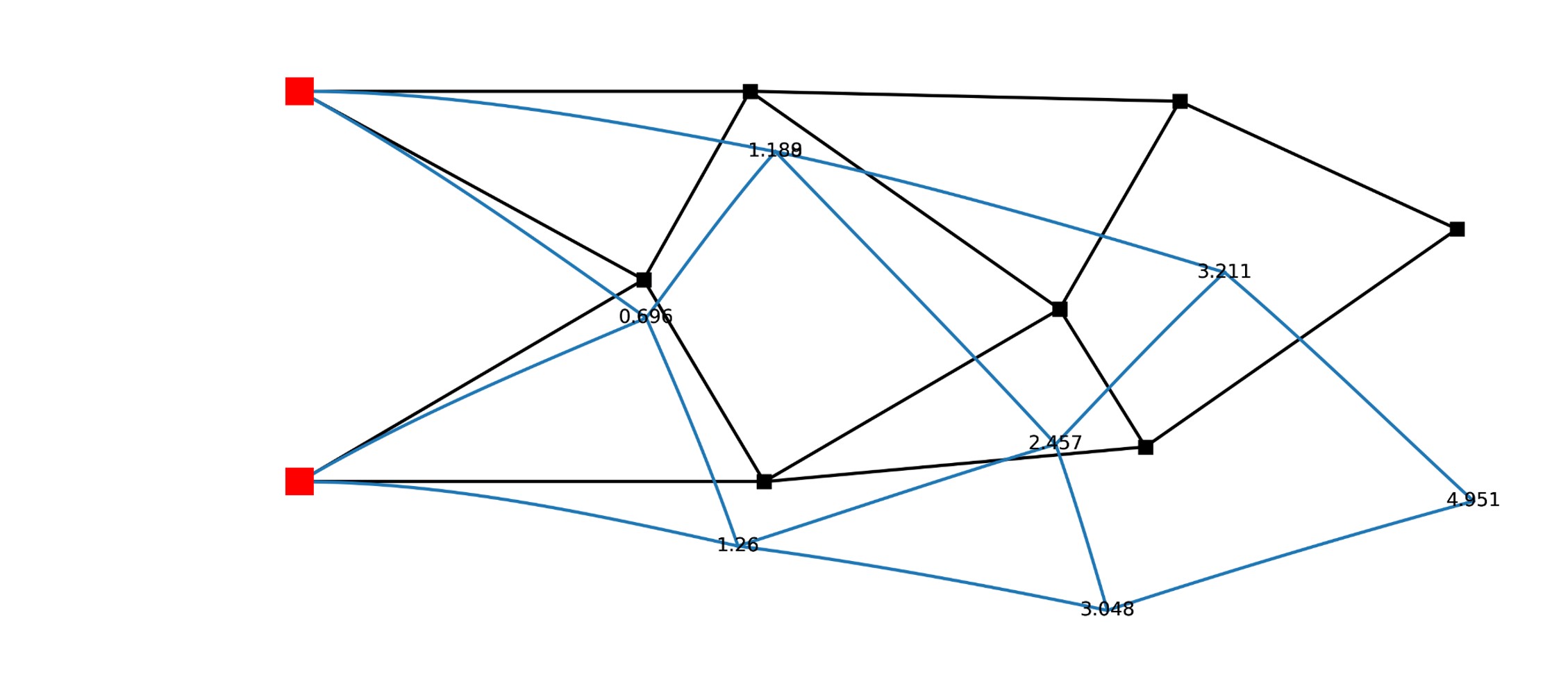}
		\label{fig:cant_strana_f}
	} \\
	\centering
	\caption{Structural analysis results of cantilever example (Forces are in kN and momements are in kNm)}
	\label{fig:cant_strana}
\end{figure}
According to the design standard \ref{subsec:excess_def} for serviceability conditions, the general deflection limit for the cantilever is determined as span/180. In this case, with a span of 1500 mm, the deflection limit is calculated as 8.33 mm. Taking into account an overall safety factor of 1.35, the corresponding serviceability load is obtained by dividing the applied load (P) by 1.35. For instance, if the applied load is 100 kN, the serviceability load is 74.08 kN. Under this serviceability load, the deflection of the cantilever is calculated to be 3.88 mm. Consequently, the deflection check for the serviceability condition is satisfied.

The validation of the model can be accomplished by assessing the objective function value and analyzing the structure using the determined optimum section properties. To determine the deflection corresponding to the final objective function value, the equation $F \times U = 0.566 \times 10^6$ is used. Given that $F$ is 100 kN, the calculated deflection value $U$ amounts to 5.66 mm. Furthermore, the Prokon software package yields a corresponding deflection of 5.72 mm for the model. By comparing the calculated deflection from the objective function value to the deflection obtained from the Prokon software package, a percentage error of 1$\%$ is observed. Therefore, the model demonstrates a close agreement with the Prokon software package, with a small percentage error of 1$\%$.

\begin{table}[H]
    \centering
    \caption{Optimized member cross-sectional details and alternative sections}
    \label{table:4}
    \begin{tabular}{cccccc}
    \hline
    \textbf{Member} &
       \textbf{\begin{tabular}[c]{@{}c@{}} $\mathbf{A_{opt}}$ \\ \SI[detect-weight]{}{\centi\meter^2} \end{tabular}} &
      \textbf{\begin{tabular}[c]{@{}c@{}} $\mathbf{I_{opt}}$ \\ \SI[detect-weight]{}{\centi\meter^4} \end{tabular}} &
      \textbf{Category} &
      \textbf{\begin{tabular}[c]{@{}c@{}}Alternate \\ CHS \\ Section\end{tabular}} &
      \textbf{\begin{tabular}[c]{@{}c@{}}Alternate \\ SHS \\ Section\end{tabular}} \\ \hline
    12 & 5.42  & 2.34  & Bottom & 48.3x4  & 60x6.3 \\ \hline
    13 & 6.12  & 2.98  & Bottom & 60.3x4  & 50x6.3 \\ \hline
    14 & 11.88 & 11.23 & Bottom & 101.6x4 & 50x6.3 \\ \hline
    1  & 11.63 & 10.76 & Top    & 101.6x4 & 60x6.3 \\ \hline
    3  & 6.65  & 3.52  & Top    & 60.3x4  & 50x6.3 \\ \hline
    5  & 6.01  & 2.87  & Top    & 60.3x4  & 50x6.3 \\ \hline
    2  & 5.87  & 2.74  & Web    & 60.3x4  & 50x6.3 \\ \hline
    4  & 3.32  & 0.88  & Web    & 42.4x4  & 50x6.3 \\ \hline
    6  & 4.04  & 1.3   & Web    & 42.4x4  & 50x6.3 \\ \hline
    7  & 3.03  & 0.73  & Web    & 42.4x4  & 50x6.3 \\ \hline
    8  & 4.93  & 1.93  & Web    & 48.3x4  & 50x6.3 \\ \hline
    9  & 2.42  & 0.47  & Web    & 42.4x4  & 50x6.3 \\ \hline
    10  & 3.47  & 0.96  & Web    & 42.4x4  & 50x6.3 \\ \hline
    11 & 5.55  & 2.45  & Web    & 48.3x4  & 50x6.3 \\ \hline
    \end{tabular}
\end{table}

\begin{table}[H]
    \caption{Utilization ratios of structural members of the optimal frame}
    \label{table:5}
    \centering
    \resizebox{\columnwidth}{!}{%
    \begin{tabular}{ccccccc}
    \hline
    \textbf{Member} &
      \textbf{Axial} &
      \textbf{Buckling} &
      \textbf{Bending} &
      \textbf{Shear} &
      \textbf{\begin{tabular}[c]{@{}c@{}}Axial + \\ Bending\end{tabular}} &
      \textbf{Remarks} \\ \hline
    1  & 0.6 (C)  & 0.59 & 0.28 & 0.04    & 0.79 & OK \\ \hline
    2  & 0.47 (C) & 0.48 & 0.12 & 0.01    & 0.57 & OK \\ \hline
    3  & 0.59 (C) & 0.6  & 0.06 & 0.004   & 0.63 & OK \\ \hline
    4  & 0.43     &      & 0.08 & 0.01    & 0.53 & OK \\ \hline
    5  & 0.53 (C) & 0.54 & 0.03 & 0.00004 & 0.73 & OK \\ \hline
    6  & 0.54 (C) & 0.58 & 0.03 & 0.003   & 0.59 & OK \\ \hline
    7  & 0.37     &      & 0.07 & 0.01    & 0.47 & OK \\ \hline
    8  & 0.33 (C) & 0.33 & 0.03 & 0.004   & 0.61 & OK \\ \hline
    9  & 0.33 (C) & 0.33 & 0.08 & 0.01    & 0.41 & OK \\ \hline
    10 & 0.41 (C) & 0.41 & 0.12 & 0.02    & 0.50 & OK \\ \hline
    11 & 0.57     &      & 0.09 & 0.01    & 0.71 & OK \\ \hline
    12 & 0.61     &      & 0.08 & 0.01    & 0.59 & OK \\ \hline
    13 & 0.55     &      & 0.06 & 0.004   & 0.63 & OK \\ \hline
    14 & 0.6      &      & 0.29 & 0.04    & 0.96 & OK \\ \hline
    \end{tabular}%
    }
\end{table}

\begin{figure}[H] 
	\centering
	\subfloat[]{
		\includegraphics[width=0.49\textwidth]{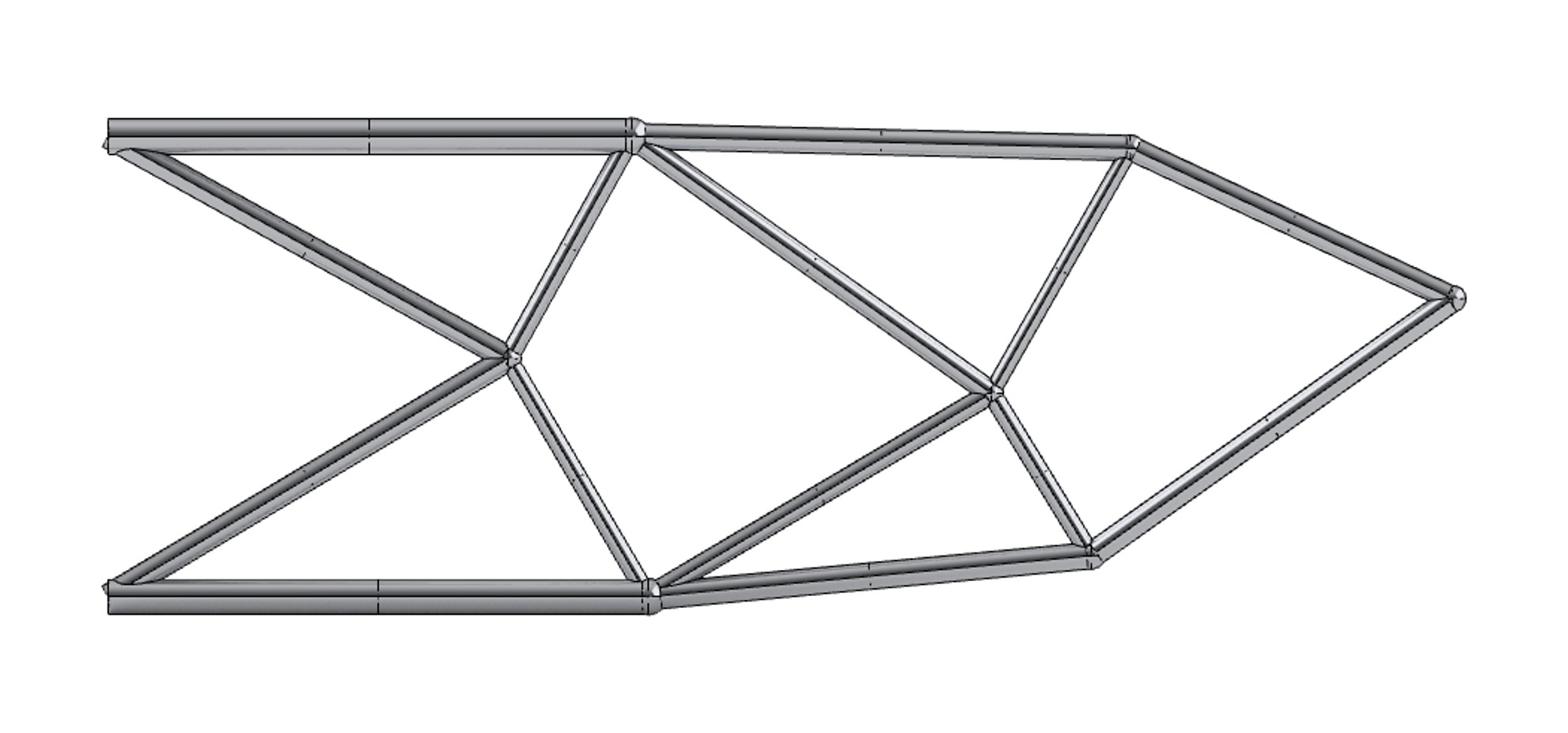}
		\label{fig:cant_ren_a}
	}
	\centering
	\subfloat[]{
		\includegraphics[width=0.49\textwidth]{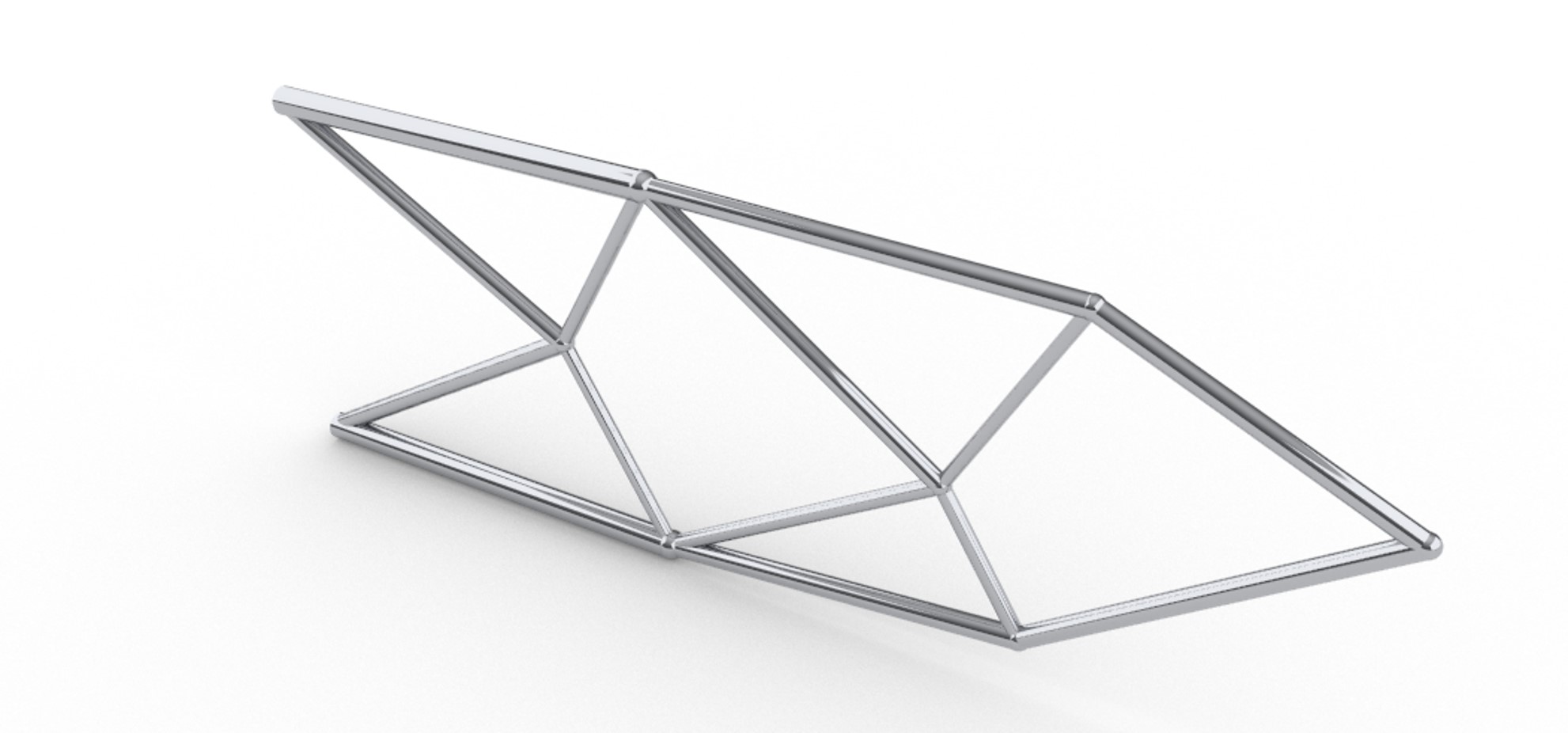}
		\label{fig:cant_ren_b} 
    }
	\centering
	\caption{Rendered views of the Optimized Cantilever}
	\label{fig:cant_ren}
\end{figure}

\subsection{A Simply Supported Beam}\label{sec:simp_beam}

\begin{figure}[H]
\centering
\includegraphics[width=0.99\textwidth]{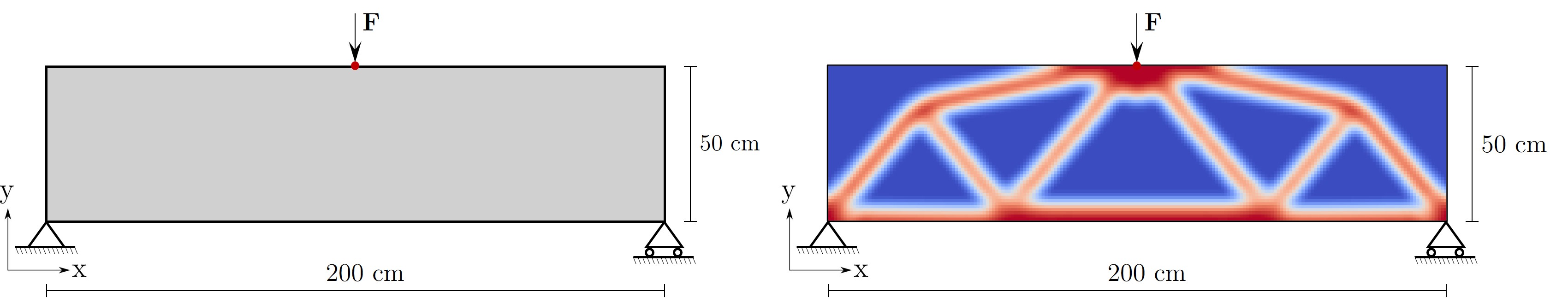}
\caption{Topology optimization of a simply supported beam example. ($\overline{E}$ = $2.1 \times 10^5$ N/mm$^2$, ${E}_{min}$ = $1 \times 10^{-9}$ N/mm$^2$, $\nu$ = $0.3$, $V_f$ = $0.3$, $p$ = $3$, $R$ = $5.4$, Max.iter = $200$)}
\label{fig:topopt_simp}
\end{figure}

\begin{figure}[H]
\centering
\includegraphics[width=0.99\textwidth]{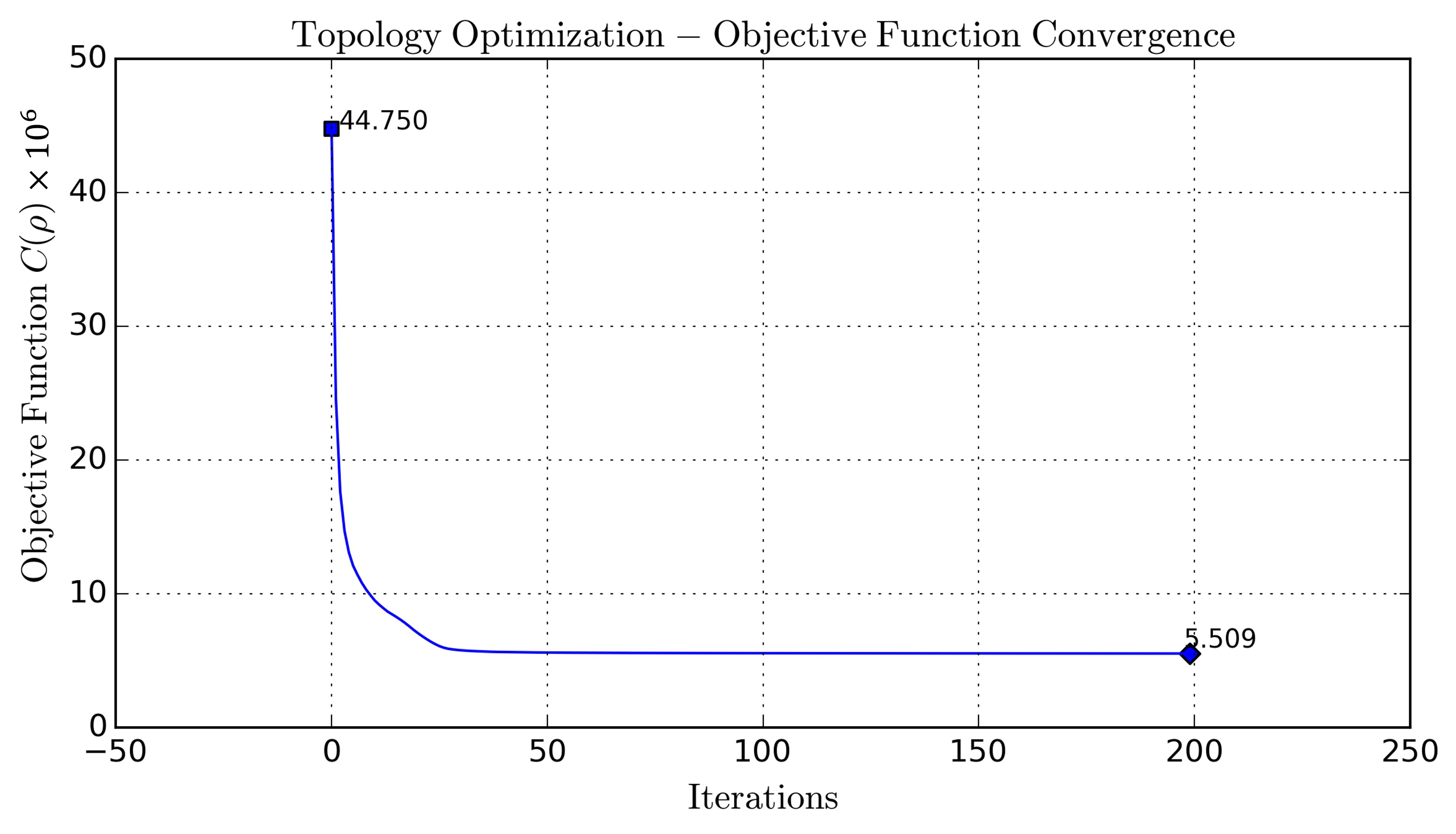}
\caption{Objective function convergence during topology optimization}
\label{fig:simp_obj_con}
\end{figure}

A simply supported beam (200 cm $\mathrm{\times}$ 50 cm) as shown in Figure \ref{fig:topopt_simp} is selected and a point load with a value of F = 100 kN is applied at the midpoint. The thickness of the plate is taken as 1 cm for size optimization.  In topology optimization, the domain was discretized into 200 $\mathrm{\times}$ 50 $\mathrm{\times}$ 1 linear hexahedral element and the objective function $C_{(\rho )}\ $ is minimized.  In this example, the optimization process incorporated specific parameters, including a volume fraction ($V_f$) of 0.3, a penalization power ($p$) of 3, a filter radius ($R$) of 5.4, and a maximum of 200 iterations. Figure \ref{fig:simp_obj_con} shows the convergence of the objective function. From the results, the 2D surface results were obtained for the skeletonization.

Similar to the previous cantilever example, the 2D binary image data obtained from the topology optimization is skeletonized and the frame model is extracted. Simply supported boundary pixels and force-applied pixels were tagged to as not to alter them during any optimization process. 

Sequential size and layout optimization are performed on the extracted frame model. In this example, the extracted frame model is directly used for size and layout optimization without merging any nodes and set to remove during layout optimization. For size optimization, the lower of the area is taken as circular areas with a lower bound radius of 0.5 cm and upper bound radius of 10 cm. For the first size optimization iteration, the initial sizes of the members are assumed to have the same cross-sectional area of $A=5.04\ cm^2$ which is found by diving the total volume $V_{\left(A\right)}=0.5\overline{V}=3000\ cm^3$ by the total length of the members. A merge ratio of 0.1 is used to remove short members after each layout optimization. 

\begin{figure}[h!] 
	\centering
        \captionsetup{justification=centering}
	\subfloat[Binary image generated for corresponding topology optimized model in Figure \ref{fig:topopt_simp}]{
		\includegraphics[width=0.49\textwidth]{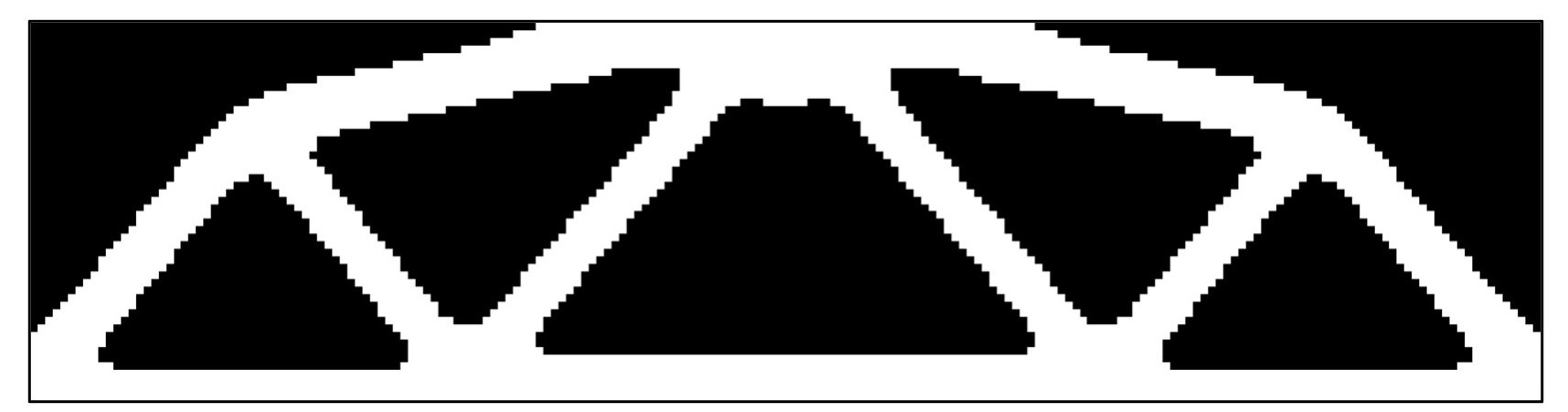}
		\label{fig:simp_skelframe_a}
	}
	\centering
	\subfloat[Skeleton obtained after tagging boundary pixels and loaded element pixel]{
		\includegraphics[width=0.49\textwidth]{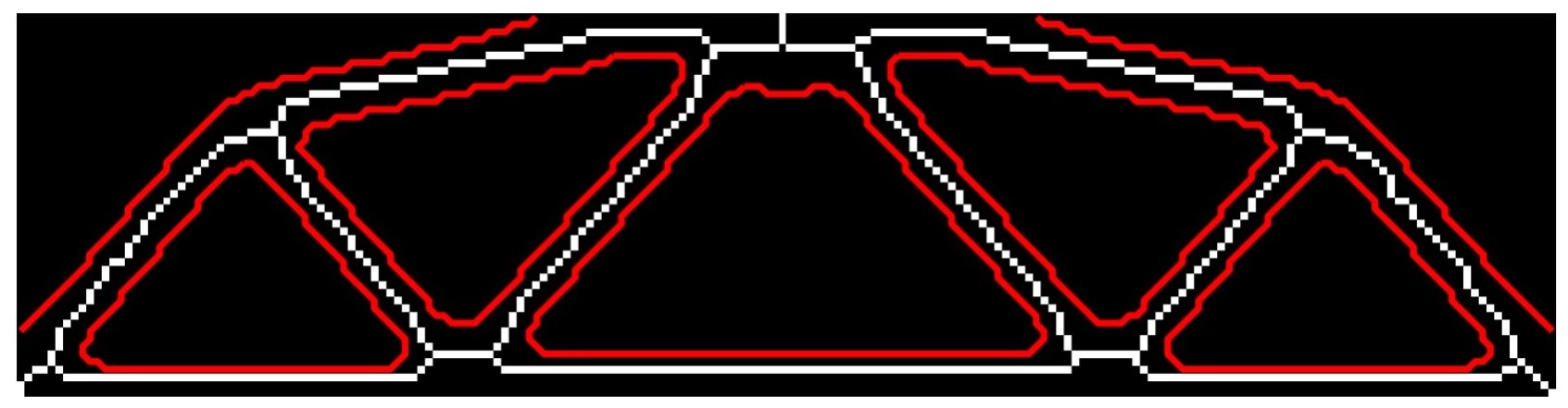}
		\label{fig:simp_skelframe_b}
	} \\
        \centering
	\subfloat[Graph building using identified featured pixels]{
		\includegraphics[width=0.49\textwidth]{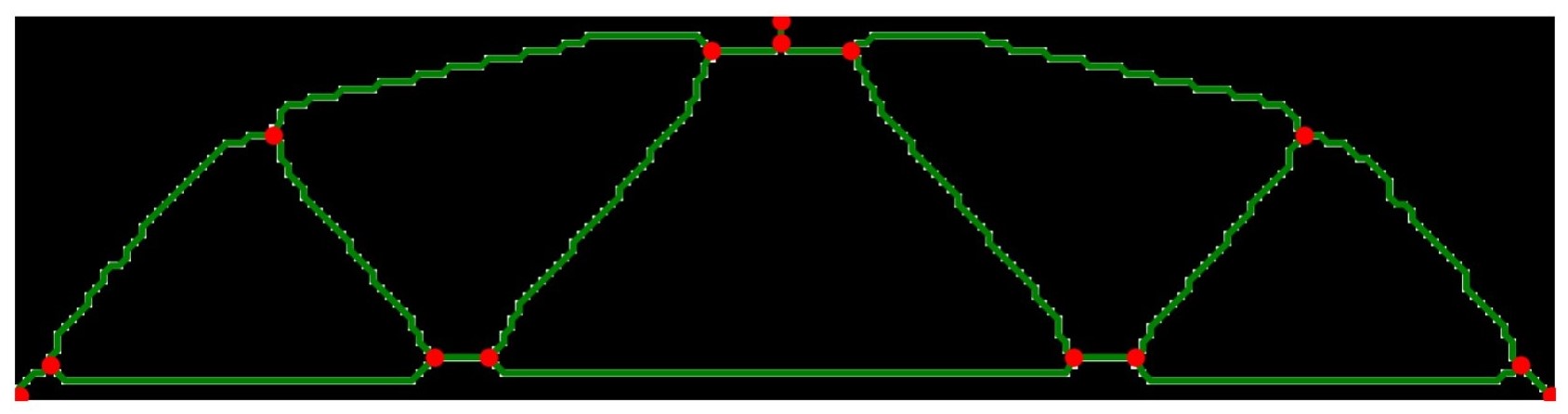}
		\label{fig:simp_skelframe_c}
	}
	\centering
	\subfloat[Initial graph model obtained from skeleton]{
		\includegraphics[width=0.49\textwidth]{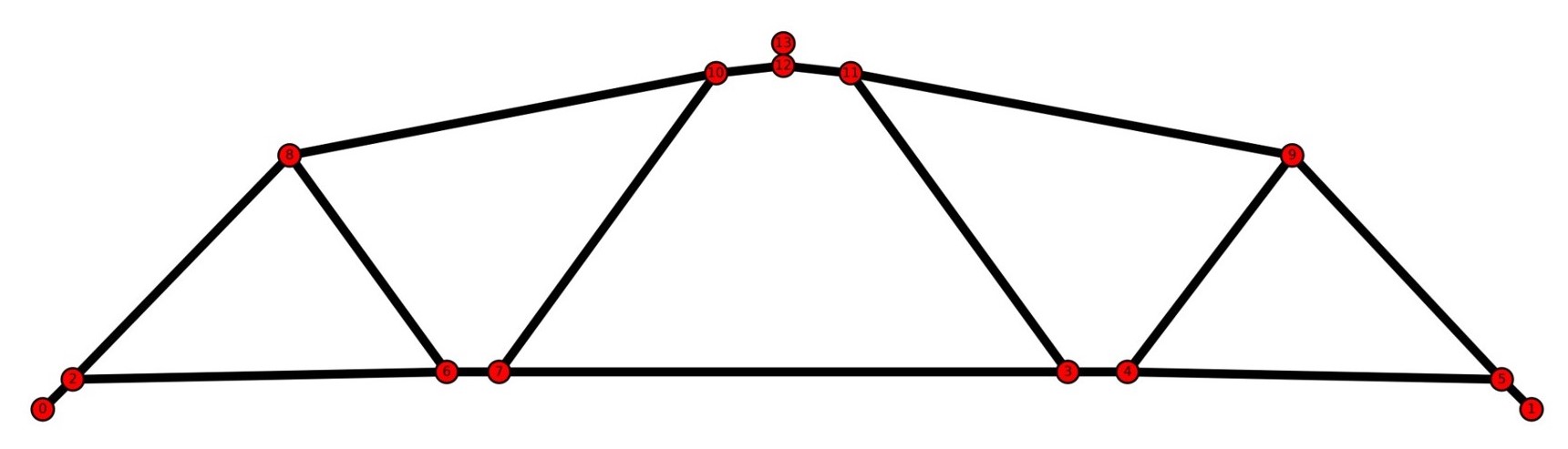}
		\label{fig:simp_skelframe_d}
	} \\ 
	\centering
	\caption{Skeletonization and frame model generation for simply supported example}
	\label{fig:simp_skelframe}
\end{figure}

\begin{figure}[H]
\centering
\includegraphics[width=0.99\textwidth]{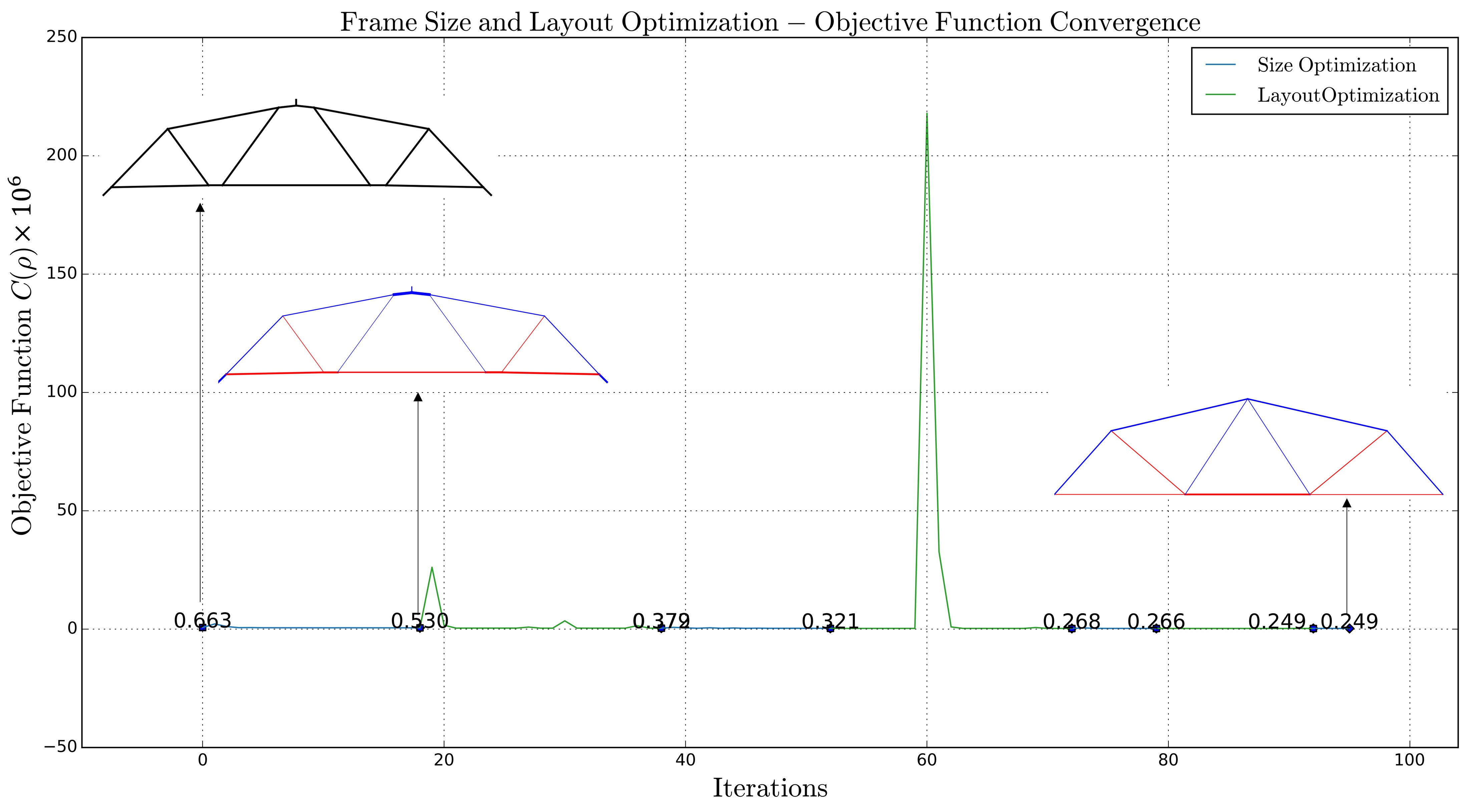}
\caption{Objective Function Convergence for Frame Size and Layout Optimization using MMA Algorithm when merging is allowed}
\label{fig:simp_mma_merge} 
\end{figure}

\begin{table}[h!]
    \centering
    \caption{Parameters used for Size and Layout Optimization}
    \label{table:7}
    \begin{tabular}{p{3.5in} c} \hline 
    \centering \textbf{Parameter} & \textbf{Value} \\ \hline 
    Volume constraint ($V_{(A,s)})$ & $\le 3000\ ${cm}$^3$ \\ \hline 
    Tolerance (${\epsilon }_{size}\ ,{\epsilon }_{layout},{\epsilon }_{frame}$) & ${10}^{-4}$ \\ \hline 
    Minimum Area ($A_{min}$) & $0.785\ ${cm}$^2$ \\ \hline 
    Maximum Area ($A_{max}$) & $\mathrm{314.160}\ ${cm}$^2$ \\ \hline 
    Bounds for Nodes & Whole design domain \\ \hline 
    Maximum Iterations for MMA Algorithm for each Size and Layout Optimization ($i_{max}$) & 20 \\ \hline 
    Merge Ratio $(\zeta )$ & 0.1 \\ \hline 
    \end{tabular}
\end{table}
Throughout the size and layout optimization, the volume of the frame structure is constrained to be equal to the volume defined by $V_f\overline{V}=3000\ cm^3$. The MMA algorithm is utilized for optimization and results are presented in Figure \ref{fig:simp_mma_merge}. The parameters used for both size and layout optimization are given in Table~\ref{table:7}.

Finally, a CAD model of the structure is generated for manufacturing with obtained node, member and cross-section details of the final optimized structure. The possibility of using Circular Hollow Sections (CHS) of thickness 4 mm as alternatives was analyzed and more suitable sections were suggested based on the final size optimization results as given in Table~\ref{table:8}.

\begin{figure}[H] 
	\centering
    \captionsetup{justification=centering}
	\subfloat[Rendered view]{
		\includegraphics[width=0.49\textwidth]{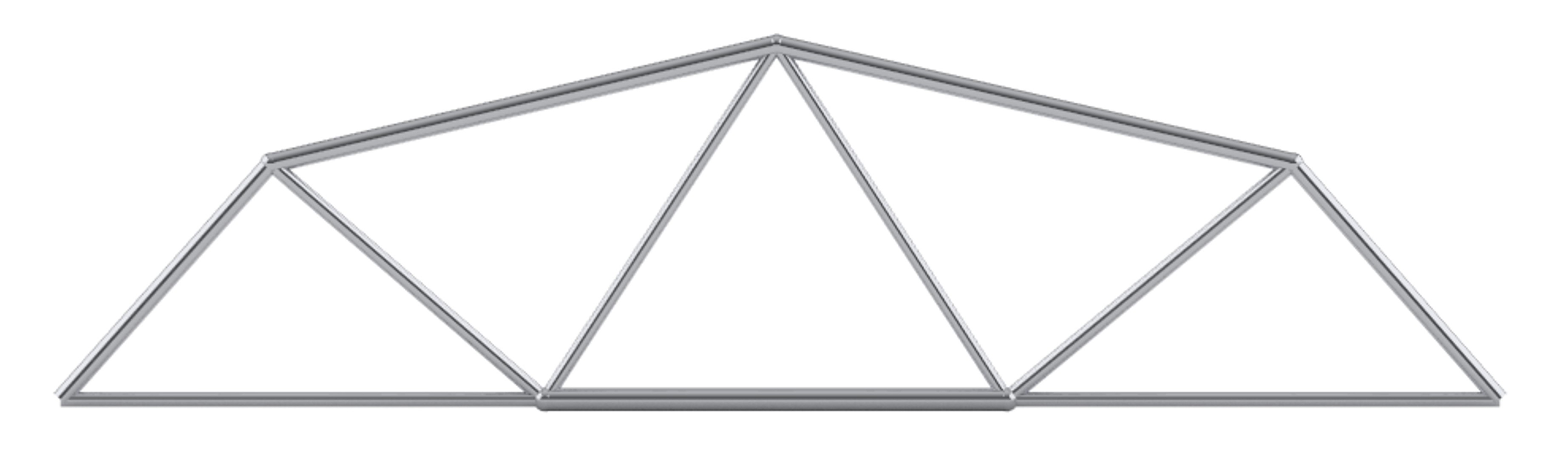}
		\label{fig:simp_rend}
	}
	\centering
	\subfloat[Axial Stresses (Tension - Red, Compression - Blue)]{
		\includegraphics[width=0.49\textwidth]{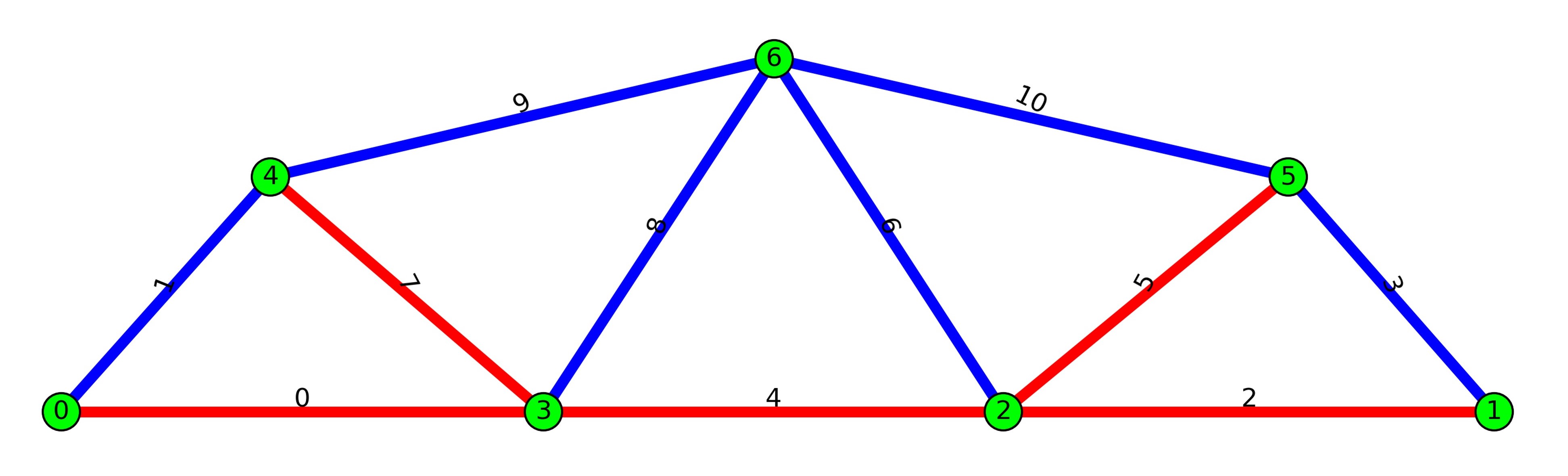}
		\label{fig:simp_axial_stress} 
    }
	\centering
	\caption{Rendered view and Axial Stresses of Members of the Simply Supported Beam}
	\label{simp_rend_and_axialstressess}
\end{figure}

\begin{table}[H]
    \centering
    \caption{Optimized member cross-sectional details and alternative sections}
    \label{table:8}
    \begin{tabular}{cccccc}
    \hline
    \textbf{Member} &
      \textbf{\textbf{${A_{opt}}$ (cm\textsuperscript{2})}} &
      \textbf{\textbf{${I_{opt}}$ (cm\textsuperscript{4})}} &
      \textbf{Category} &
      \textbf{\begin{tabular}[c]{@{}c@{}}Alternate \\ CHS \\ Section\end{tabular}} \\ \hline
    1 & 3.39 & 0.91 & Bottom & 42.4x4 \\ \hline 
    3 & 3.30 & 0.87 & Bottom & 42.4x4 \\ \hline 
    5 & 7.77 & 4.8 & Bottom & 76.1x4 \\ \hline 
    2 & 5.07 & 2.04 & Top & 48.3x4 \\ \hline 
    4 & 4.99 & 1.98 & Top & 48.3x4 \\ \hline 
    10 & 6.28 & 3.14 & Top & 60.3x4 \\ \hline 
    11 & 6.26 & 3.12 & Top & 60.3x4 \\ \hline 
    6 & 3.72 & 1.1 & Web & 42.4x4 \\ \hline 
    7 & 2.81 & 0.63 & Web & 42.4x4 \\ \hline 
    8 & 3.64 & 1.05 & Web & 42.4x4 \\ \hline 
    9 & 2.85 & 0.65 & Web & 42.4x4 \\ \hline 
    \end{tabular}
\end{table}

\begin{table}[H]
    \caption{Utilization ratios of structural members of the optimal frame}
    \label{table:9}
    \centering
    \resizebox{\columnwidth}{!}{%
    \begin{tabular}{ccccccc}
    \hline
    \textbf{Member} & \textbf{Axial} & \textbf{Buckling} & \textbf{Bending} & \textbf{Shear} & \textbf{\begin{tabular}[c]{@{}c@{}}Axial + \\ Bending\end{tabular}} & \textbf{Remarks} \\ \hline
    1         & 0.39      &           & 0.12      & 0.01      &     0.50      & OK \\ \hline
    2         & 0.51 (C)  & 0.52      & 0.03      & 0.001     &     0.46      & OK \\ \hline
    3         & 0.39      &           & 0.12      & 0.008     &     0.48     & OK \\ \hline
    4         & 0.50 (C)  & 0.51      & 0.03      & 0.002     &     0.46      & OK \\ \hline
    5         & 0.47      &           & 0.06      & 0.0001    &     0.46      & OK \\ \hline
    6         & 0.42      &           & 0.06      & 0.003     &     0.45      & OK \\ \hline
    7         & 0.32      &           & 0.06      & 0.003     &     0.32      & OK \\ \hline
    8         & 0.41      &           & 0.06      & 0.003     &     0.42      & OK \\ \hline
    9         & 0.32 (C)  & 0.35      & 0.05      & 0.004     &     0.33      & OK \\ \hline
    10        & 0.50 (C)  & 0.52      & 0.12      & 0.008     &     0.53     & OK \\ \hline
    11        & 0.50 (C)  & 0.52      & 0.12      & 0.008     &     0.53      & OK \\ \hline
    \end{tabular}%
    }
\end{table}

\begin{figure}[H] 
	\centering
	\subfloat[Structure]{
		\includegraphics[width=0.49\textwidth]{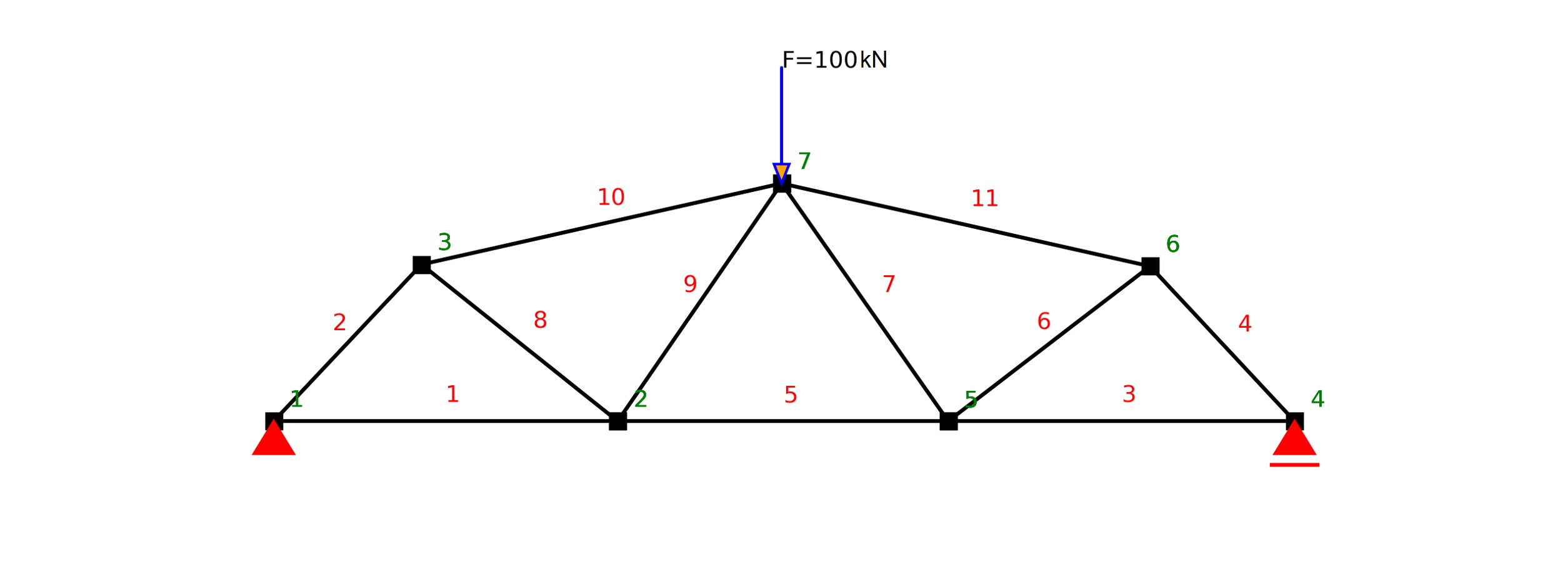}
		\label{fig:simp_strana_a}
	}
	\centering
	\subfloat[Reaction Forces]{
		\includegraphics[width=0.49\textwidth]{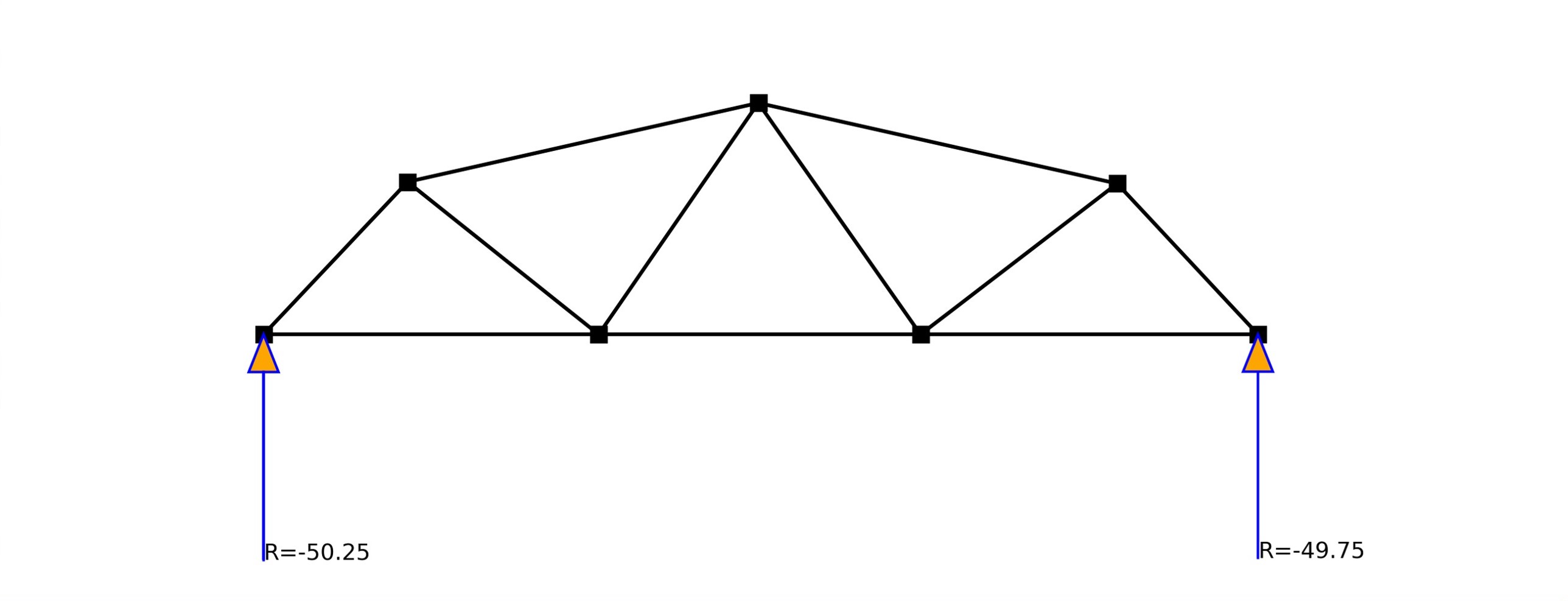}
		\label{fig:simp_strana_b}
	} \\[-3pt]
    \centering
	\subfloat[Axial Forces]{
		\includegraphics[width=0.49\textwidth]{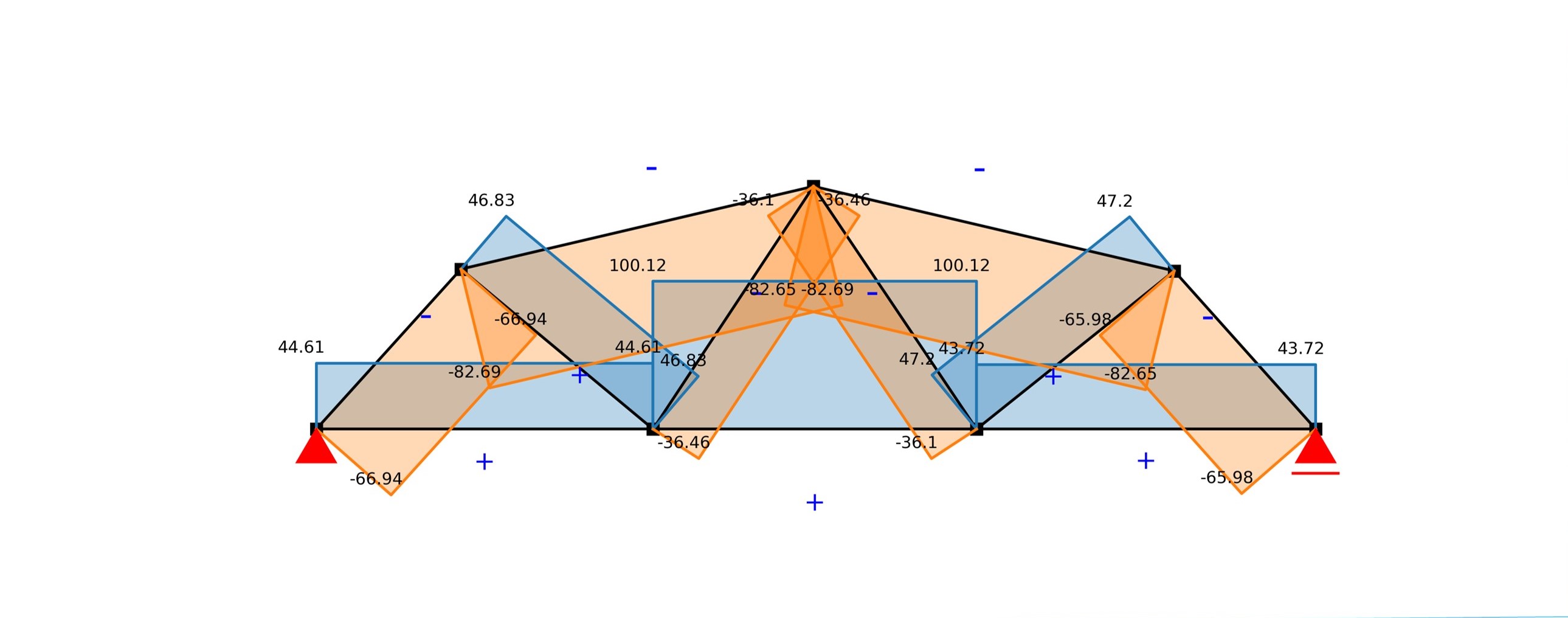}
		\label{fig:simp_strana_c}
	}
	\centering
	\subfloat[Bending Moments]{
		\includegraphics[width=0.49\textwidth]{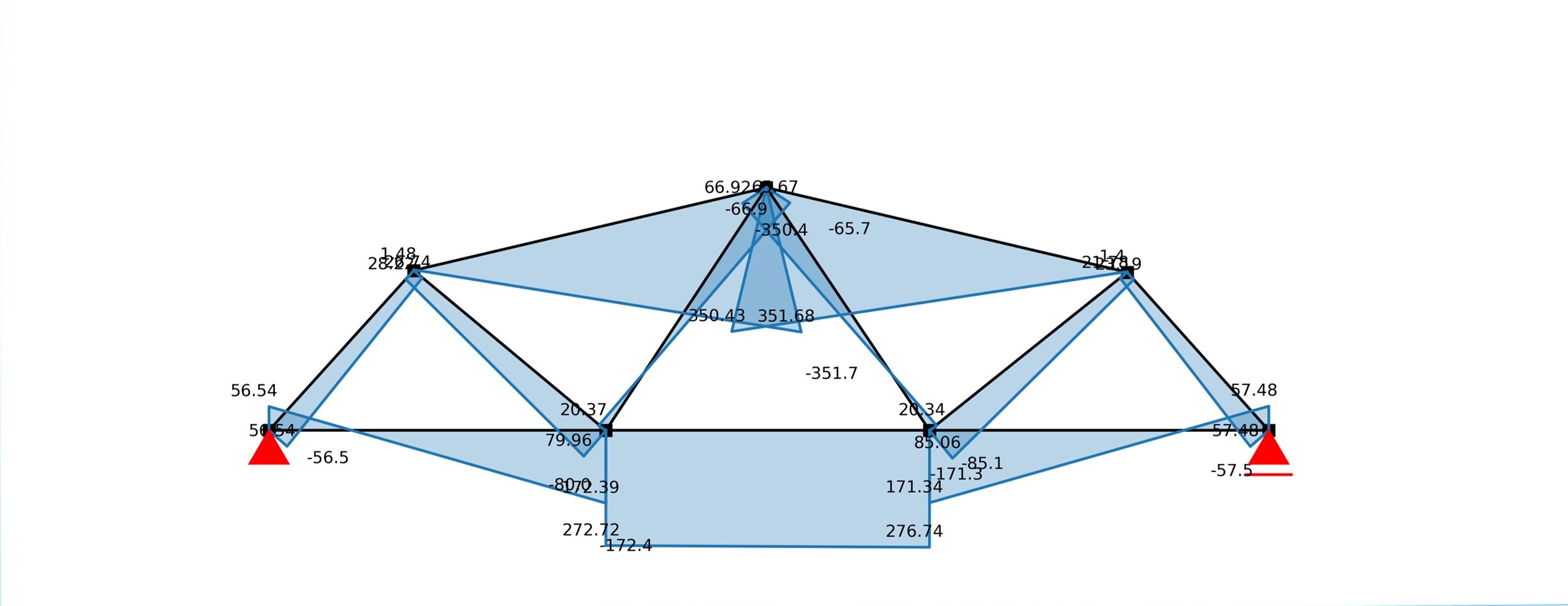}
		\label{fig:simp_strana_d}
	} \\[-3pt] 
    \centering
	\subfloat[Shear Forces]{
		\includegraphics[width=0.49\textwidth]{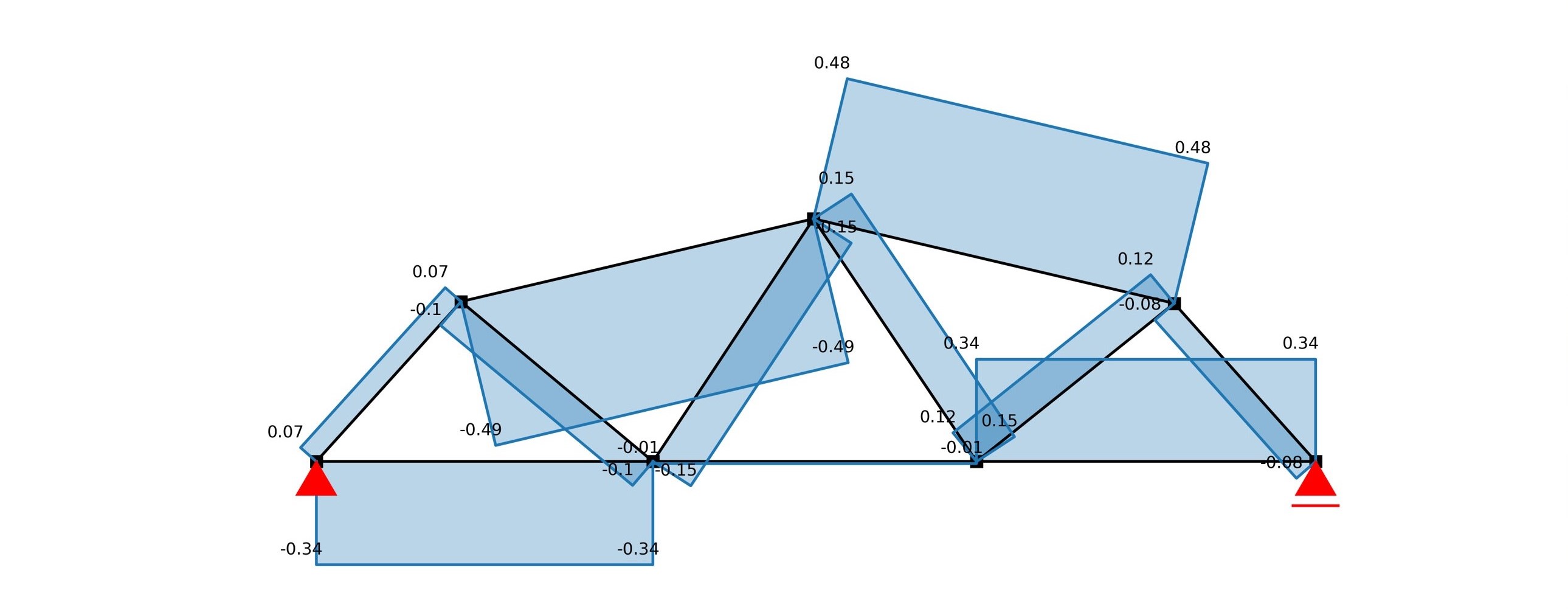}
		\label{fig:simp_strana_e}
	}
	\centering
	\subfloat[Deflection]{
		\includegraphics[width=0.49\textwidth]{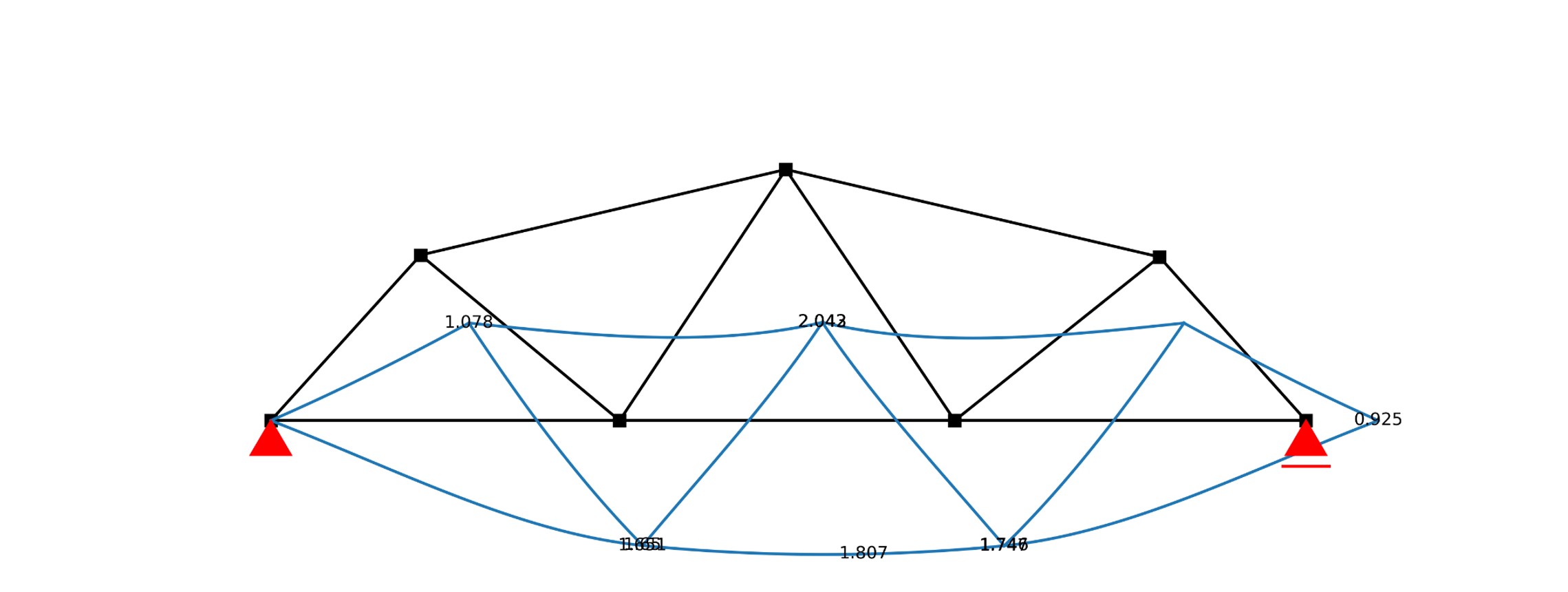}
		\label{fig:simp_strana_f}
	} \\ 
	\centering
	\caption{Structural analysis results of simply supported example (Forces are in kN and momements are in kNm)}
	\label{fig:simp_strana}
\end{figure}
For the selected sections, the structural analysis and design was done similar to the previous example. The analysis results and the utilization ratios according to Eurocode 3 is presented below in Figure \ref{fig:simp_strana} and Table~\ref{table:9}.

\subsection{Similar Cantilever and Simply Supported Examples}\label{sec:sim_examples}
\begin{figure}[H]
\centering
\includegraphics[width=0.90\textwidth]{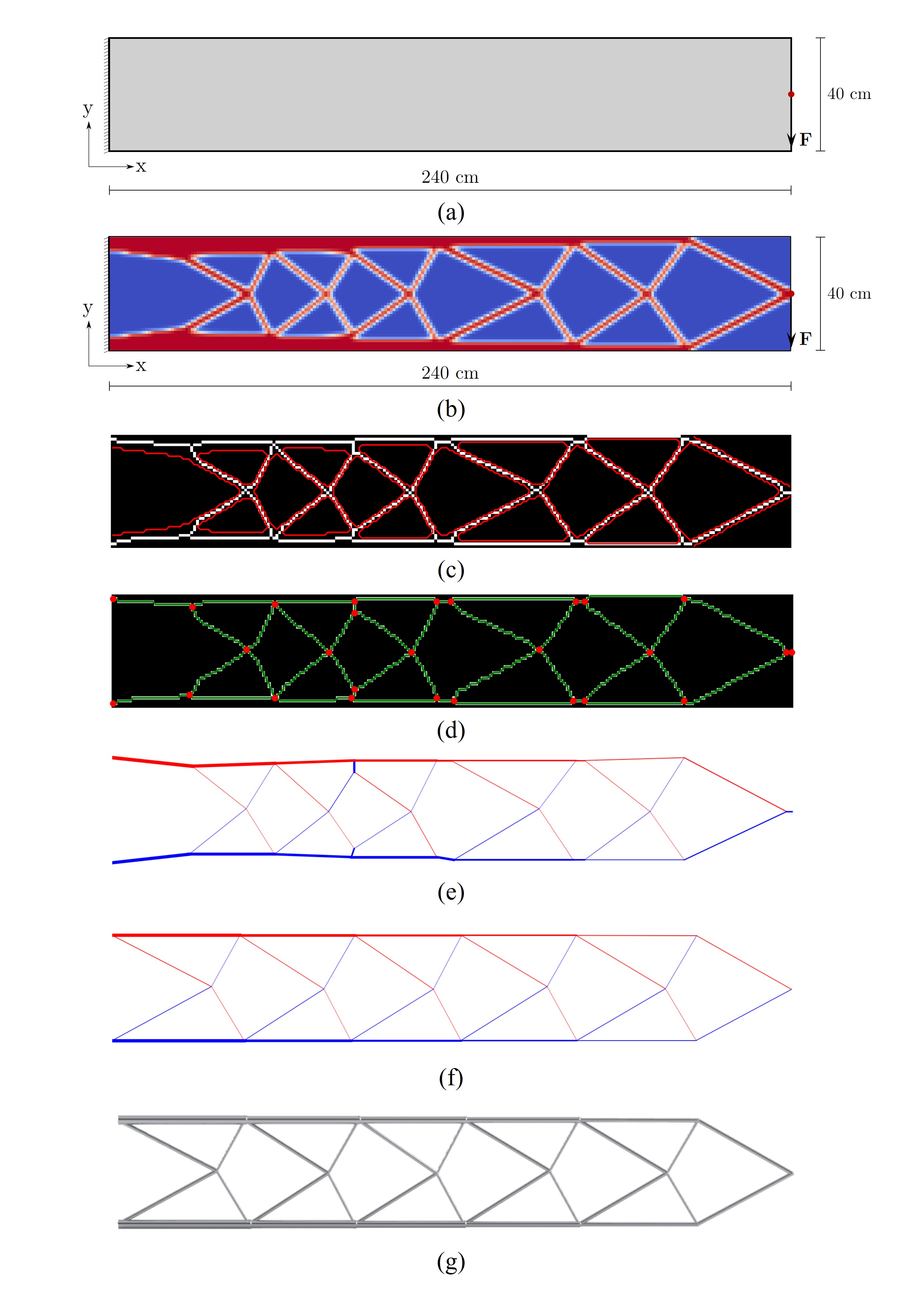}
\caption{Cantilever Example II. (a) Problem, (b) Topology optimized model for volume fraction $V_{f}=0.3$, (c) Skeletonization, (d) Graph model generation, (e) Initial frame model, (f) Final optimized frame model, (g) Generated CAD model}
\label{fig:cant_ex_2}
\end{figure}

\begin{figure}[H]
\centering
\includegraphics[width=0.90\textwidth]{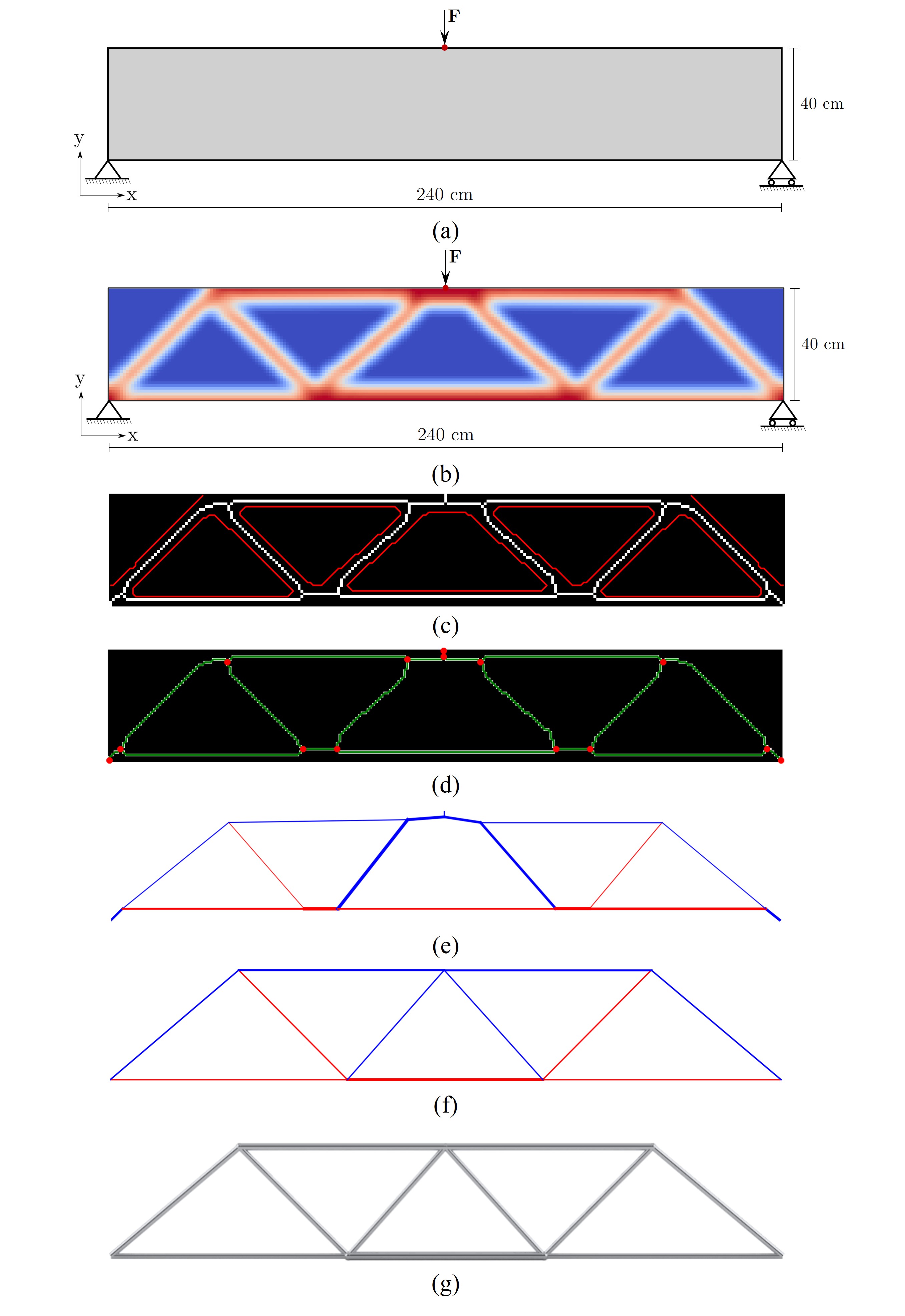}
\caption{Simply Supported Beam Example II. (a) Problem, (b) Topology optimized model for volume fraction $V_{f}=0.3$, (c) Skeletonization, (d) Graph model generation, (e) Initial frame model, (f) Final optimized frame model, (g) Generated CAD model}
\label{fig:simp_ex_2}
\end{figure}

Similar to the above examples in Sections \ref{sec:cant_plate} and \ref{sec:simp_beam} , various depth/span ratio of the design domain is optimized and visual outputs of them are presented here. Figure \ref{fig:cant_ex_2} shows the cantilever example and Figure \ref{fig:simp_ex_2} shows the simply supported beam example. For both these examples, the same material properties are utilized as used in above examples. For topology optimization of the cantilever example II, the filter radius is set to 1.2. Other than that, for both examples optimization parameters in Table \ref{table:5} are utilized. After the initial frame model extraction, the MMA algorithm is directly used to optimize the frame model.
%


\section{Conclusions}
\label{sec:conclusion}

In conclusion, this research addressed the lack of unified frameworks for generating ready-to-manufacture parametric CAD models of structurally viable designs. A novel approach was presented, which involved topology optimization, skeletonization, frame extraction, size and layout optimization, and CAD model generation with structural integrity assessment. The workflow was implemented in a single parametric platform Rhino-Grasshopper. The presented approach demonstrated the ability to generate structurally optimized and manufacturable CAD models efficiently. The structural performance of the final model can be assessed according to the standard codes of practice for verification. Overall, this research provides a significant contribution towards developing a practical workflow for generating structurally optimized CAD models.

This research study outcome can initially assist the structural designer for identify the optimal load path, geometric design and structural design and it can be further developed in various aspects such as: connection detailing, redundant members for alternative load paths, 3D structural problems and long term effects in steel (corrosion, relaxation, etc.)


\section*{Acknowledgment}
Financial assistance of the Senate Research Committee of University of Moratuwa (Research Grant No. SRC/LT/2020/01) is highly appreciated towards the completion of this research.


\section*{Data availability Statement}
Open-source implementations will be made available upon the publication of the paper.


\bibliographystyle{unsrt}
\bibliography{design-informed}

\end{document}